\documentclass[10pt,twoside]{article}
\usepackage{amsmath,amsfonts,amssymb,amsthm,amscd}
\usepackage{color,cite}
\usepackage{doc}
\usepackage{graphicx}
\newtheorem{theorem}{Theorem}

\newtheorem{remark}{Remark}

\newcommand{\HYP}{\mathbb{H}^3}

\newcommand{\ba}{\mathbf{a}}

\newcommand{\cA}{\mathcal{A}}

\begin{document}

\title{Non-fundamental trunc-simplex tilings and their optimal hyperball packings and coverings in hyperbolic space \\ I. For families F1-F4}

\author{Emil Moln\'{a}r$^1$, Milica Stojanovi\'{c}$^{2}$, Jen\H{o} Szirmai$^1$ \\ \\
$^1$Budapest University of Tehnology and Economics Institute of Mathematics, \\Department of Geometry, H-1521 Budapest, Hungary\\ emolnar@math.bme.hu, \ szirmai@math.bme.hu \\
\\$^{2}$Faculty of Organizational Sciences, University of Belgrade\\11040 Belgrade, Serbia \\ milica.stojanovic@fon.bg.ac.rs\\}

\date{}

\maketitle

\pagestyle{myheadings} \markboth{\rm Moln\'{a}r, Stojanovi\'{c}, Szirmai}{\rm Non-fundamental trunc-simplex tilings ...}


\begin{abstract}
Supergroups of some hyperbolic space groups are classified as a continuation of our former works. Fundamental domains will be integer parts of 
truncated tetrahedra belonging to families F1 - F4, for a while, by the notation of E. Moln\'{a}r et al. in $2006$. 
As an application, optimal congruent hyperball packings and coverings to the truncation base planes with their very good densities are computed. 
This covering density is better than the conjecture of L.~Fejes~T\'oth for balls and horoballs in $1964$.
\end{abstract}






\section{Introduction} \label{sec:1}

\subsection{Short history}
Hyperbolic space groups are isometry groups, acting discontinuously on hyperbolic 
3-space $\mathbb{H}^3$ with compact fundamental domains. 
It will be investigated some series of such groups by looking for their fundamental domains. Face pairing 
identifications on a given polyhedron give us generators and relations 
for a space group by the algorithmically generalized Poincar\'e Theorem \cite{BM}, \cite{M92}, \cite{Mo11} (as in subsection 1.2).

The simplest fundamental domains are 3-simplices (tetrahedra) and their integer parts by inner symmetries. 
In the process of classifying the fundamental simplices, 64 combinatorially different face pairings of fundamental simplices were determined 
\cite{Z}, \cite{MP94}, and also 35 solid transitive non-fundamental simplex identifications \cite{MPS06}. I. K. Zhuk \cite{Z} 
classified Euclidean and hyperbolic fundamental simplices of finite volume. An algorithmic procedure was given 
by E. Moln\'{a}r and I. Prok \cite{MP88}. In \cite{MP94,MPS97,MPS06} the authors summarized all those results, arranging identified 
simplices into 32 families. Each of them is characterized by the so-called maximal series of simplex tilings, i.e. maximal symmetry groups with smallest 
fundamental domains. The other members are derived by symmetry breaking of this maximal one. These principles are also discussed by $D$-symbols in \cite{MPS06}, see also 1.2. 

Some complete cases of supergroups with fundamental truncated simplices (short\-ly trunc-simplices) are discussed in \cite{M90,S97,S10,S11, S13,S14,S17}. 
As a classification, supergroups of the groups with fundamental trunc-simplices are given in \cite{S19}. Especially, supergroups of the groups with 
fundamental simplices belonging to families F1-F4, by the notations given in \cite{MPS06}, are investigated in \cite{M90,S10,S14,S17}. 
For that purpose, the stabilizer group of the corresponding vertex figure is analyzed. The new method, for such analysis, established e.g. in the summary works
\cite{S17,S19}, is also confirmed by the results given in \cite{LMV18}.

\subsection{Preliminaries}
Generators and relations for a space group $G$ with a given polyhedron $P$ as a fundamental domain can be obtained by the Poincar\`{e} theorem, 
applying Poincar\`{e} algorithm. In this paper such polyhedra will be simplices and trunc-simplices, for overview see e.g. \cite{S17,S19}, and here mainly their integer parts. 

It is necessary to consider all {\it face pairing identifications} of such domains. These are isometries which generate a group $G$ and induce subdivision of 
faces and directed edge segments of $P$ into equivalence classes, such that an edge segment does not contain two $G$-equivalent points in its interior. 
The Poincar\`{e} algorithm (see \cite{BM}, \cite{M92}) gives us for each edge segment class one cycle transformation of the 
form $c = g_r \dots g_2 g_1$, where $g_i$, $i = 1,2, \ldots r$ are face pairing identifications. Since each of these cycle transformations 
is a rotation of order $\nu$, the cycle relations are of the form $\left(g_r \dots g_2 g_1\right)^{\nu} = 1$. The Poincar\`{e} 
theorem guarantees that these cycle relations, together with relations, $g_i^2 = 1$ to the occasional involutive 
generators $g_i = g_i^{-1}$, form a complete set of defining relations for $G$.

\begin{remark} As an intuitive picture (e.g. from \cite{M92}) we shall use in the following, the space tiling $P^G$ by fundamental 
domains (polyhedra) represents the space group $G$ itself after having chosen an arbitrary identity domain $P$. 
For any generator $g \in  G$ domain $P$ will be adjacent with the $g$-image domain $P^g$ along face $f_g$ denoted simply either 
by $g$ on the side (half space) of $P^g$, or by $g^{-1}$ on the side of $P$ (we prefer the first). Involutive generators 
$g = g^{-1}$ we write just on the corresponding face. Thus for any $\gamma \in G$, domain $P^{\gamma}$ will 
be adjacent with domain $P^{g \gamma }$ along an image face $f_g^{\gamma }$, signed simply either by $g$ on side 
of $P^{g \gamma }$, or by $g^{-1}$ on the side of $P^{\gamma }$. So domain $P^{\gamma }$ is obtained by 
$\gamma = g_r \dots g_2 g_1$ from $P$ through the path in domains $P$, $P^{g_1}$, $P^{g_2 g_1}$, $\ldots $ , $P^{\gamma }$  
entered finely through $g_r$. Imagine also the group graph by coloured directed edges with so many colours as many generators 
we have. From $g$ and $g^{-1}$ we use only $g$, as to direction of the corresponding graph edge.
\end{remark}

{\bf Classification all non-fundamental trunc-simplex tilings} means the following finite process (7 steps as a natural sketchy generalization of our earlier works 
by $D-$symbol method or barycentric subsimplex orbits, 
see e.g. \cite{MPS06} and \cite{Mo11}, you can skip this difficult part at the first read).
\begin{enumerate}
\item We consider a 3-simplex $T$ (as tetrahedron) with its barycentric 
subdivision into 24 subsimplices. Any subsimplex has a 3-centre (solid centre, first formally, later also metrically), a 
2-centre in a simplex face, a 1-centre 
in an edge of the provious face, a 0-centre as a vertex of the previous edge, as usual. The simplex has first the piece-wise linear (PL) topology. 
Any subsimplex has a $3-$face. $2-, 1-, 0-$faces, opposite to the corresponding centres. Then $2-,1-,0-$ (involutive (involutory)) adjacency operations $\sigma^i$ 
can be introduced. Then
occasional symmetry operations come acting on points, segments $\dots$, preserving the subsimplices and incidences with adjacency operations. E.g. we can write by convention
\begin{equation}
\sigma^iC^g:=(\sigma^iC)^g=\sigma^i(C^g), \tag{1.1}
\end{equation}
that means for any subsimplex $C$: the adjacency operation $\sigma^i$ ($ i\in \{0,1,2\}$ first) and mapping $g$ commute. So, $\sigma^i$'s act on left, symmetries act on right.
\item 
The $\sigma^i$ operation, as above ($i\in \{0,1,2\}$ at the first considerations),
can be considered as a local reflection in the $i-$side face of any subsimples $C$.
Then any symmetry mapping $g$ (if it exists) can be given by a starting subsimplex 
and its image subsimplex and by simplex-wise extension, so that (1.1) holds further 
step-by-step, as requirement. The axiomatically defined symmetric $D-$matrix function 
$(m^{ij}=m^{ji})$ and the consequences for $D-$symbol ``theory" (e.g. in \cite{MPS06}, \cite{Mo11})
can guarantee the extension of $g$ first to the first simplex $T$, then 
-- introducing $\sigma^3$ operation for the tile adjacency in the next step -- onto the whole tiling $\mathcal{T}$.
\item
The face pairing generator symmetry mappings of $T$, 
i.e. $g_1,g_2,\dots$ (may be involutive one as well), introduce 
also $\sigma^3$ adjacency operations and generate a simplex tiling
$\mathcal{T}={T}^G$ under a symmetry group $G$, finitely generated by the previous 
$g_1,g_2,\dots$. The $D-$matrix function has free parameters yet for $m^{23}=m^{32}$ 
entries, so we shall
have infinite series of tilings for a given $D-$diagram with certain free rotational orders. The $D$-matrix function $_kD\mapsto m^{ij}(_kD)$ of the $D$-symbol says, by requirement,
that the $ij-$edge at the meeting of $\sigma^i$, $\sigma^j$ side faces of any subsimplex $_kC$ in the $D-$set of the $D-$symbol, schould be surrounded by $2m^{ij}(_kD)$
subsimplices in the tiling $\mathcal{T}={T}^G$. (Note, e.g. for the $3-$simplex tiling $\mathcal{T}$ that $m^{01}=m^{12}=3$ constant, and 
$m^{02}=m^{03}=m^{13}=2$ constant, because of the barycentric subdivision of $\mathcal{T}$.)
\item
Around the half-edge $G-$equivalence classes we get subsimplices providing defining relations for $G$ (besides the involution relations), just in the sense of Poincar\'e.
However, it is our main consideration, we can take the $D^0$ subsymbols, 
obtained by leaving $\sigma^0$ operations from the $D-$diagram 
(and the $D-$matrix entries), 
the components describing the fundamental domains of the vertex stabilizer 
subgroups for the vertex equivalence classes induced, as a cutting $2-$surface.
The curvature formula (\cite{MPS06}, \cite{Mo11}) says whether a vertex stabilizer 
(depending on the free parameters) is either a finite spherical group 
(of positive curvature), or Euclidean group (of zero curvature), or hyperbolic one 
(of negative curvature). 
\item
In this paper we concentrate on the last cases, where the parameters involve only outer vertices of hyperbolic space $\mathbb{H}^3$.
Than we cut (truncate) this vertex(ices) orthogonally by polar plane(s) in $\mathbb{H}^3$. The vertex stabilizer group acts then, preserving the two half spaces of the polar plane(s).
Introducing half-space changing isometries to this $\mathbb{H}^2-$stabilizer, preserving the original simplex tiling $\mathcal{T}$, so we get a space group with compact fundamental
domain. The procedure induces face pairing at truncating vertex polars as well.
\item
This new face pairings can be formulated by a new adjacency operation 
$\sigma^{0 *}$ commuting with other $\sigma^j$'s $(j\in \{1,2,3\})$ but not with $\sigma^0$.
Computer can also help, but we have made careful case-by-case discussions, first for later experiences.
\item
A $D-$symbol isomorphism leads to equivariant group extension. Morphism does it onto 
a smaller $D-$diagram, or opposite procedures, lead to symmetry breaking members in the same family,
where the free parameters need extra attention (see \cite{MPS06} and \cite{Mo11}).
\end{enumerate}

\subsection{Results}

In Section 2 we illustrate our method with some cases from family F1, then the leading representative group series will be detailed from each family (Sections 2-5),
sometimes with characteristic examples, also accompanied with figures. Finally we have proved the following summary 

\begin{theorem}
\label{th:mth}
For integer parts of the trunc-simplices of families F1 - F4 there are given the 73 group extensions all together. Namely,

\noindent F1: Representing series  $^{*233} \Gamma (u) = \: ^{\bar{4}3m}_{\ \ 24} \Gamma (u)$ have 14 extensions,

\noindent F2: Representing series $^{2*2} \Gamma _3 (u, 2v) = \: ^{\bar{4}2m}_{\ \ \ 8} \Gamma (u, 2v)$ have 21 extensions,

\noindent F3: Representing series  $^{*33} \Gamma (2u, v) =  \: ^{3m}_{\ \ 6} \Gamma (2u, v)$ have 17 extensions,

\noindent F4: Representing series  $^{*22} \Gamma_1 (u, 2v, w) =  \: ^{mm2}_{\ \ \ \ 4} \Gamma_1 (u, 2v, w)$ have 21 extensions.
\end{theorem}

There are given all the space group series in figures and by group presentations, namely, the face pairing generators and defining relations to 
the edge-segment equivalence classes in Sections 2-5. The other non-fundamental series for families F5 - F12 will be published later.

The truncation plane of a trunc-simplex will be a natural base plane for a hypersphere (constant distance surface 
from the base plane in both sides), whose ball (the inner part of the hypersphere) with its images under the symmetry group of the corresponding trunc-simplex
tiling, can fil (pack) the whole space $\mathbb{H}^3$ without common interior point, or can cover $\mathbb{H}^3$ without gap, respectively.

As a by-product, we have determined the optimal, densest packing for each family from 
F1-F4, moreover, the optimal thinnest (loosest) covering for families F1-F2.
Now this is straightforward, because the Coxeter-Schl\"afli (extended) reflection groups 
$\{\overline{u},\overline{v},\overline{w}\}$ 
occur for these families. See tables to Section 6. 
However, the coverings for families F3-F4 deserve a separate publication with worser densities 
than the optimal covering at F2.

\begin{theorem}
\label{th:pac}
Among the hyperball packings and coverings to the truncation base planes of hyperballs in the trunc-simplex tilings in our families 
F1-F4 there is a densest packing in F1 to group $^{*233} \Gamma$ ($u=7$) (as complete Coxeter group $\{3, 3, 7\}$) 
with density $\approx 0.82251$ (in Fig.~\ref{fig:F1-1-2}); and there is a loosest covering in F2 to group $^{2*2} \Gamma$ 
($u=7, 2v=6$) (extended $\{7, 3, 7\}$) with density $\approx 1.26829$ (in Fig.~\ref{fig:F2-1-2} and Fig.~25).
\end{theorem}
We conjecture (see also investigations \cite{Sz13-4, Sz17-1,Sz17-2,Sz18,Sz19-4,Sz19} of the third author in References) that we have found 
among all arrangements of congruent hyperballs in $\mathbb{H}^3$ the densest hyperball packing with density $\approx 0.82251$ above 
(related to Dirichlet-Voronoi cell subdivision) to the regular trunc-simplex group $\{3, 3, 7\}$; furthermore the loosest hyperball covering with density $\approx 1.26829$ above, 
to the extended trunc-orthoscheme group $\{7,3,7\}$. This last result is published first time here (surprising a little bit). This covering density is less
than that in the conjecture of L.~Fejes~T\'oth \cite{FTL} for balls and horoballs.
\section{Family F1} \label{sec:2}
The maximal series of family F1 is characterized by the simplex tiling whose group 
series is $^{\bar{4}3m}_{\ \ 24} \Gamma (u)$ with previous (crystallographic) 
notation from \cite{MPS06}, and now with the new (orbifold) notation 
$^{*233} \Gamma (u)$ (see \cite{LMV18}). It has the trivial extension 
by the reflection $a_0 = \bar{m}_0$ (Fig.~1:F1-1)
$$^{*233} \Gamma (u) = ( m_0, m_1, m_2, m_3 - 1 = m_0^2 = m_1^2 = m_2^2 = m_3^2 = (m_0 m_1)^3 = (m_0 m_2)^2 = $$
$$ = (m_0 m_3)^2 = (m_1 m_2)^3 = (m_1 m_3)^2 = (m_2 m_3)^u; \ 6 < u \in \mathbb{N} ) $$
$$\mathrm{with} \ \mathrm{extenion} \ ( a_0 - 1 = a_0^2 \stackrel{1}{=} (a_0 m_1)^2 \stackrel{2}{=} (a_0 m_2)^2 \stackrel{3}{=} (a_0 m_3)^2 ) .$$

\vspace{3mm}

{\bf Other non-fundamental cases in F1:}

\vspace{3mm}

\noindent $\bullet \  ^{233}\Gamma _1(2u) = (r, r_2, m_3 - 1 = r^3 = r_2^2 = m_3^2 = m_3 r m_3 r^{-1} = (r r_2)^3 = (m_3 r_2 m_3 r_2)^u, \\ 3 < u \in \mathbb{N}) = \ ^{23}_{12}\Gamma _1(2u)$, with 2 extensions

\noindent {\bf (1)} $(\bar{m}_0, \bar{m}_1 - 1 = \bar{m}_0^2 = \bar{m}_1^2 \stackrel{1}{=} \bar{m}_0 r^{-1} \bar{m}_1 r \stackrel{2}{=} \bar{m}_0 r_2 \bar{m}_1 r_2 \stackrel{3}{=} (\bar{m}_0 m_3)^2 \stackrel{4}{=} (\bar{m}_1 m_3)^2)$ 

\noindent {\bf(2)} $(\bar{s} - 1 \stackrel{1}{=} (\bar{s} r)^2 \stackrel{2}{=} (\bar{s} r_2)^2 \stackrel{3}{=} m_3 \bar{s} m_3 \bar{s}^{-1})$.
\

\begin{figure}[htbp]
	\centering
\includegraphics[width=11cm]{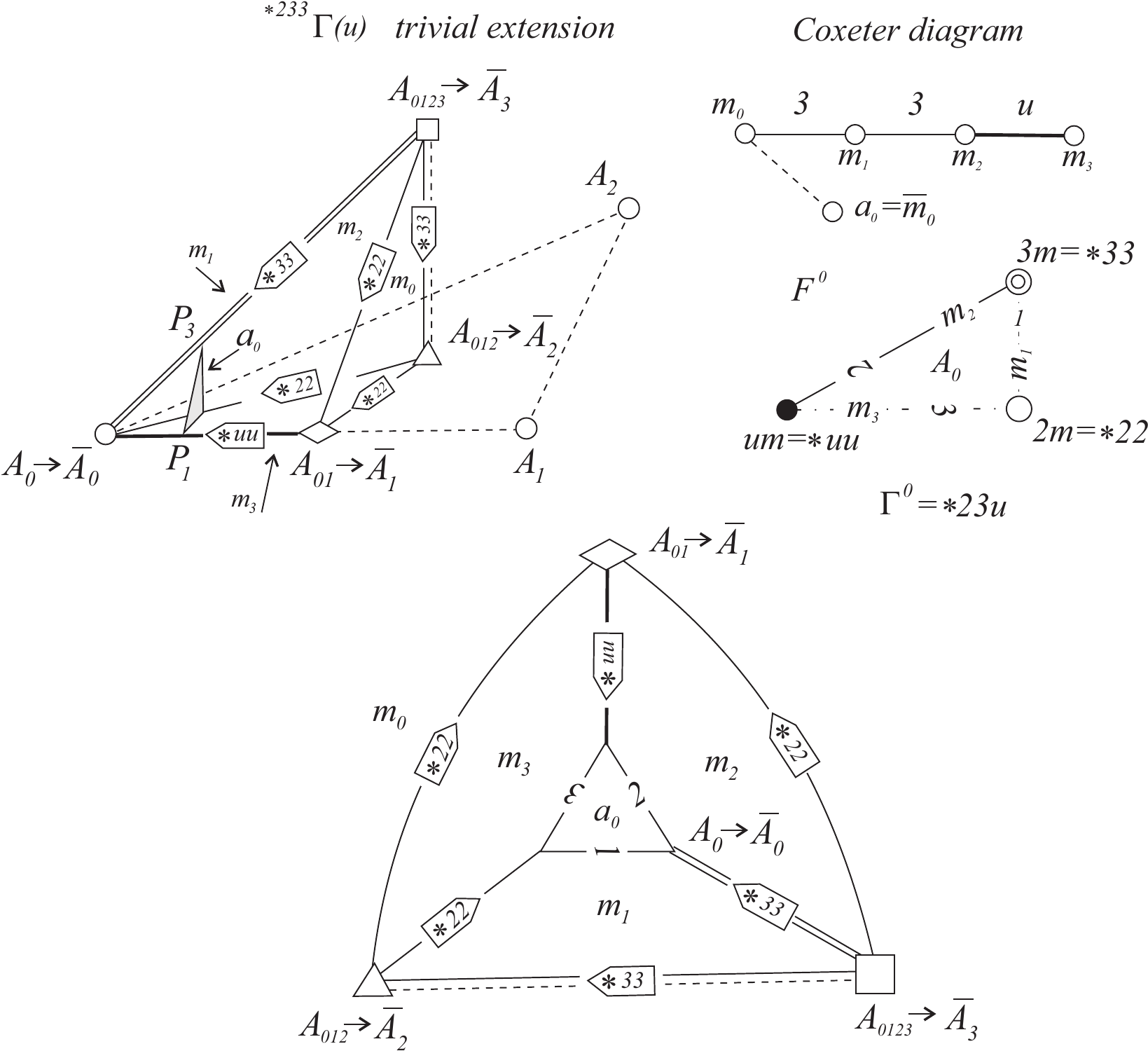}  \includegraphics[width=11cm]{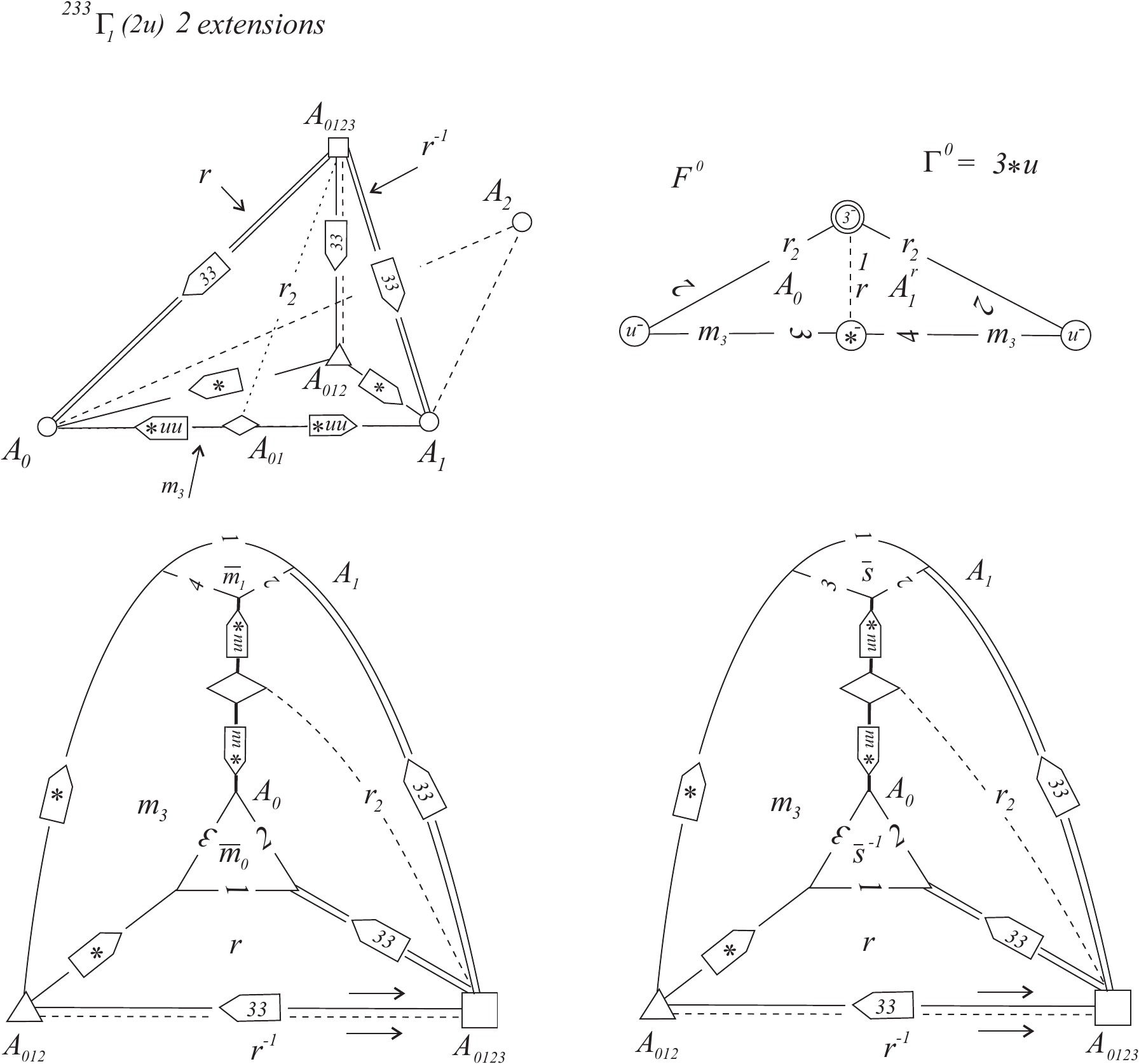}
\caption{F1-1 to $ ^{*233}\Gamma(u)$ and F1-2 to $ ^{233}\Gamma_1(2u)$}
\label{fig:F1-1-2}
\end{figure}

The pictures are given in Fig.~\ref{fig:F1-1-2}: F1-2 
with $\Gamma ^0(A_0, A_1) = 3*u$ as $\mathbb{H}^2$ group with half-turn extension.
\vspace{3mm}

\noindent $\bullet \ ^{233}\Gamma _2(u) = (r, r_2, r_3 - 1 = r^3 = r_2^2 = r_3^2 = (r r_2)^3 = (r r_3)^2 = (r_2 r_3)^u, 6 < u \in \mathbb{N}) = \ ^{23}_{12}\Gamma _2(u)$, with 2 extensions

\noindent {\bf (1)} $(\bar{m}_0, \bar{m}_1 - 1 = \bar{m}_0^2 = \bar{m}_1^2 \stackrel{1}{=} \bar{m}_0 r^{-1} \bar{m}_1 r \stackrel{2}{=} \bar{m}_0 r_2 \bar{m}_1 r_2 \stackrel{3}{=} \bar{m}_0 r_3 \bar{m}_1 r_3)$ 

\noindent {\bf(2)} $(\bar{s} - 1 \stackrel{1}{=} (\bar{s} r)^2 \stackrel{2}{=} (\bar{s} r_2)^2 \stackrel{3}{=} (\bar{s} r_3)^2)$.

\begin{figure}[htbp]
	\centering
\includegraphics[width=11cm]{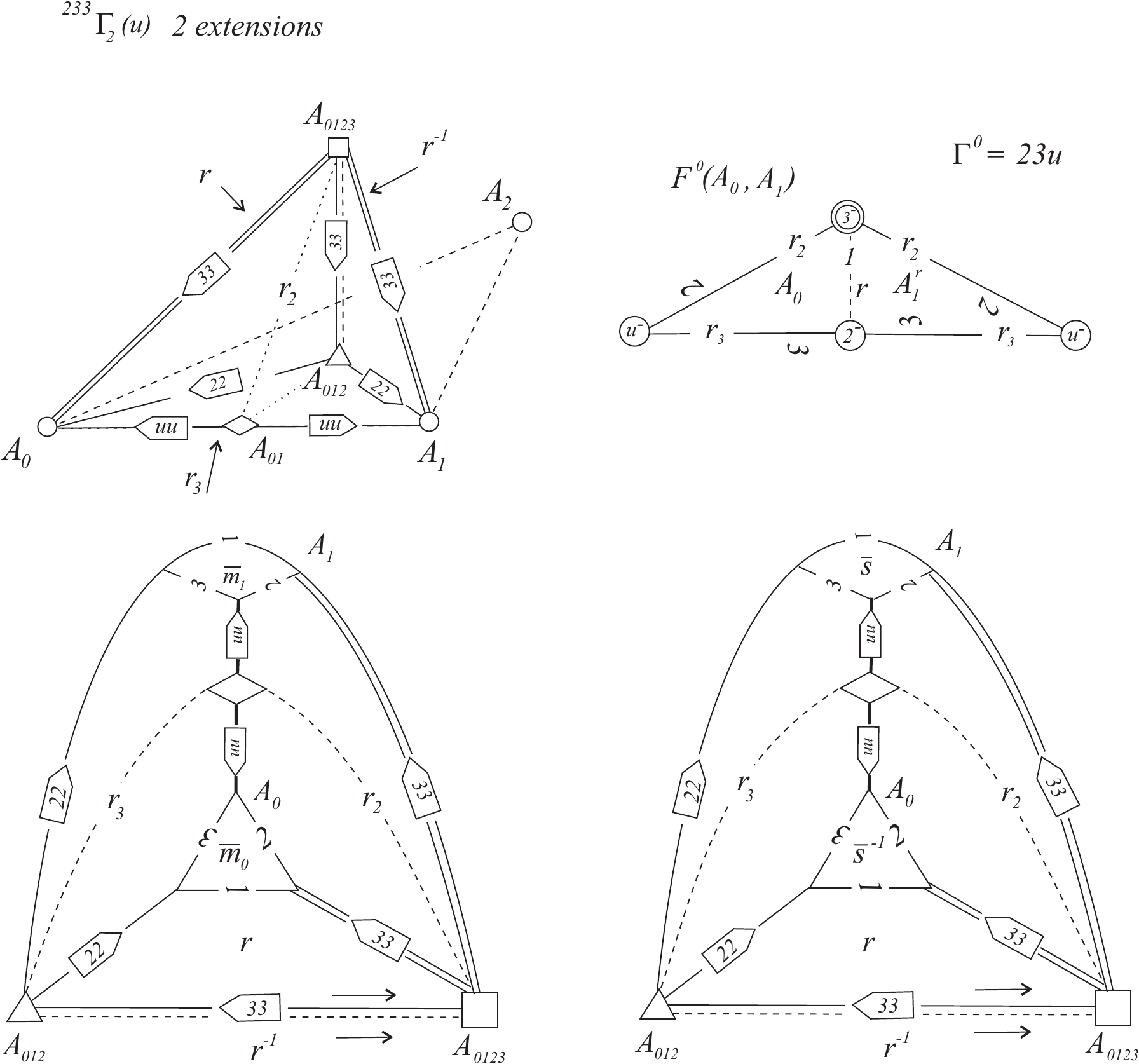} \includegraphics[width=11cm]{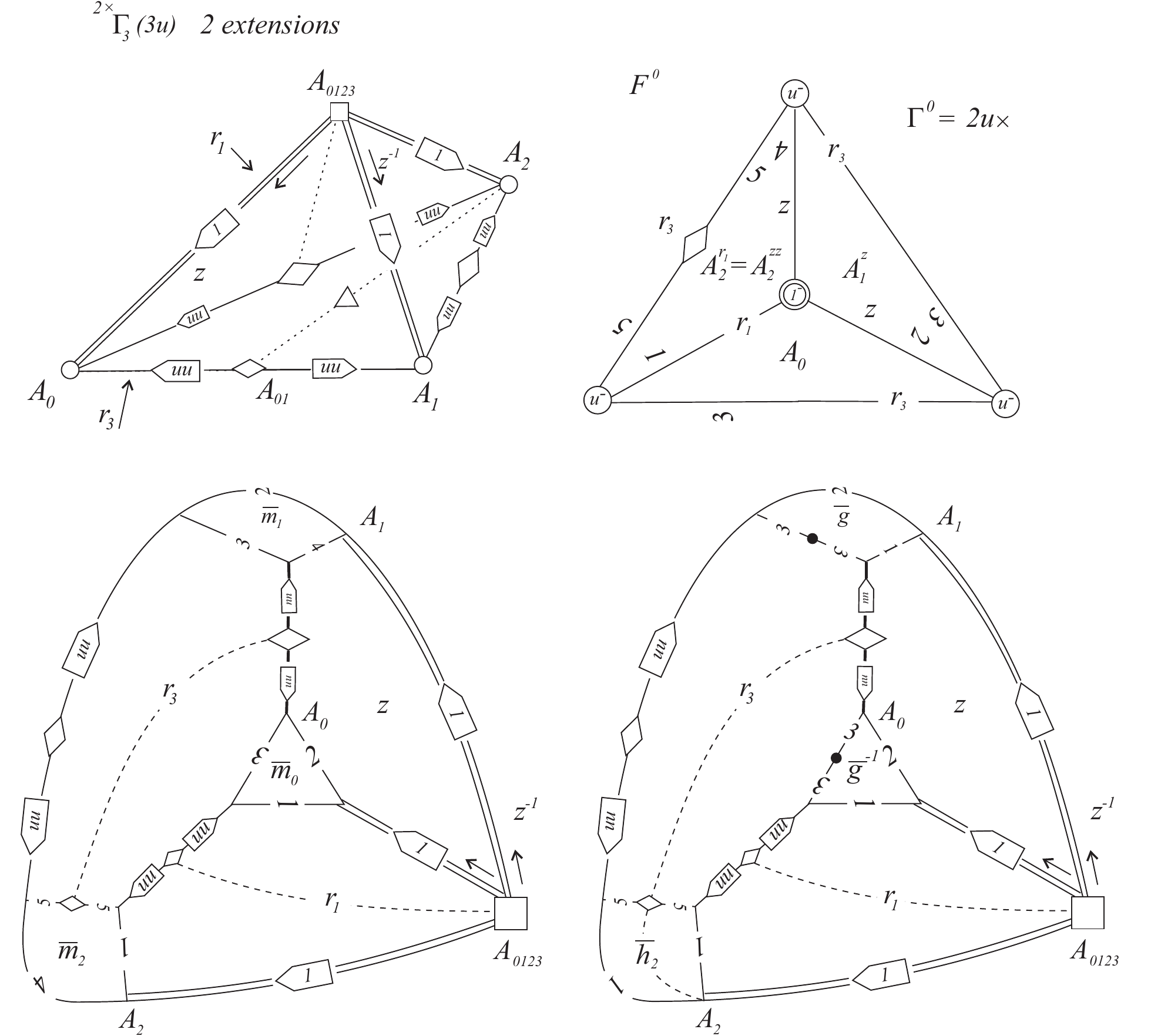} 
\caption{F1-3 to $ ^{233}\Gamma_2(u)$ and F1-4 to $ ^{2\times}\Gamma_3(3u)$}
\label{fig:F1-3-4}
\end{figure}

In Fig.~\ref{fig:F1-3-4}: F1-3 
the pictures are given with $\Gamma ^0 = 23u$ as $\mathbb{H}^2$ group with half-turn extension.

\vspace{3mm}

\noindent $\bullet \ ^{2 \times }\Gamma _3(3u) = (z, r_1, r_3 - 1 = r_1^2 = r_3^2 = z z r_1 = (r_1 r_3 z r_3 z^{-1} r_3)^u, 2 < u \in \mathbb{N}) = \ ^{\bar{4}}_{4}\Gamma _3(3u)$, with 2 extensions

\noindent {\bf (1)} $(\bar{m}_0, \bar{m}_1, \bar{m}_2 - 1 = \bar{m}_0^2 = \bar{m}_1^2 = \bar{m}_2^2 \stackrel{1}{=} \bar{m}_0 r_1 \bar{m}_2 r_1 \stackrel{2}{=} \bar{m}_0 z^{-1}  \bar{m}_1 z \stackrel{3}{=} \bar{m}_0 r_3 \bar{m}_1 r_3 \stackrel{4}{=} \bar{m}_1 z^{-1}  \bar{m}_2 z \stackrel{5}{=} (\bar{m}_2 r_3)^2)$ 

\noindent {\bf(2)} $(\bar{g}, \bar{h}_2 - 1 = \bar{h}_2^2 \stackrel{1}{=} \bar{g} z^{-1} \bar{h}_2 r_1 \stackrel{2}{=} (\bar{g} z)^2 \stackrel{3}{=} (\bar{g} r_3)^2 \stackrel{5}{=} 
(\bar{h}_2 r_3)^2 )$.

The pictures are given in Fig.~\ref{fig:F1-3-4}: F1-4 
with $\Gamma ^0 = 2u \times $ as $\mathbb{H}^2$ group with half-turn extension.

\vspace{3mm}

\noindent $\bullet \ ^{33}\Gamma _5(4u) = (r, r_2, r_3 - 1 = r^3 = r_2^2 = r_3^2 = (r r_2)^2 = (r_3 r r_3 r^{-1} r_3 r_2)^u, \frac{3}{2} < u \in \mathbb{N}) = \ ^{3}_{3}\Gamma _5(4u)$, with 2 extensions

\noindent {\bf (1)} $(\bar{m}_0, \bar{m}_1, \bar{m}_2 - 1 = \bar{m}_0^2 = \bar{m}_1^2 = \bar{m}_2^2 \stackrel{1}{=} \bar{m}_0 r^{-1} \bar{m}_1 r \stackrel{2}{=} \bar{m}_0 r_2 \bar{m}_1 r_2 \stackrel{3}{=} \\ \bar{m}_0 r_3 \bar{m}_2 r_3 \stackrel{4}{=} (\bar{m}_1 r_3)^2 \stackrel{5}{=} \bar{m}_2 r \bar{m}_2 r^{-1})$ 

\noindent {\bf(2)} $(\bar{h}_0, \bar{h}_1, \bar{h}_2 - 1 = \bar{h}_0^2 = \bar{h}_1^2 = \bar{h}_2^2 \stackrel{1}{=} \bar{h}_0 r_2 \bar{h}_1 r \stackrel{3}{=} \bar{h}_0 r_3 \bar{h}_2 r_3  \stackrel{4}{=} (\bar{h}_1 r_3)^2 \stackrel{5}{=} (\bar{h}_2 r)^2 )$. 

\begin{figure}[htbp]
	\centering
\includegraphics[width=11cm]{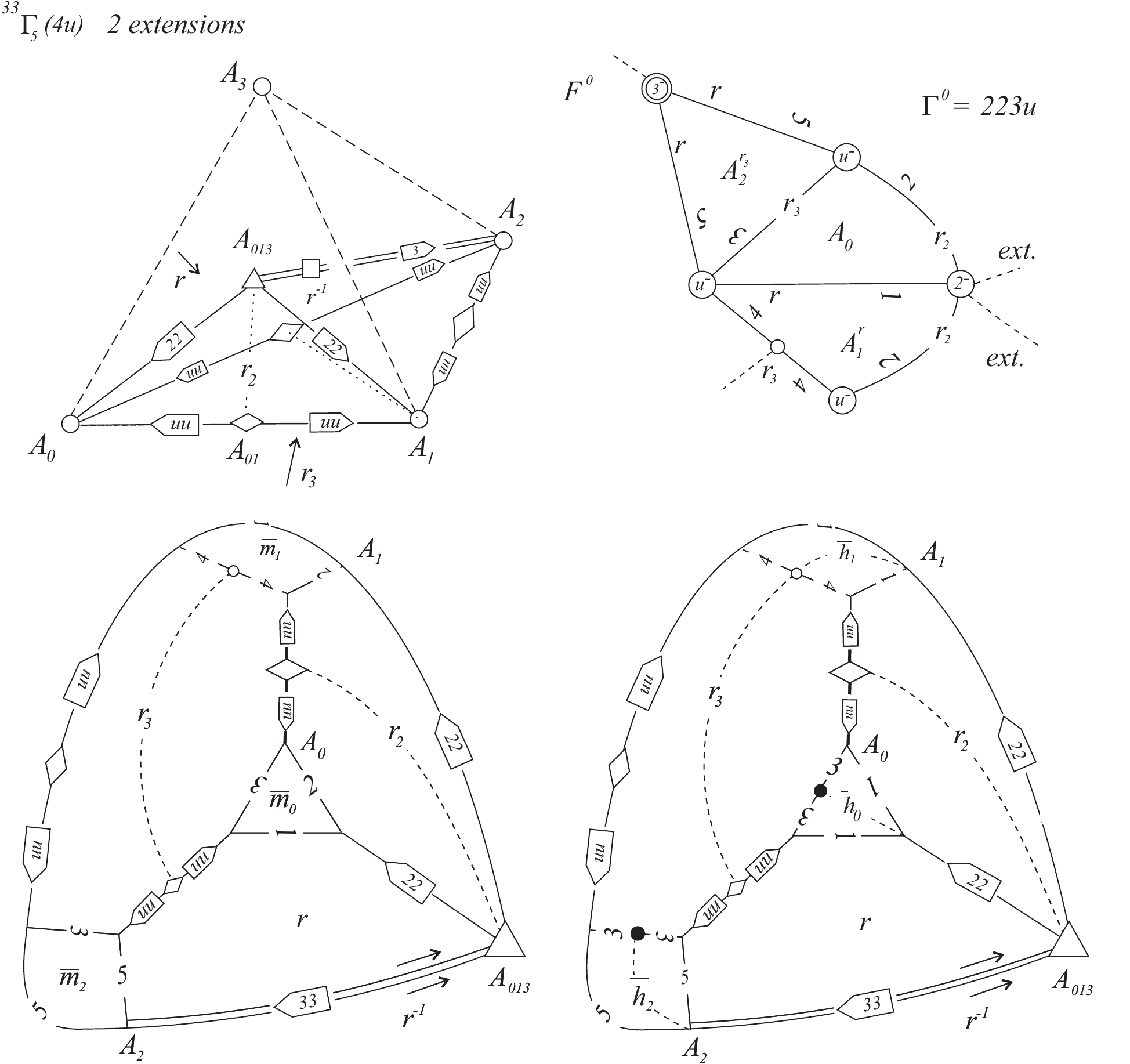} \includegraphics[width=11cm]{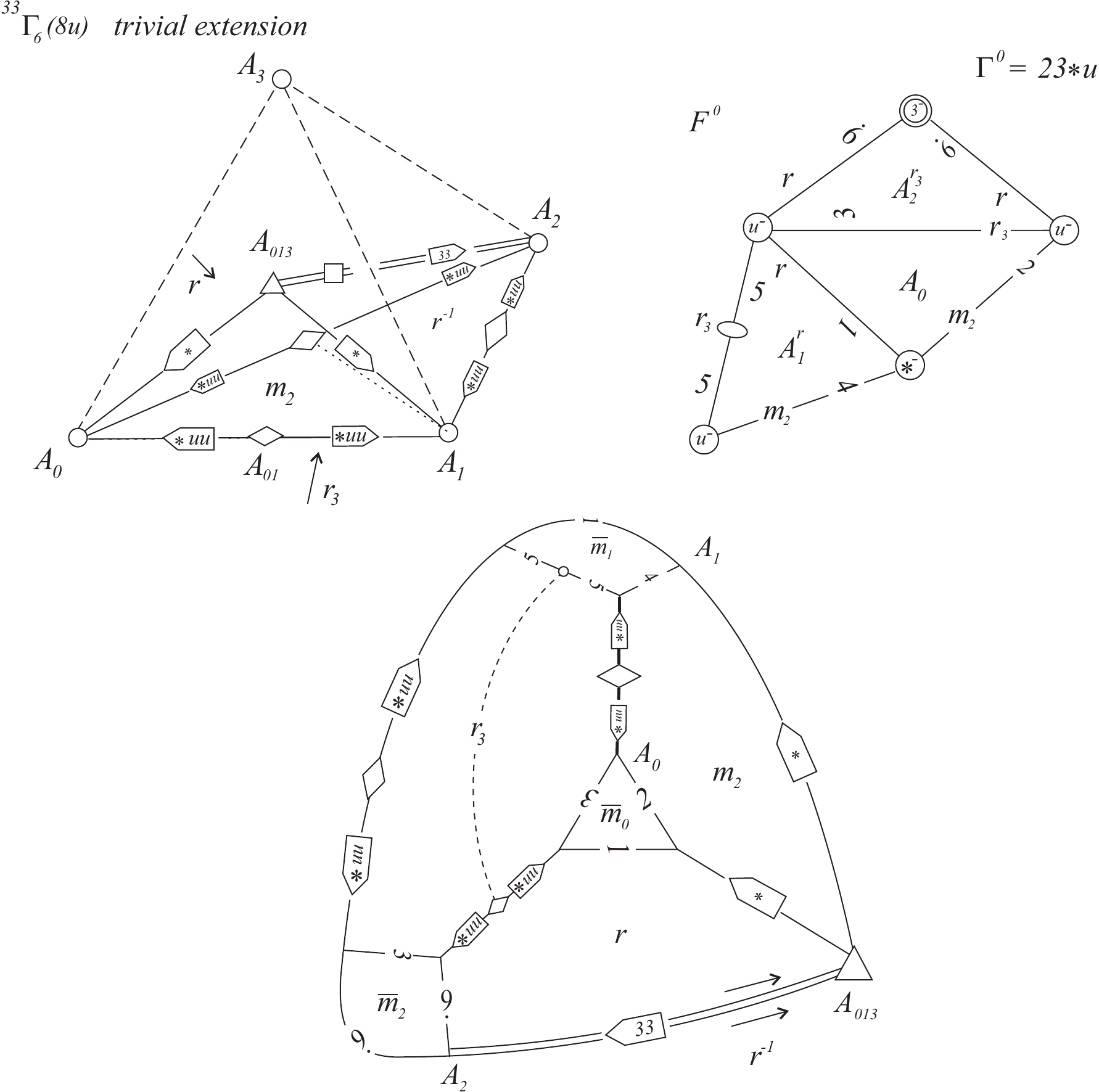}
\caption{F1-5 to $ ^{33}\Gamma_5(4u)$ and F1-6 to $ ^{33}\Gamma_6(8u)$}
\label{fig:F1-5-6}
\end{figure}

The pictures with $\Gamma ^0 = 223 $ as $\mathbb{H}^2$ group with two half-turn extensions are given in Fig.~\ref{fig:F1-5-6}: F1-5.

\vspace{3mm}

\noindent $\bullet \ ^{33}\Gamma _6(8u) = (r, m_2, r_3 - 1 = r^3 = m_2^2 = r_3^2 = m_2 r m_2 r^{-1} = (m_2 r_3 r r_3 r^{-1} r_3)^{2u}, \frac{3}{4} < u \in \mathbb{N}) = \ ^{3}_{3}\Gamma _6(8u)$, with trivial extension

\noindent {\bf (1)} $(\bar{m}_0, \bar{m}_1, \bar{m}_2 - 1 = \bar{m}_0^2 = \bar{m}_1^2 = \bar{m}_2^2 \stackrel{1}{=} \bar{m}_0 r^{-1} \bar{m}_1 r \stackrel{2}{=} 
(\bar{m}_0 m_2)^2 \stackrel{3}{=} \bar{m}_0 r_3 \bar{m}_3 r_3 \stackrel{4}{=} (\bar{m}_1 m_2)^2 \stackrel{5}{=} (\bar{m}_1 r_3)^2 \stackrel{6}{=}  \bar{m}_2 r \bar{m}_2 r^{-1})$.

See Fig.~\ref{fig:F1-5-6}: F1-6
for the pictures with $\Gamma ^0 = 23*u$. 

\vspace{3mm}

\begin{figure}[htbp]
	\centering
\includegraphics[width=14cm]{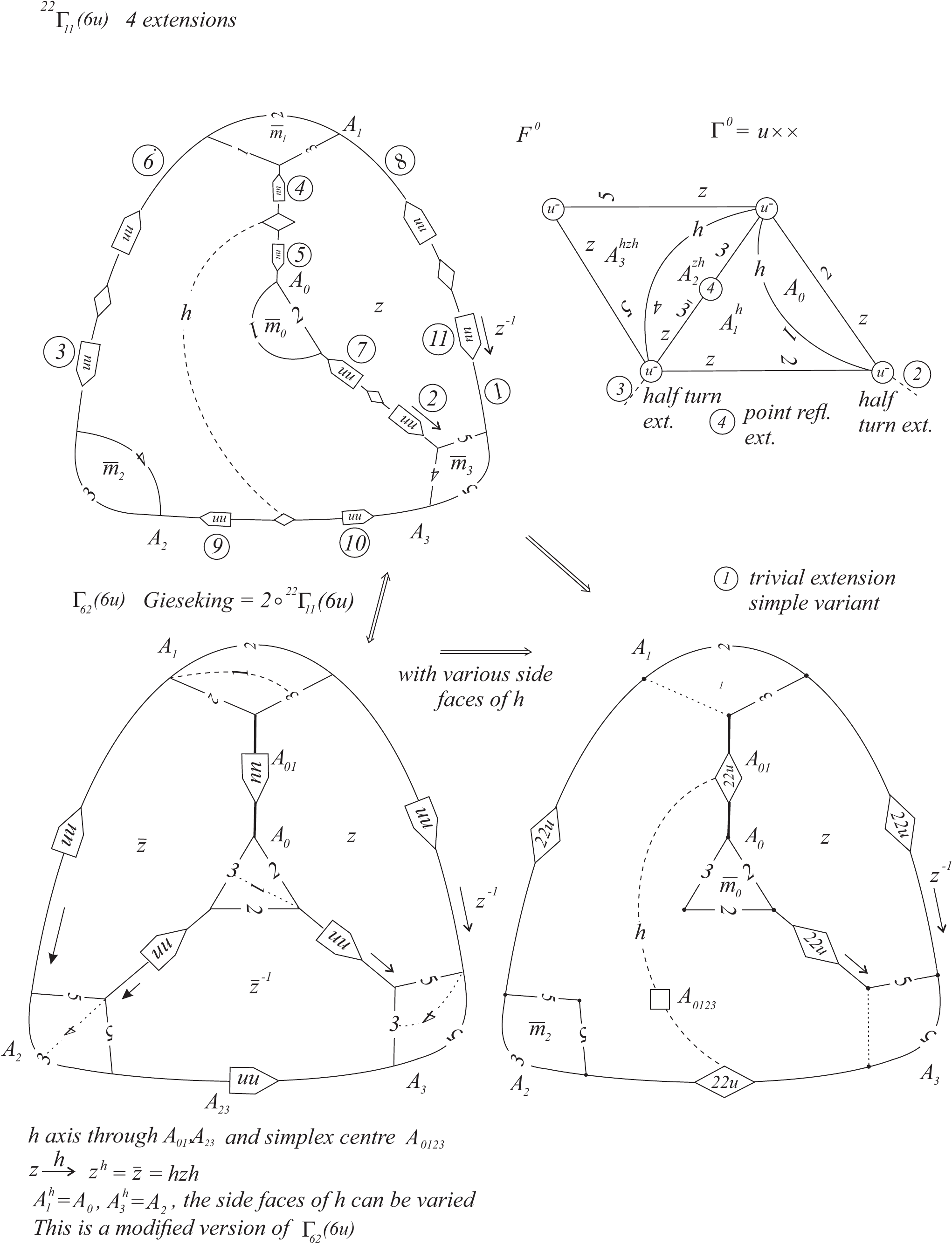} 
\caption{F1-$7_1$ to $ ^{22}\Gamma_{11}(6u)$ as Gieseking orbifolds}
\label{fig:F1-7}
\end{figure}

\begin{figure}[htbp]
	\centering
\includegraphics[width=14cm]{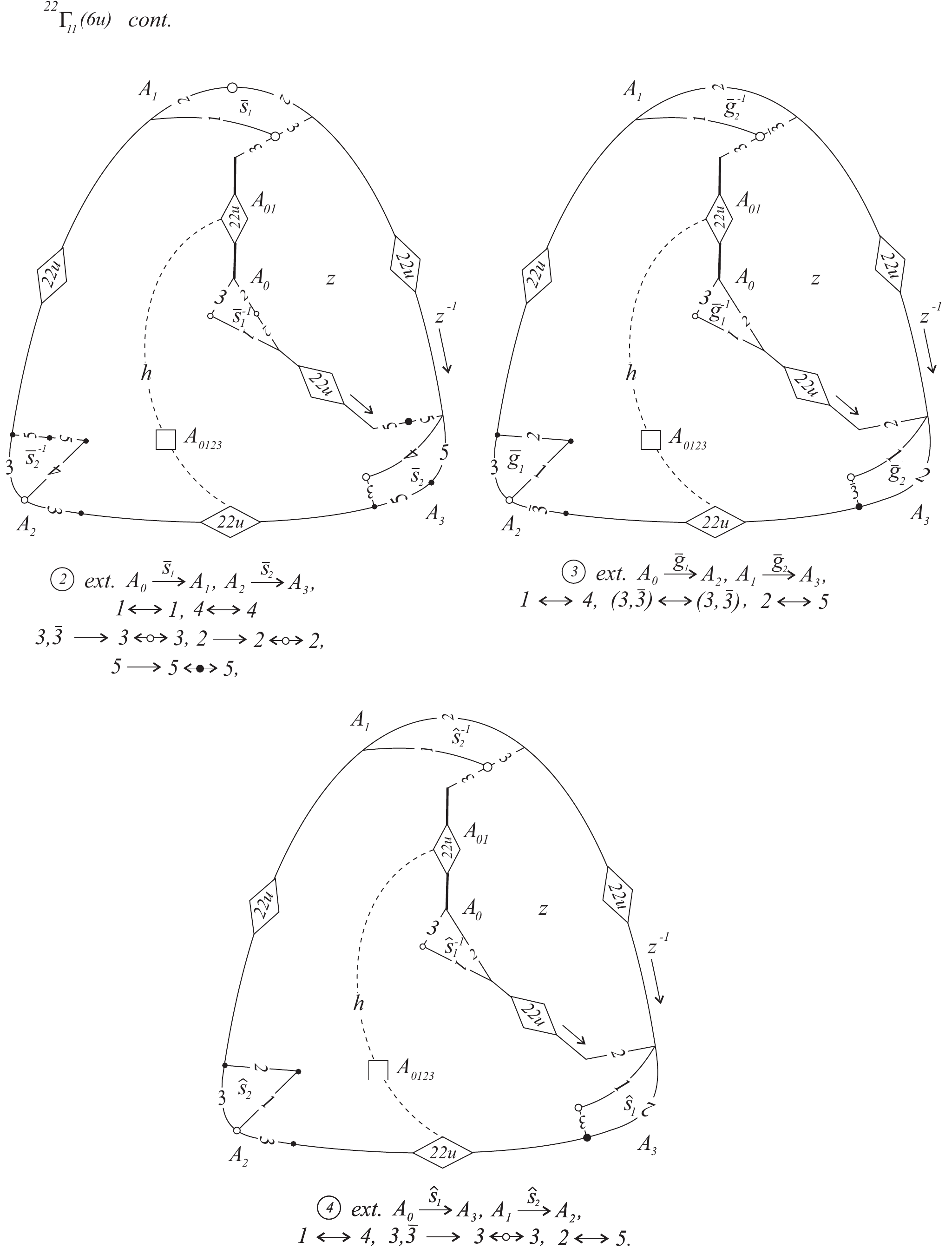} 
\caption{F1-$7_2$ to $ ^{22}\Gamma_{11}(6u)$, first version of ext. (2)-(4)}
\label{fig:F1-7-c}
\end{figure}

\begin{figure}[htbp]
	\centering
\includegraphics[width=14cm]{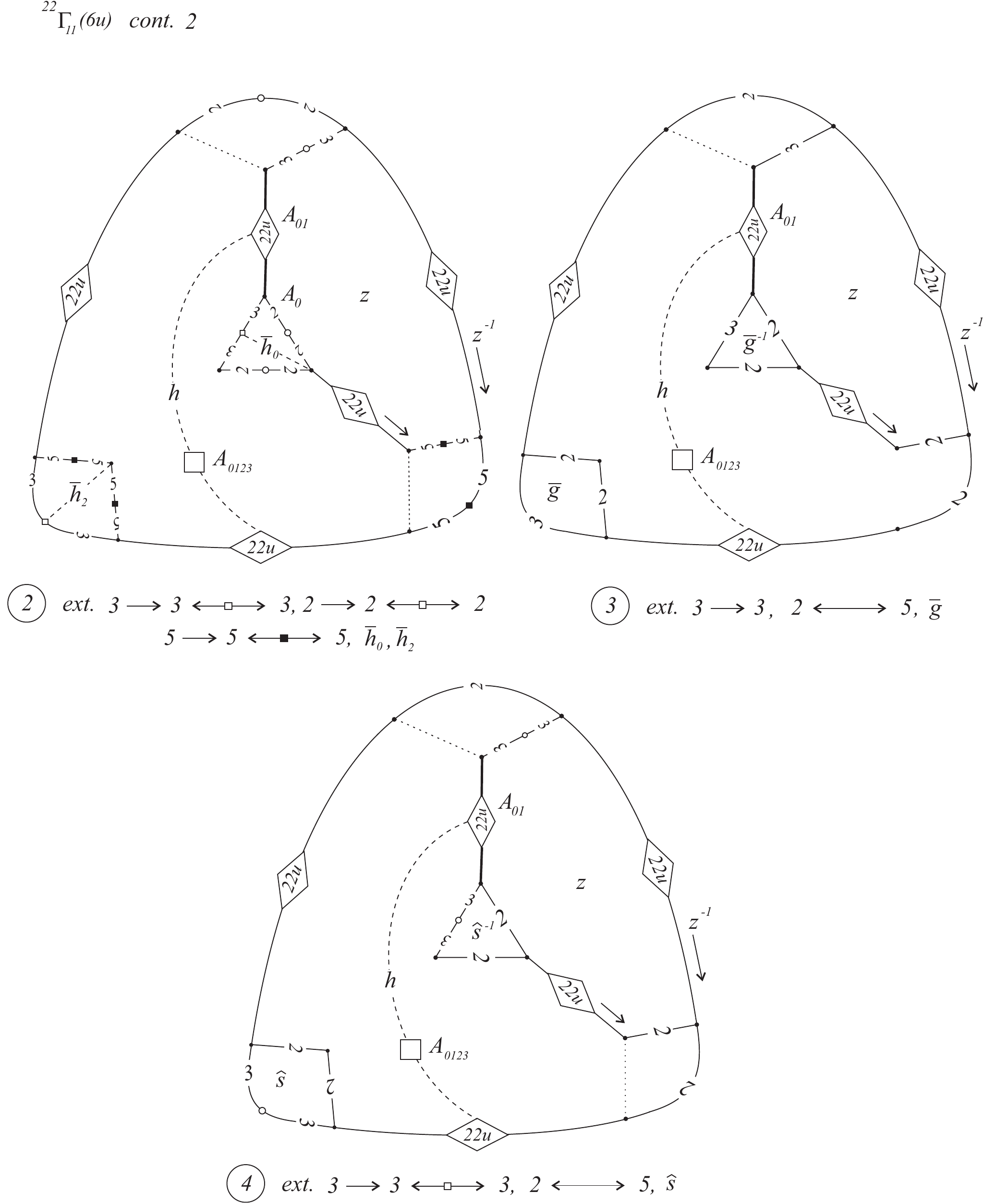} 
\caption{F1-$7_3$ to $ ^{22}\Gamma_{11}(6u)$, second simpler version of ext. (2)-(4)}
\label{fig:F1-7-c2}
\end{figure}

\noindent $\bullet \ ^{22}\Gamma _{11}(6u) = (z, h - 1 = h^2 = (h z z h z h z^{-1} h z^{-1} z^{-1})^u, \ 1 < u \in \mathbb{N}) = \ ^{2}_{2}\Gamma _{11}(6u)$, 

\noindent with 4 extensions, each in two variants, given in Fig.~\ref{fig:F1-7} - \ref{fig:F1-7-c2} to F1-7. These so-called Gieseking orbifolds are derived from the famous non-orientable manifold with an infinite vertex class, originally if $u = 1$ (cusp, but $1<u$ now). 
The half-turn extension $h$ allows various half simplices and fundamental domains, e.g. for the trivial extension as follows (Fig.~\ref{fig:F1-7}: F1-$7_1$).

\noindent {\bf (1)} $(\bar{m}_0, \bar{m}_1, \bar{m}_2, \bar{m}_3 - 1 = \bar{m}_0^2 = 
\bar{m}_1^2 = \bar{m}_2^2 = \bar{m}_3^2 \stackrel{1}{=} \bar{m}_0 h \bar{m}_1 h 
\stackrel{2}{=} \bar{m}_0 z^{-1} \bar{m}_1 z \stackrel{3}{=} \bar{m}_1 z^{-1} \bar{m}_2 z 
\stackrel{4}{=} \bar{m}_2 h \bar{m}_3 h \stackrel{5}{=} \bar{m}_3 z \bar{m}_3 z^{-1} ).$

Fig.~\ref{fig:F1-7-c2}: F1-$7_3$ shows the simpler version of extensions (2)-(4) by (additional) presentations:

\noindent {\bf (2)} $(\bar{h}_0, \bar{h}_2 - 1 = \bar{h}_0^2 = \bar{h}_2^2 \stackrel{2}{=} (\bar{h}_0 z^{-1} h)^2  \stackrel{3}{=} \bar{h}_0 h z^{-1} \bar{h}_2 z h \stackrel{5}{=} (\bar{h}_2 h z h)^2)$,

\noindent {\bf (3)} $(\bar{g} - 1 \stackrel{2}{=} \bar{g} h z h \bar{g}^{-1} z^{-1} h \stackrel{3}{=} (\bar{g} z h)^2)$,

\noindent {\bf (4)} $(\hat{s} - 1 \stackrel{2}{=} \hat{s} h z h \hat{s}^{-1} h z 
\stackrel{3}{=} (\hat{s} z h)^2 $).


\subsection{The 4 extensions to Gieseking orbifolds as fundamental trunc-simplex tilings}
\begin{figure}[htbp]
	\centering
\includegraphics[width=13cm]{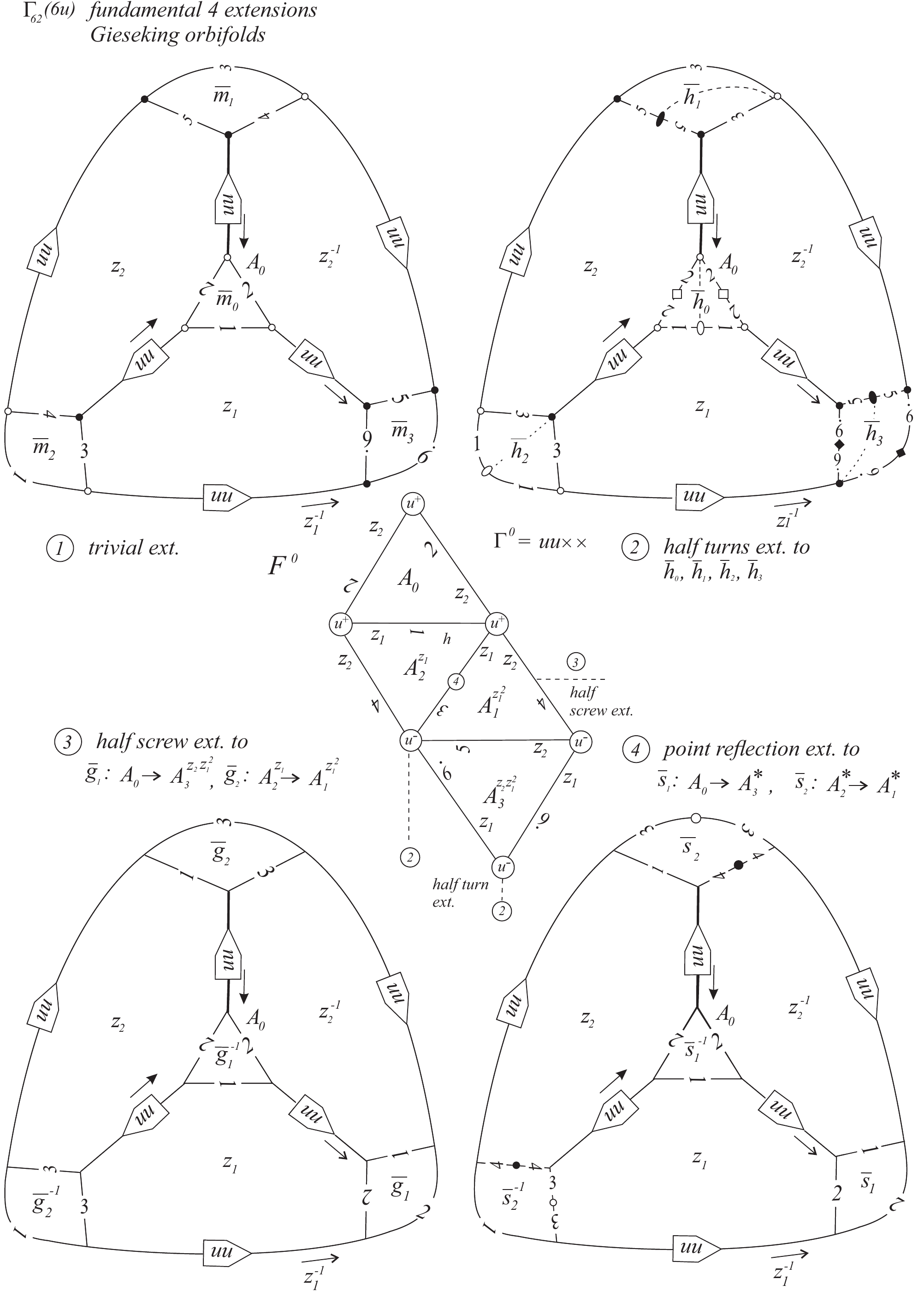} 
\caption{F1-7+ to $ \Gamma_{62}(6u)$, the $4$ extensions to Gieseking fundamental orbifolds}
\label{fig:F1-7+}
\end{figure}

Our Fig.~\ref{fig:F1-7+}: F1-7+ to $\Gamma_{62}(6u)$ shows this famous simplex 
face pairing, as repetition from \cite{S19}, with orientation reversing transforms 
$$
\mathbf{z}_1:z_1^{-1}\rightarrow z_1~\text{fixing vertex $A_3$ and}
~\mathbf{z}_2:z_2^{-1}\rightarrow z_2~\text{fixing $A_0$}
$$

All vertices lie at the infinity (absolute) of $\mathbb{H}^3$, iff the ``rotation order'' parameter at the edge equivalence class is $u=1$. 
At the polar plane of $A_0$ (say, e.g.) we look the fundamental domain $F^0$ of the vertex class stabilizer subgroup $\Gamma^0=uu\!\!\times\!\!\times$ 
that is the Euclidean Klein bottle group $\mathbf{4.pg}$ (with two cross caps $\times$, 
or projective planes), iff $u=1$. 
Moreover, we get interesting ``surgery effects" that lead e.g. to the famous 
Fomenko-Matveev-Weeks
hyperbolic space form of minimal volume \cite{V}.

Otherwise, iff $1<u \in \mathbb{N}$ we get $\mathbb{H}^2$ plane group (area of $F^0$ is $\frac{4\pi}{u}$ less then $4\pi$ \cite{LMV18}).
Then the generators, reversing the half-spaces of the polar plane $a_0$, but preserving the original simplex tiling $\mathcal{T}$, lead to $4$ group
extensions with compact fundamental domains as trunc-simplices (or octahedra $\mathcal{O}^1_{62},\dots \mathcal{O}^4_{62})$.

Presentation of the original simplex tiling group is 
\begin{equation}
\Gamma_{62}(6u)=(z_1,z_2~-~1=(z_1z_2z_2z_1^{-1}z_2^{-1}z_1)^u;~1<u\in \mathbb{N}). \notag
\end{equation}
\begin{enumerate}
\item[(\bf{1})] Then comes the trivial extension with plane reflections:																		  
\begin{equation}
\begin{gathered}
(\overline{m}_0,\overline{m}_1,\overline{m}_2,\overline{m}_3~-~1=\overline{m}_0^2=\overline{m}_1^2=\overline{m}_2^2=\overline{m}_3^2 \stackrel{1}{=} 
\overline{m}_0 z_1^{-1}\overline{m}_2z_1\stackrel{2}{=}\overline{m}_0 z_2\overline{m}_0 z_2^{-1} 
\stackrel{3}{=} \\ \stackrel{3}{=}
\overline{m}_1 z_1\overline{m}_2z_1^{-1}\stackrel{4}{=}\overline{m}_1 z_2\overline{m}_2z_2^{-1}\stackrel{5}{=}\overline{m}_1 z_2^{-1}\overline{m}_3z_2
\stackrel{6}{=}\overline{m}_3 z_1\overline{m}_3z_1^{-1}). \notag
\end{gathered}
\end{equation}
\item[(\bf{2})] 
The half turns extension, changing the half spaces at $F^0$ (inducing point reflections at edges $2,6$, moreover, a $2-$knot $u-$segment combination):
\begin{equation}
\begin{gathered}
(\overline{h}_0,\overline{h}_1,\overline{h}_2,\overline{h}_3~-~1=\overline{h}_0^2=\overline{h}_1^2=\overline{h}_2^2=\overline{h}_3^2 \stackrel{1}{=} 
\overline{h}_0 z_1^{-1}\overline{h}_2z_1\stackrel{2}{=}(\overline{h}_0 z_1)^2 \stackrel{3}{=} \\ \stackrel{3}{=}
\overline{h}_1 z_1\overline{h}_2z_2^{-1}\stackrel{5}{=}\overline{h}_1 z_2^{-1}\overline{h}_3z_2 \stackrel{6}{=} (\overline{h}_3 z_1)^2. \notag
\end{gathered}
\end{equation}
\item[(\bf{3})] The most interesting half-screw extension, changing the half-spaces at $F^0$ yields a ``non-orientable $u-$knot" (by glide reflections $\overline{g}_1,\overline{g}_2$):
\begin{equation}
\begin{gathered}
(\overline{g}_1,\overline{g}_2~-~1 \stackrel{1}{=} 
\overline{g}_1 z_2\overline{g}_2^{-1} z_1 \stackrel{2}{=}\overline{g}_1 z_1g_1^{-1}z_2 \stackrel{3}{=} \overline{g}_2 z_1\overline{g}_2z_2). \notag
\end{gathered}
\end{equation}
\item[(\bf{4})] The $4$th extension is a point reflection at edges $3,4$ of $F^0$ (inducing screw motions $\overline{s}_1,\overline{s}_2$ and $u-$knot):
\begin{equation}
\begin{gathered}
(\overline{s}_1,\overline{s}_2~-~1\stackrel{1}{=} 
\overline{s}_1 z_2\overline{s}_2^{-1} z_1 \stackrel{2}{=}\overline{s}_1 z_1s_1^{-1}z_2^{-1} \stackrel{3}{=} (\overline{s}_2 z_1)^2 \stackrel{4}{=} (\overline{s}_2 z_2)^2). \notag
\end{gathered}
\end{equation}
We think that these phenomena are remarkable!
\end{enumerate}

\section{Family 2}
In family F2 the maximal series is characterized by the simplex tiling whose group 
series is $^{2*2} \Gamma (u, 2v)$ (or $^{\bar{4}2m}_{\ \ \ 8} \Gamma (u, 2v)$ by 
notations in \cite{MPS06}). Remember that in case $u = 2v$ we order this tiling 
(and similarly 
the later ones) to family F1 with parameter $\overline{u}=u=2v$, 
since so many simplices (i.e. $2\overline{u}$ barycentric subsimplices) 
surround every simplex edge.

\begin{equation}
\begin{gathered}
^{2*2} \Gamma (u, 2v) = ( m_0, m_2, m_3, r_1 - 1 = m_0^2 = m_2^2 = m_3^2 = r_1^2 = (m_0 m_2)^2 = (m_0 m_3)^2 = \\
= m_0 r_1 m_2 r_1 = (m_2 m_3)^u = (m_3 r_1 m_3 r_1)^v, \ u \neq 2v \in \mathbb{N}, \ \dfrac{1}{u} + \dfrac{1}{v} < \dfrac{1}{2} ) \notag
\end{gathered}
\end{equation}
\noindent with the trivial extension 

$\ (\bar{m}_0, \bar{m}_2 - 1 = \bar{m}_0^2 = \bar{m}_2^2 \stackrel{1}{=} \bar{m}_0 r_1 \bar{m}_2 r_1 \stackrel{2}{=} (\bar{m}_0 m_2)^2
\stackrel{3}{=} (\bar{m}_0 m_3)^2 \stackrel{4}{=} (\bar{m}_2 m_0)^2 \stackrel{5}{=} (\bar{m}_2 m_3)^2).$

\vspace{3mm}

\begin{figure}[htbp]
	\centering
\includegraphics[width=110mm]{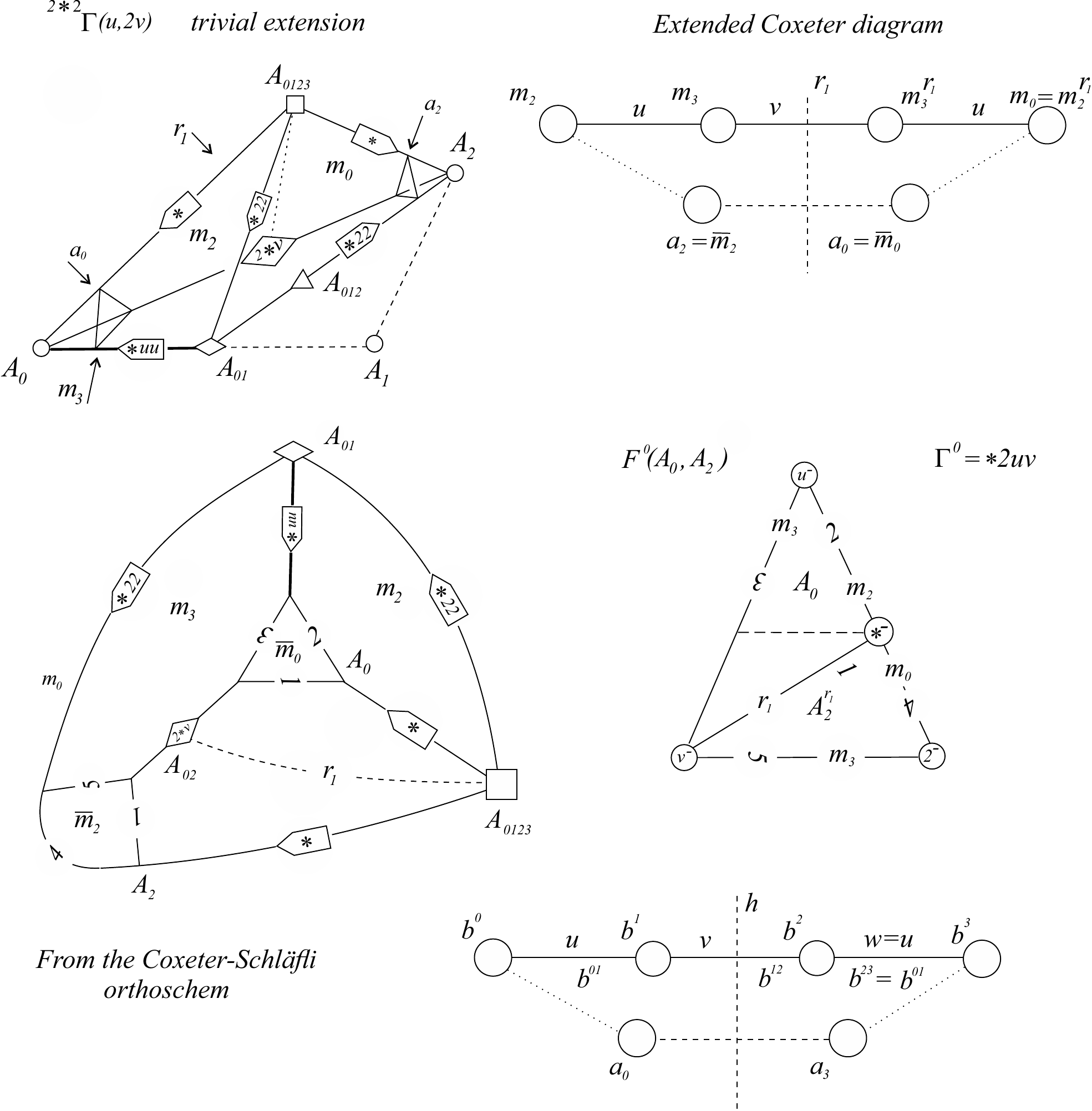} \includegraphics[width=110mm]{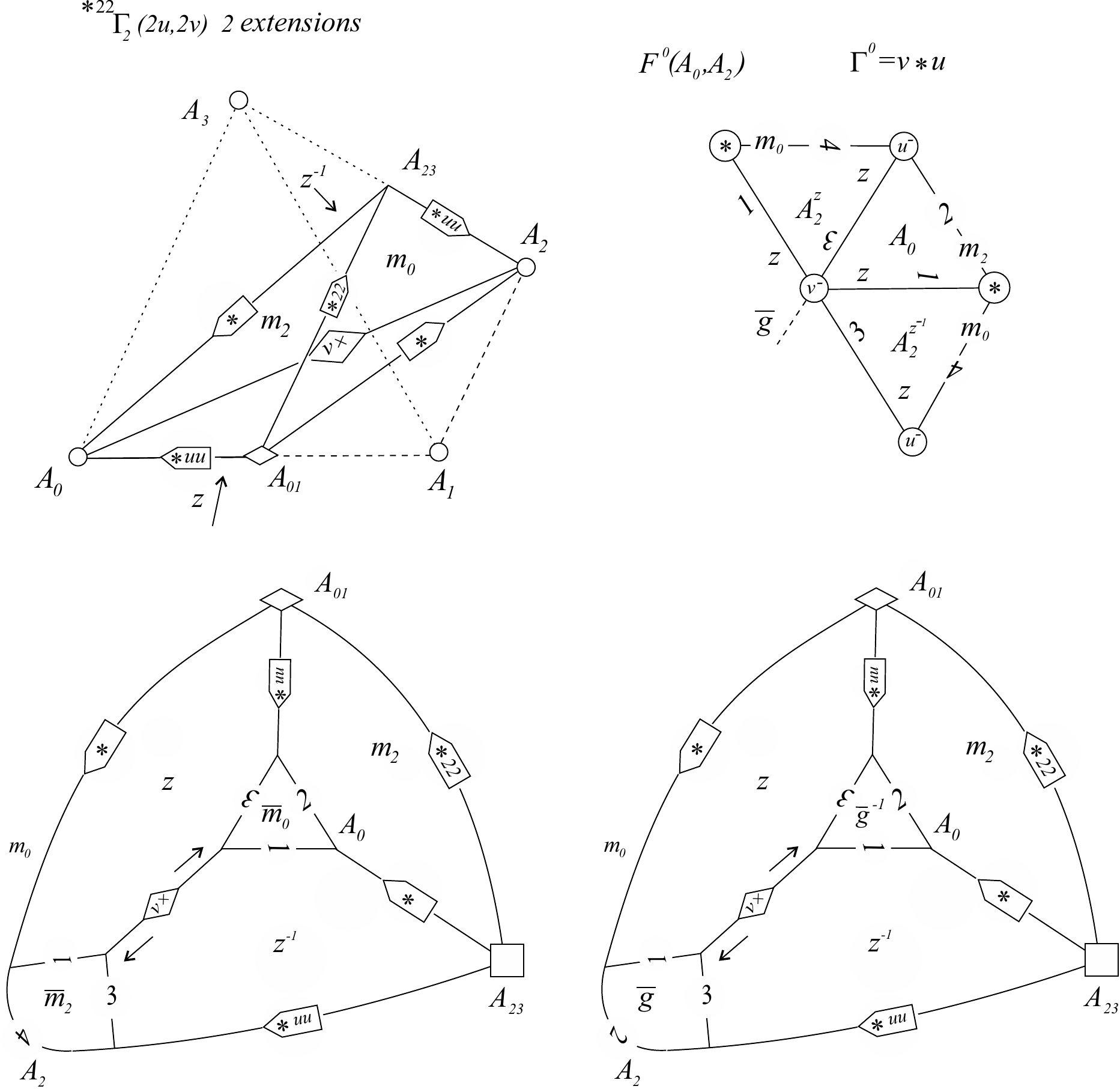}
\caption{F2-1 to $ ^{2*2}\Gamma(u,2v)$ and F2-2 to $ ^{*22}\Gamma_2(2u,2v)$}
\label{fig:F2-1-2}
\end{figure}

\noindent Pictures are given in Fig.~\ref{fig:F2-1-2}: F2-1
with $\Gamma ^0 (A_0, A_2) = *2uv$ as $\mathbb{H}^2$ group.

\vspace{3mm}

{\bf Other non-fundamental cases in F2}

\vspace{3mm}

\noindent $\bullet \  ^{*22}\Gamma _2 (2u, 2v)  = (m_0, m_2, z - 1 = m_0^2 = m_2^2 = (m_0 m_2)^2 = m_2 z m_0 z^{-1} = \\
(m_0 z m_2 z^{-1})^u = (z z)^v, \ u \neq v \in \mathbb{N}, \ \dfrac{1}{u} + \dfrac{2}{v} < 1) = \ ^{mm2}_{\: \ \ \ \ 4} \Gamma_2 (2u, 2v) \\$ with 2 extensions in Fig.~\ref{fig:F2-1-2}: F2-2

\noindent {\bf (1)} $(\bar{m}_0, \bar{m}_2 - 1 = \bar{m}_0^2 = \bar{m}_2^2 \stackrel{1}{=} \bar{m}_0 z \bar{m}_2 z^{-1} \stackrel{2}{=} (\bar{m}_0 m_2)^2 \stackrel{3}{=} \bar{m}_0 z^{-1} \bar{m}_2 z  \stackrel{4}{=} (\bar{m}_2 m_0)^2 )$ 

\noindent {\bf (2)} $(\bar{g} - 1 \stackrel{1}{=} (\bar{g} z^{-1})^2 \stackrel{2}{=} m_2 \bar{g} m_0 \bar{g}^{-1} \stackrel{3}{=} (\bar{g} z)^2 ).$

\noindent $\Gamma ^0 = v*u$ is a $\mathbb{H}^2$ group with half-turn extension yielding glide reflection $\bar{g}$.

\vspace{3mm}

\noindent $\bullet \  ^{2\times} \Gamma _1 (2u, 2v) = (r_2, m_3, z - 1 = r_2^2 = m_3^2 = z z r_2 =
(m_3 r_2 m_3 r_2)^u = (m_3 z m_3 z^{-1})^v, \ u \neq v \in \mathbb{N}, \ \dfrac{1}{u} + \dfrac{2}{v} < 1) = \ ^{\bar{4}}_{4} \Gamma _1 (2u, 2v) $ with 2 extensions

\noindent {\bf (1)} $(\bar{m}_0, \bar{m}_1, \bar{m}_2 - 1 = \bar{m}_0^2 = \bar{m}_1^2 = \bar{m}_2^2 \stackrel{1}{=} \bar{m}_0 z^{-1} \bar{m}_2 z \stackrel{2}{=} \bar{m}_0 r_2 \bar{m}_1 r_2 \stackrel{3}{=} (\bar{m}_0 m_3)^2 \stackrel{4}{=} \bar{m}_1 z \bar{m}_2 z^{-1} \stackrel{5}{=} (\bar{m}_1 m_3)^2 \stackrel{6}{=} (\bar{m}_2 m_3)^2)$ 

\noindent {\bf (2)} $(\bar{h}_2, \bar{s} - 1 = \bar{h}_2^2 \stackrel{1}{=} \bar{s} z \bar{h}_2 z \stackrel{2}{=} (\bar{s} r_2)^2 \stackrel{3}{=} m_3 \bar{s} m_3 \bar{s}^{-1} 
\stackrel{6}{=} (m_3 \bar{h}_2)^2 ).$

\begin{figure}[htbp]
	\centering
\includegraphics[width=110mm]{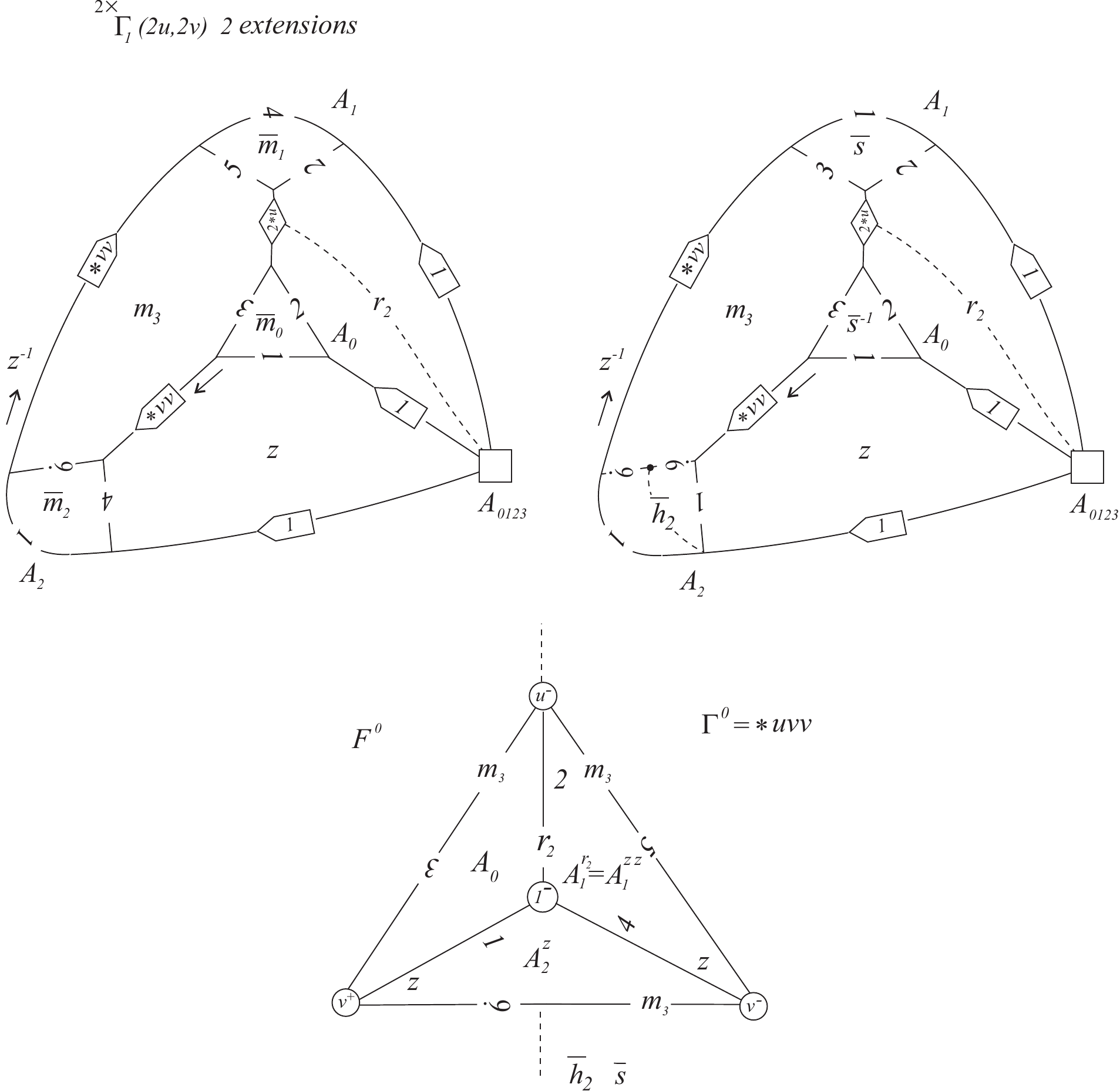} \includegraphics[width=110mm]{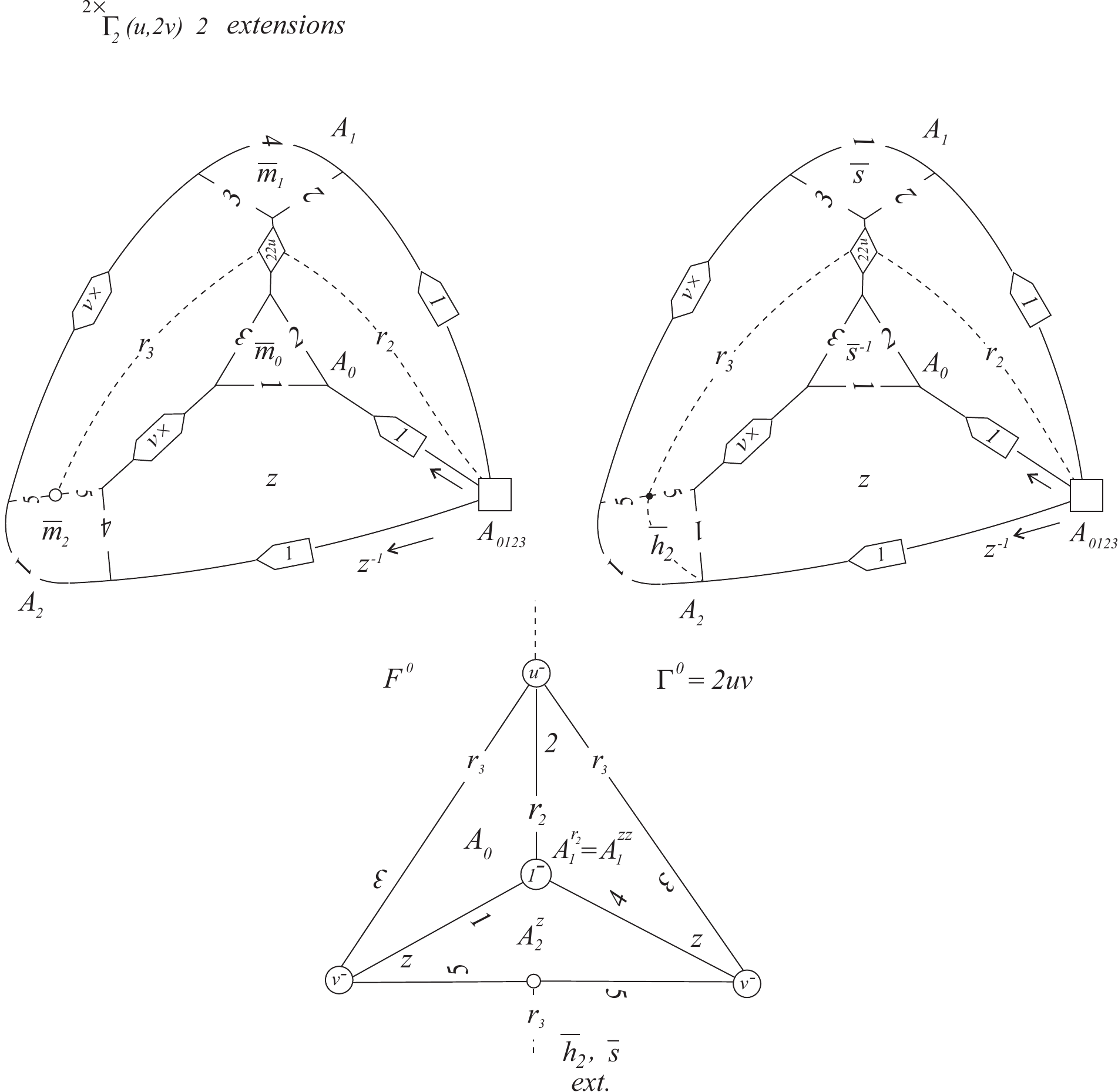}
\caption{F2-3 to $ ^{2\times}\Gamma_1(2u,2v)$ and F2-4 to $ ^{2\times}\Gamma_2(u,2v)$}
\label{fig:F2-3-4}
\end{figure}

\noindent Pictures are given in Fig.~\ref{fig:F2-3-4}: F2-3
with $\Gamma ^0 = *uvv$ as $\mathbb{H}^2$ group with half-turn extension.

\vspace{3mm}

\noindent $\bullet \  ^{2\times} \Gamma _2 (u, 2v) = (r_2, r_3, z - 1 = r_2^2 = r_3^2 = z z r_2 =
(r_2 r_3)^u = (z r_3 z r_3)^v, \ u \neq 2v  \in \mathbb{N}, \ \dfrac{1}{u} + \dfrac{1}{v} < \dfrac{1}{2}) = \ ^{\bar{4}}_{4} \Gamma _2 (u, 2v) $ with 2 extensions

\noindent {\bf (1)} $(\bar{m}_0, \bar{m}_1, \bar{m}_2 - 1 = \bar{m}_0^2 = \bar{m}_1^2 = \bar{m}_2^2 \stackrel{1}{=} \bar{m}_0 z^{-1} \bar{m}_2 z \stackrel{2}{=} \bar{m}_0 r_2 \bar{m}_1 r_2 \stackrel{3}{=} \bar{m}_0 r_3 \bar{m}_1 r_3 \stackrel{4}{=} \bar{m}_1 z \bar{m}_2 z^{-1} \stackrel{5}{=} (\bar{m}_2 r_3)^2)$ 

\noindent {\bf (2)} $(\bar{h}_2, \bar{s} - 1 = \bar{h}_2^2 \stackrel{1}{=} \bar{s} z \bar{h}_2 z \stackrel{2}{=} (\bar{s} r_2)^2 \stackrel{3}{=} 
(\bar{s} r_3)^2 \stackrel{5}{=} (\bar{h}_2 r_3)^2 ).$

\noindent Pictures are given in Fig.~\ref{fig:F2-3-4}: F2-4
with $\Gamma ^0 = 2uv$ as $\mathbb{H}^2$ group with half-turn extension.

\vspace{3mm}

\noindent $\bullet \  ^{222} \Gamma _2 (u, 2v) = 
(r_0, r_1 r_2, r_3 - 1 = r_0^2 = r_1^2 = r_2^2 = 
r_3^2 = r_0 r_1 r_2 = (r_2 r_3)^u = (r_0 r_3 r_1 r_3)^v, \ u \neq 2v \in  \mathbb{N}, \ \dfrac{1}{u} + \dfrac{1}{v} < \dfrac{1}{2}) = \ ^{222}_{\: \ \ 4} \Gamma _2 (u, 2v)$ with 2 extensions

\noindent {\bf (1)} $(\bar{m}_0, \bar{m}_1, \bar{m}_2 - 1 = \bar{m}_0^2 = \bar{m}_1^2 = \bar{m}_2^2 \stackrel{1}{=} \bar{m}_0 r_1 \bar{m}_2 r_1 \stackrel{2}{=} \bar{m}_0 r_2 \bar{m}_1 r_2 \stackrel{3}{=} \\ \bar{m}_0 r_3 \bar{m}_1 r_3 \stackrel{4}{=} \bar{m}_1 r_0 \bar{m}_2 r_1 \stackrel{5}{=} (\bar{m}_2 r_3)^2 )$

\noindent {\bf (2)} $(\bar{h}_2, \bar{s} - 1 = \bar{h}_2^2 \stackrel{1}{=} 
\bar{s} r_0 \bar{h}_2 r_1 \stackrel{2}{=} (\bar{s} r_2)^2 \stackrel{3}{=} 
(\bar{s} r_3)^2 \stackrel{5}{=} (\bar{h}_2 r_3)^2 ).$

\noindent See Fig.~\ref{fig:F2-5-6}: F2-5 with $\Gamma ^0 = 2uv$ as $\mathbb{H}^2$ group with half-turn extension.

\subsection{Truncated half-simplices from F2}

\vspace{3mm}

\noindent $\bullet \  ^{22}\Gamma _7(2u, 4v) = (r_2, r_3, h - 1 = r_2^2 = r_3^2 = h^2 = (r_2 r_3)^u = (r_2 h r_3 h r_3 h r_2 h)^v, \\ u \neq 2v \in \mathbb{N}, \ \dfrac{1}{u} + \dfrac{1}{v} < 1 ) = ^2_2\Gamma _7(2u, 4v)$ with 4 extensions

\noindent {\bf (1)} $(\bar{m}_0, \bar{m}_1, \bar{m}_2, \bar{m}_3  - 1 = \bar{m}_0^2 = \bar{m}_1^2 = \bar{m}_2^2 = \bar{m}_3^2 \stackrel{1}{=} \bar{m}_0 h \bar{m}_3 h \stackrel{2}{=} \bar{m}_0 r_2 \bar{m}_1 r_2 \stackrel{3}{=} \bar{m}_0 r_3 \bar{m}_1 r_3 \stackrel{4}{=} \bar{m}_1 h \bar{m}_2 h \stackrel{5}{=} (\bar{m}_2 r_3)^2 \stackrel{6}{=} (\bar{m}_3 r_2)^2 )$

\noindent {\bf (2)} $(\bar{g}_1, \bar{g}_2 = 1 \stackrel{1}{=} \bar{g}_1 h \bar{g}_2 h \stackrel{2}{=} \bar{g}_1 r_3 \bar{g}_1 r_2 \stackrel{5}{=} r_3 \bar{g}_2 r_2 \bar{g}_2^{-1} )$

\noindent {\bf (3)} $(\bar{s}_1, \bar{s}_2 = 1 \stackrel{1}{=} \bar{s}_1 h \bar{s}_2 h \stackrel{2}{=} (\bar{s}_1 r_2)^2 \stackrel{3}{=} (\bar{s}_1 r_3)^2 \stackrel{5}{=} r_3 \bar{s}_2 r_2 \bar{s}_2^{-1} )$

\noindent {\bf (4)} $(\bar{h}_0, \bar{h}_1, \bar{h}_2, \bar{h}_3  - 1 = 
\bar{h}_0^2 = \bar{h}_1^2 = \bar{h}_2^2 = \bar{h}_3^2 \stackrel{1}{=} 
\bar{h}_0 h \bar{h}_3 h \stackrel{2}{=} \bar{h}_0 r_2 \bar{h}_1 r_3 
\stackrel{4}{=} \bar{h}_1 h \bar{h}_2 h \stackrel{5}{=} (\bar{h}_2 r_3)^2  
\stackrel{6}{=} (\bar{h}_3 r_2)^2 ).$

\begin{figure}[htbp]
	\centering
\includegraphics[width=100mm]{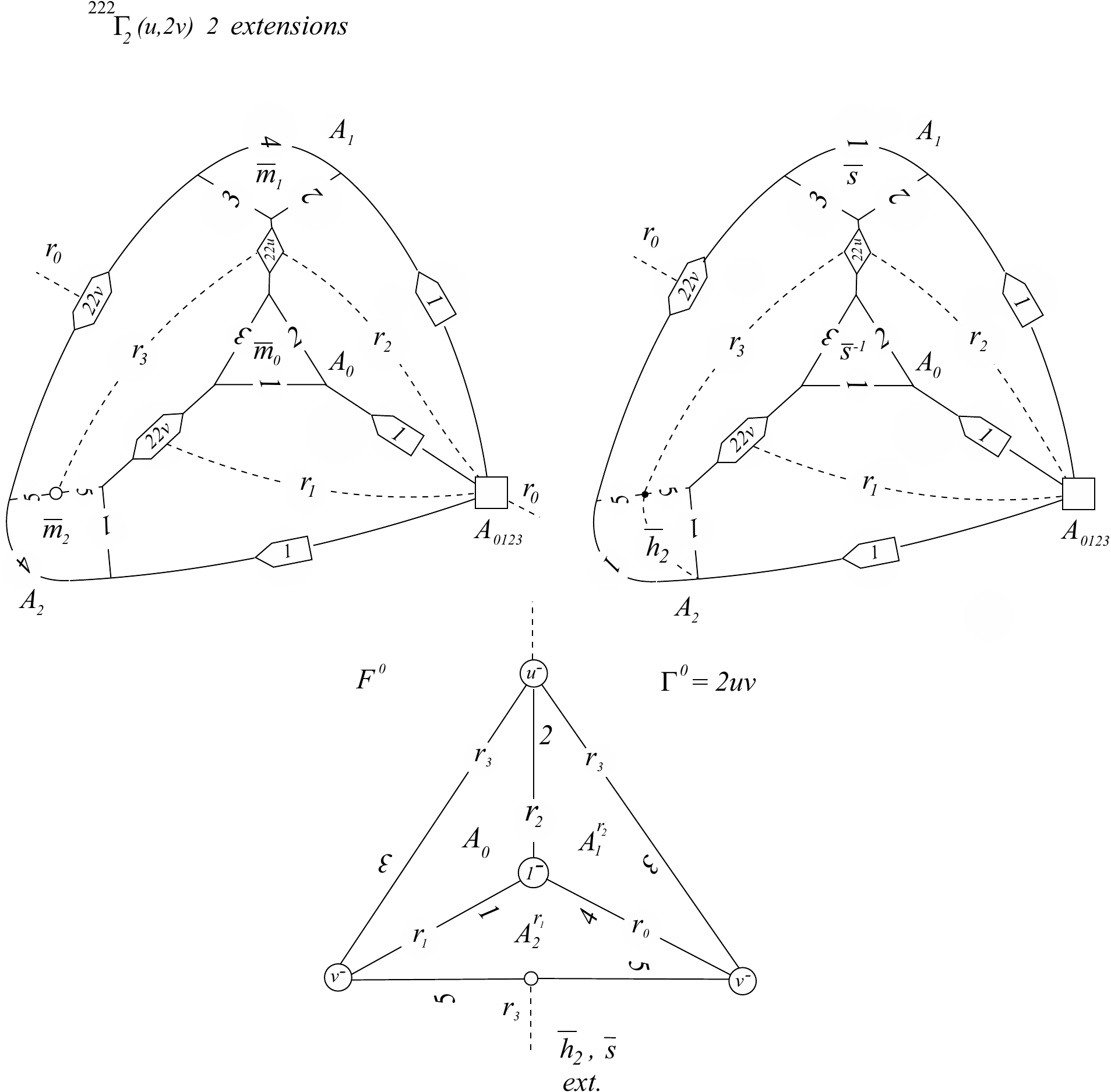} \includegraphics[width=100mm]{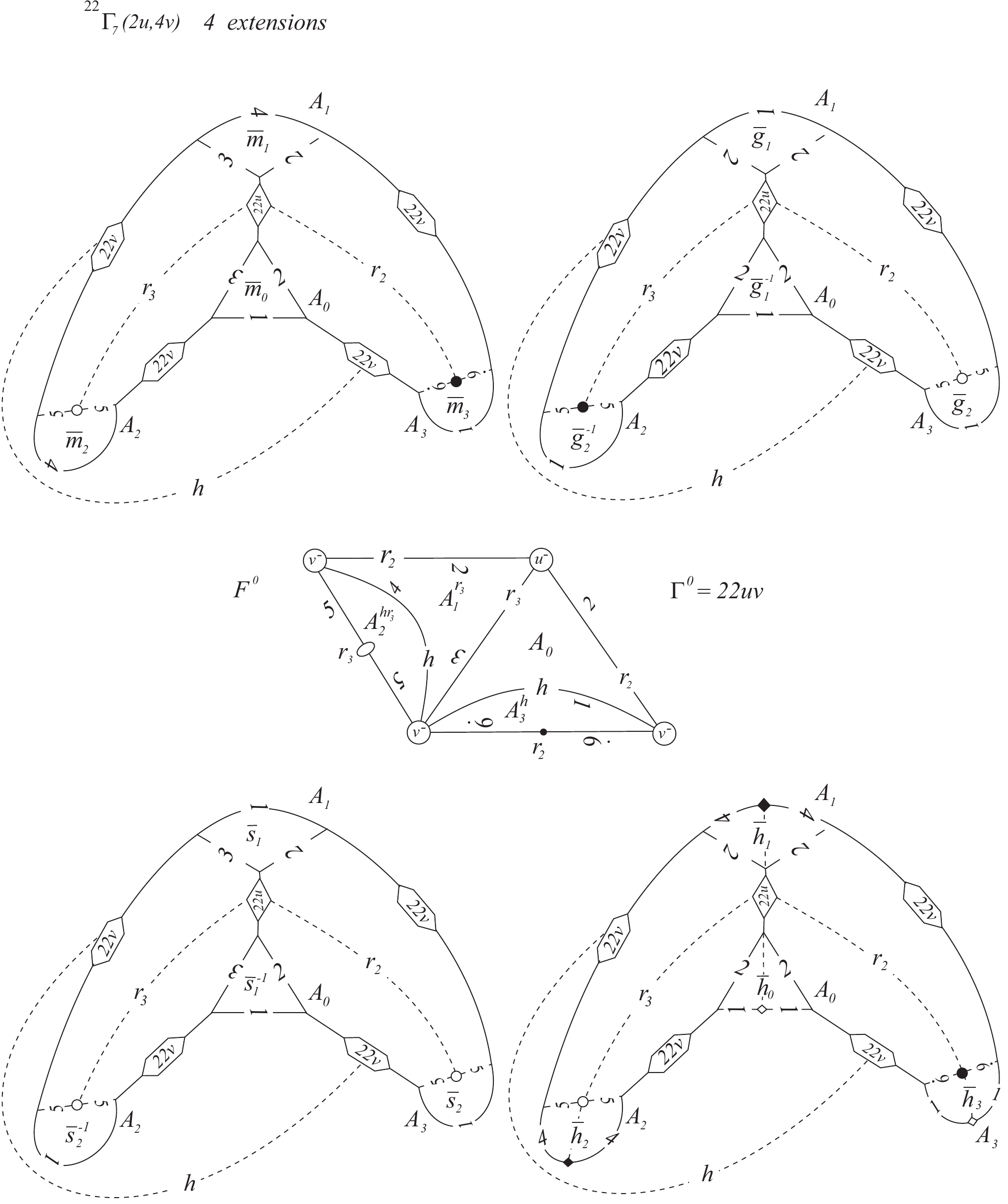}
\caption{F2-5 to $ ^{222}\Gamma_2(u,2v)$ and F2-$6_1$ to $ ^{22}\Gamma_7(2u,4v)$}
\label{fig:F2-5-6}
\end{figure}

\begin{figure}[htbp]
	\centering
\includegraphics[width=130mm]{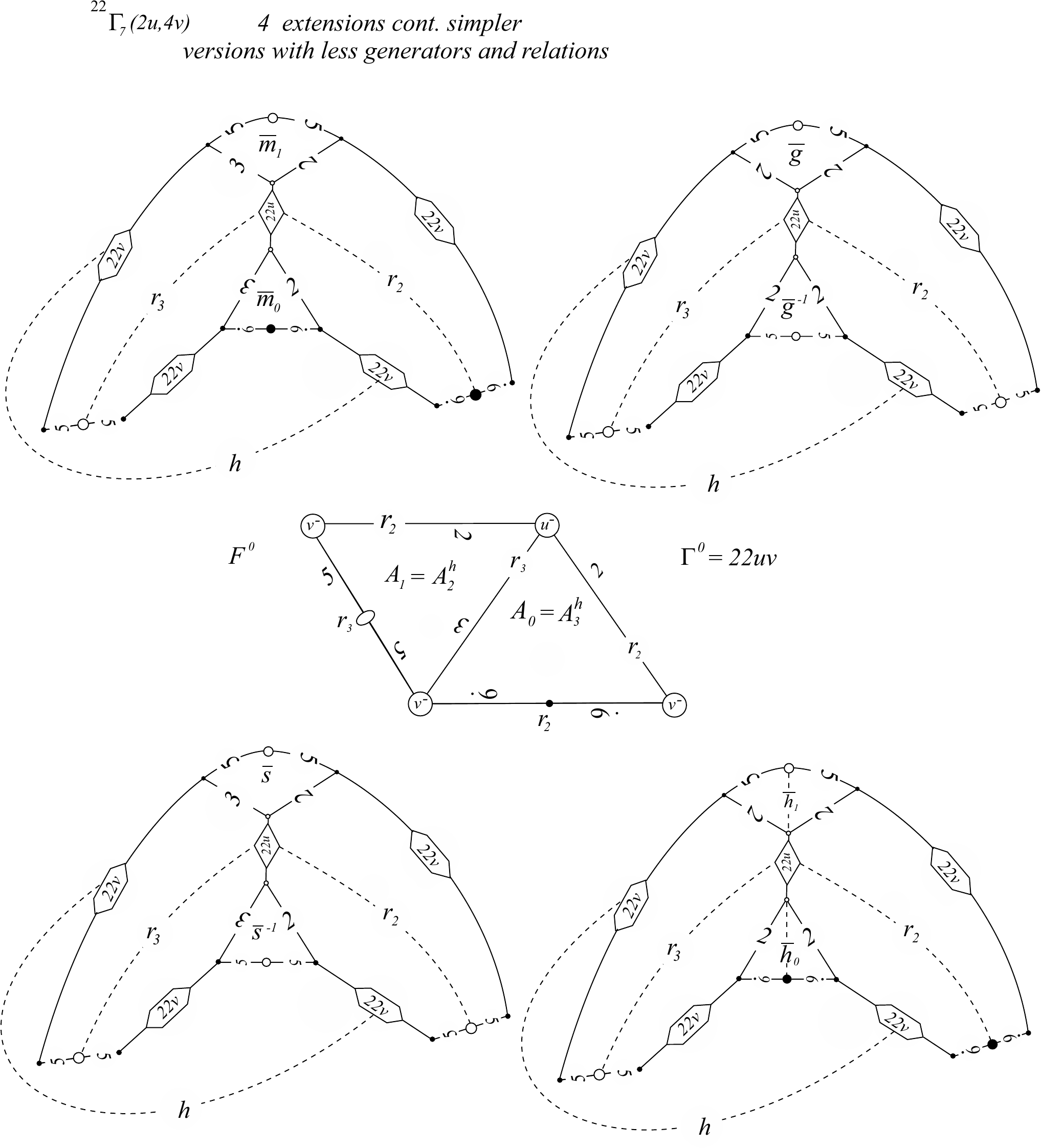} 
\caption{F2-$6_2$ to $ ^{22}\Gamma_7(2u,4v)$ simpler version}
\label{fig:F2-6cont}
\end{figure}

Look at Fig.~\ref{fig:F2-5-6}: F2-$6_1$ for these presentations with $\Gamma ^0 = 22uv$ as $\mathbb{H}^2$ group. At Fig.~\ref{fig:F2-6cont} : F2-$6_2$ you look pictures for the simpler versions. E.g. for the third extension (by a half-turn about edge 3, so 2 as well, yielding screw motion $\bar{s}$) we get in addition:

\noindent {\bf (3)} $(\bar{s} - 1  \stackrel{2}{=} (\bar{s} r_2)^2 \stackrel{3}{=} 
(\bar{s} r_3)^2 \stackrel{5}{=} \bar{s} h r_3 h \bar{s}^{-1} h r_2 h ).$

\vspace{3mm}

\noindent $\bullet \  ^{22}\Gamma _8(u, 4v) = (r, h - 1 = h^2 = r^u = (h r h r h r^{-1} h r^{-1})^v, \ 2 < u \neq 4v \in \mathbb{N}, \ \dfrac{2}{u} + \dfrac{1}{v} < 1 ) = ^2_2\Gamma _8(u, 4v)$ with 4 extensions

\noindent {\bf (1)} $(\bar{m}_0, \bar{m}_1, \bar{m}_2, \bar{m}_3  - 1 = \bar{m}_0^2 = \bar{m}_1^2 = \bar{m}_2^2 = \bar{m}_3^2 \stackrel{1}{=} \bar{m}_0 h \bar{m}_3 h \stackrel{2}{=} \bar{m}_0 r \bar{m}_0 r^{-1} \stackrel{3}{=} \bar{m}_1 h \bar{m}_2 h \stackrel{4}{=} \bar{m}_1 r \bar{m}_1 r^{-1} \stackrel{5}{=} \bar{m}_2 r^{-1} \bar{m}_3 r )$

\noindent {\bf (2)} $(\bar{s}_1, \bar{s}_2 = 1 \stackrel{1}{=} \bar{s}_1 h \bar{s}_2 h \stackrel{2}{=} \bar{s}_1 r \bar{s}_1^{-1} r^{-1} \stackrel{5}{=} (r \bar{s}_2)^2 )$

\noindent {\bf (3)} $(\bar{h}_0, \bar{h}_1, \bar{h}_2, \bar{h}_3  - 1 = \bar{h}_0^2 = \bar{h}_1^2 = \bar{h}_2^2 = \bar{h}_3^2 \stackrel{1}{=} \bar{h}_0 h \bar{h}_3 h \stackrel{2}{=} (\bar{h}_0 r)^2 \stackrel{3}{=} \bar{h}_1 h \bar{h}_2 h \stackrel{4}{=} (\bar{h}_1 r)^2  \stackrel{5}{=} \bar{h}_2 r^{-1} \bar{h}_3 r)$

\noindent {\bf (4)} $(\bar{g}_1, \bar{g}_2 = 1 \stackrel{1}{=} 
\bar{g}_1 h \bar{g}_2 h \stackrel{2}{=} \bar{g}_1 r \bar{g}_1^{-1} r 
\stackrel{5}{=} (r \bar{g}_2)^2 ).$

\begin{figure}[htbp]
	\centering
\includegraphics[width=130mm]{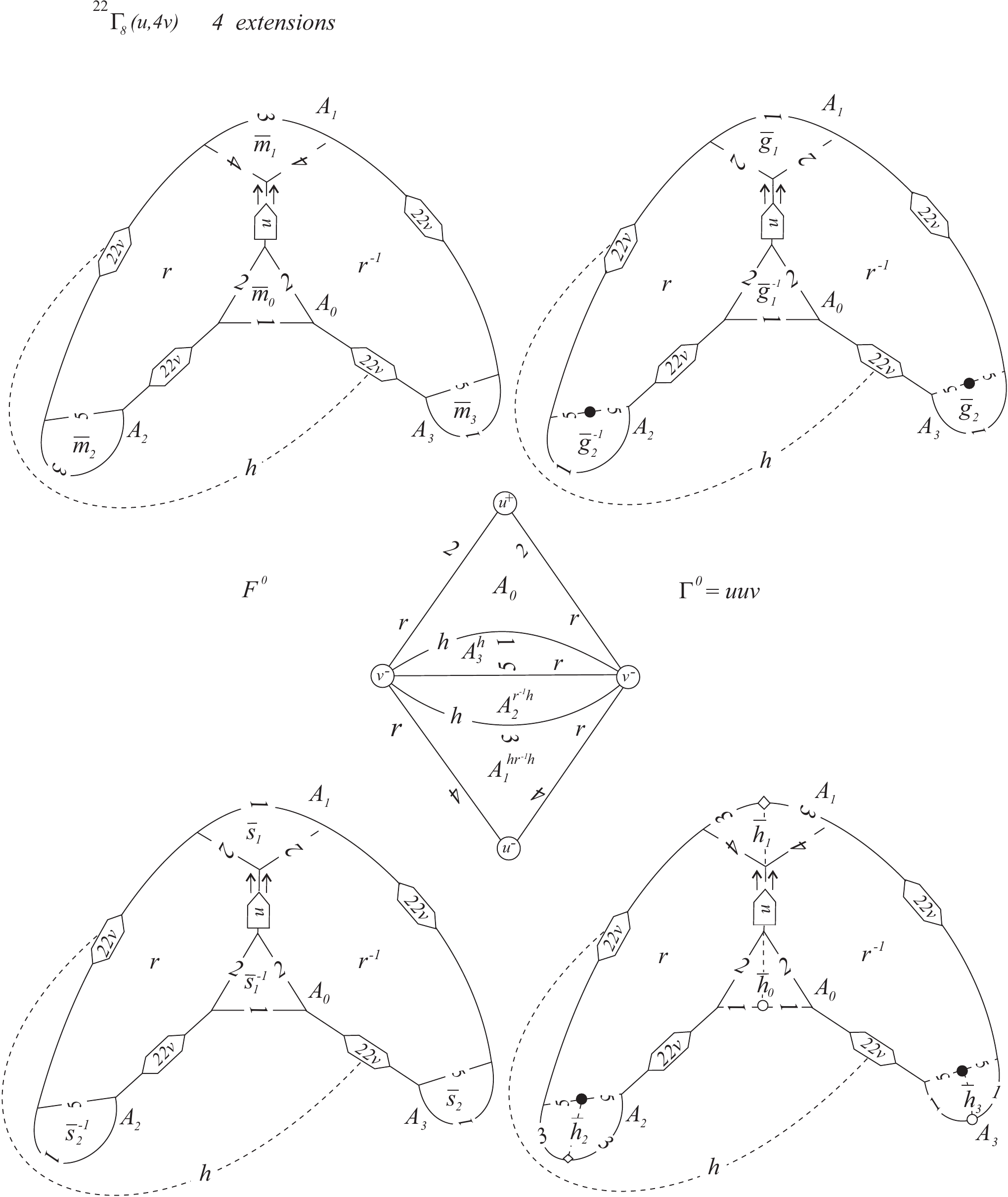}
\caption{F2-7 to $ ^{22}\Gamma_8(u,4v)$, simpler version is not pictured}
\label{fig:F2-7}
\end{figure}

Look at Fig.~\ref{fig:F2-7}: F2-7
with $\Gamma ^0 = uuv$ as $\mathbb{H}^2$ group. Simpler presentation is not pictured.

\vspace{3mm}

\noindent $\bullet \  ^{22}\Gamma _9(2u, 2v) = (s, h - 1 = h^2 = (s h s^{-1} h)^u = (s s h)^v = 1, \ 1 < u \in \mathbb{N}, 1 < v \in \mathbb{N}, \ \dfrac{1}{u} + \dfrac{2}{v} < 1 )  = ^2_2\Gamma _9(2u, 2v)$ with 2 extensions

\noindent {\bf (1)} $(\bar{m}_0, \bar{m}_1, \bar{m}_2, \bar{m}_3  - 1 = \bar{m}_0^2 = \bar{m}_1^2 = \bar{m}_2^2 = \bar{m}_3^2 \stackrel{1}{=} \bar{m}_0 h \bar{m}_1 h \stackrel{2}{=} \bar{m}_0 s^{-1} \bar{m}_3 s \stackrel{3}{=} \bar{m}_1 s \bar{m}_2 s^{-1} \stackrel{4}{=} \bar{m}_1 s^{-1} \bar{m}_2 s \stackrel{5}{=} \bar{m}_2 h \bar{m}_3 h )$

\noindent {\bf (2)} $(\bar{g}_1, \bar{g}_2 = 1 \stackrel{1}{=} \bar{g}_1 h \bar{g}_2^{-1} h \stackrel{2}{=} \bar{g}_1 s \bar{g}_2 s \stackrel{3}{=} (\bar{g}_2 s)^2 ).$

\begin{figure}[htbp]
	\centering
\includegraphics[width=130mm]{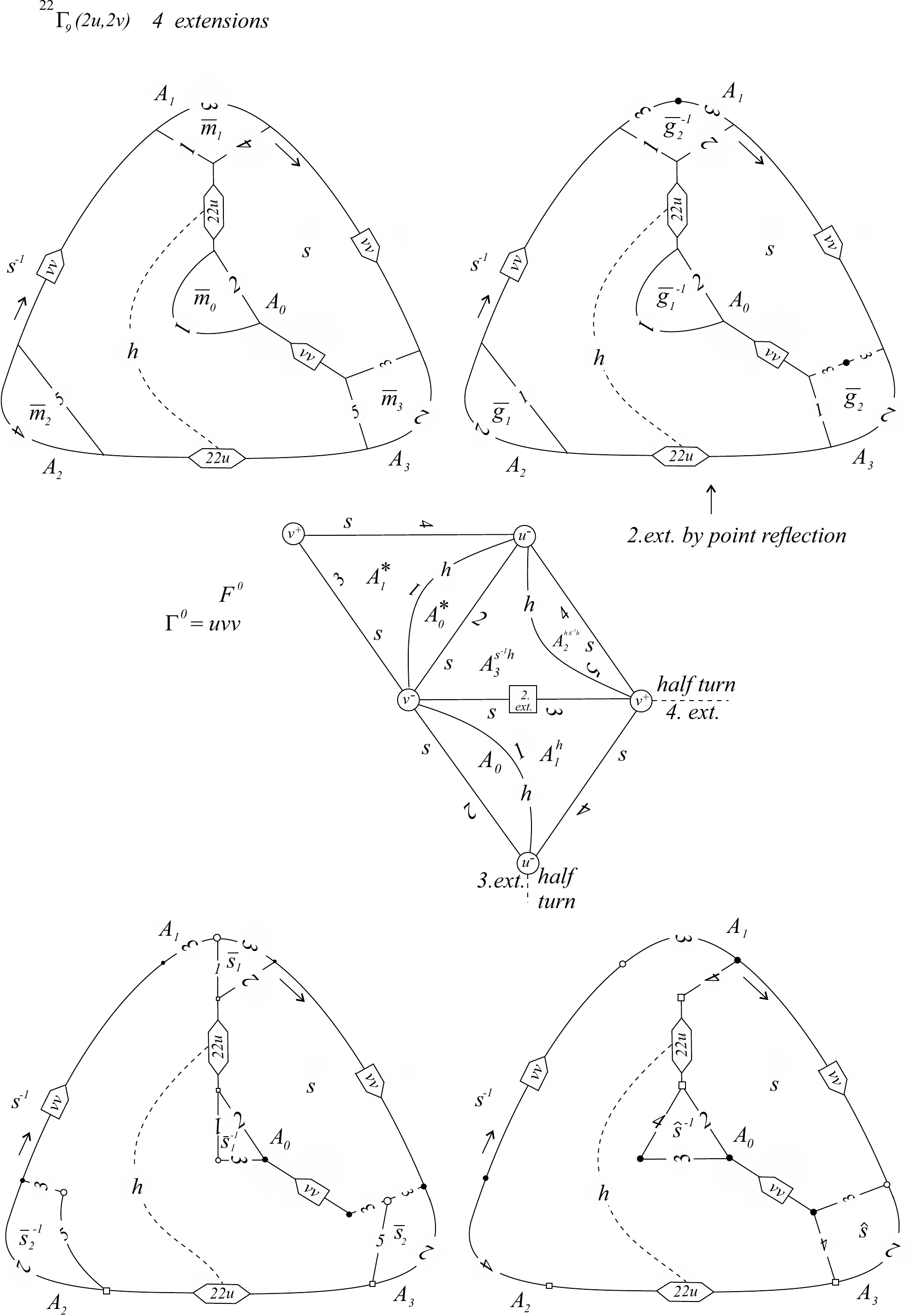}
\caption{F2-$8_1$ to $ ^{22}\Gamma_9(2u,2v)$ first version}
\label{fig:F2-8}
\end{figure}

\begin{figure}[htbp]
	\centering
\includegraphics[width=130mm]{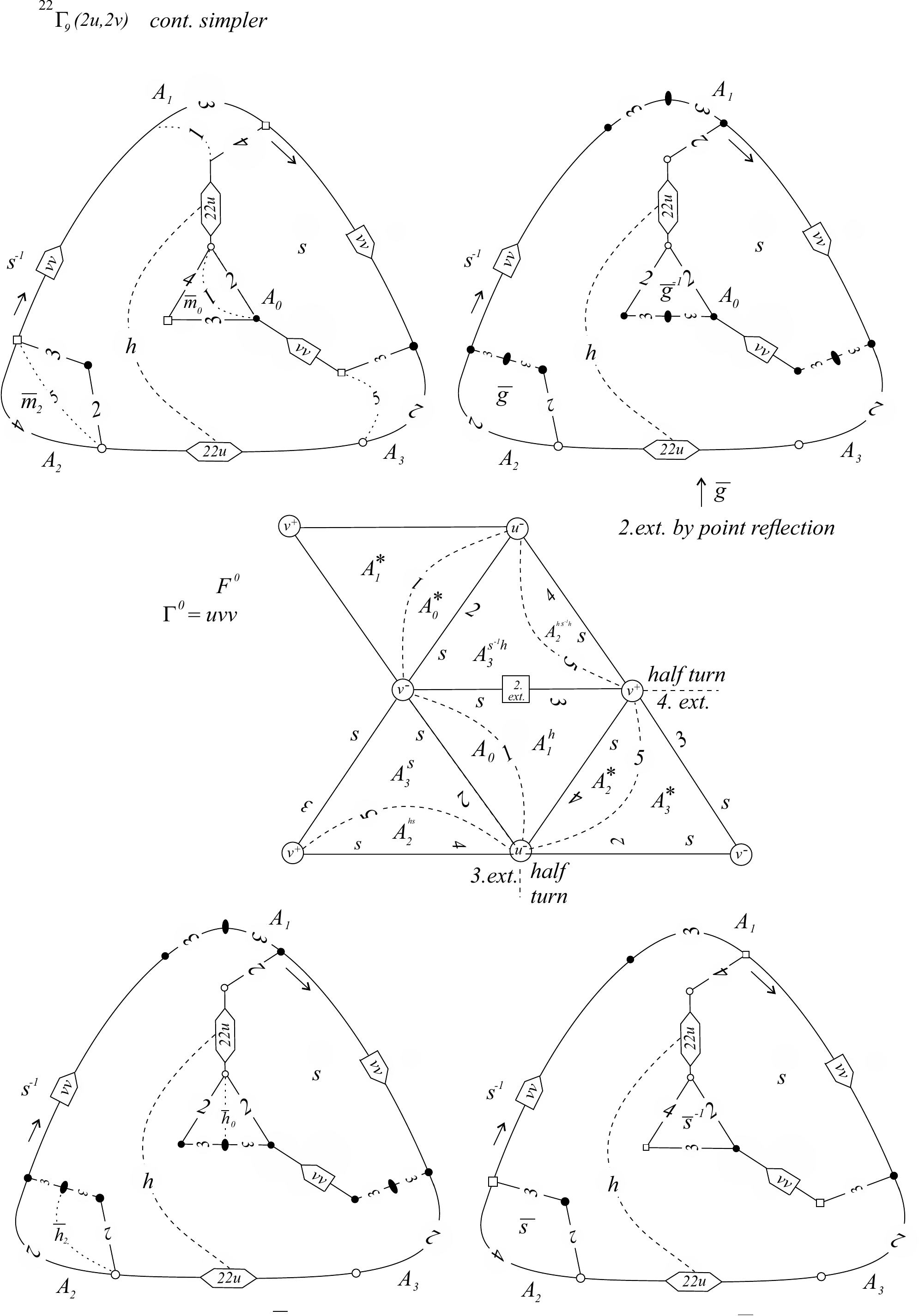}
\caption{F2-$8_2$ to $ ^{22}\Gamma_9(2u,2v)$ second simpler version}
\label{fig:F2-8cont}
\end{figure}

Look at Fig.~\ref{fig:F2-8}: F2-$8_1$ for these presentations with $\Gamma ^0 = uvv$ as $\mathbb{H}^2$ group. For the simpler second version look at Fig.~\ref{fig:F2-8cont} F2-$8_2$. For instance the second extension occurs with point reflection at edge 3, yielding glide reflection $\bar{g}$ as follows in addition:

\noindent {\bf (2)} $(\bar{g} - 1 \stackrel{2}{=} \bar{g} s h \bar{g} h s \stackrel{3}{=} (\bar{g} h s^{-1} h)^2 ) $.


\subsection{The $8$ extensions for fundamental trunc-simplex tiling and group $\Gamma_{58}(2u, 4v)$}
These top number extensions are pictured in Fig.~15 and 16 as to $F2-8+$. Maybe, the numbered edges at truncation are more convenient than in the original publication \cite{S17}. 
The stabilizer subgroup $\Gamma^0=uuvv$ of the vertex class are pictured in Fig.~15, where every vertex triangle has one $u$-point $(u^+$ or $u^-)$ 
and $v^+$-,$v^-$-point as vertices to the spatial rotation axes, respectively. Thus the $8$ possible correspondences are upper bound, indeed, that can be realized.
\begin{equation}
\Gamma_{58}(2u,4v)=(s_1,s_2~-~1=(s_1s_2^{-1})^u=
(s_1s_1s_2 s_2)^v,~u\ne2v\in\mathbb{N},~\frac{1}{u}+\frac{1}{v}<1). \notag
\end{equation}
\begin{enumerate}
\item[(\bf{1})] The trivial extension is for $\mathcal{O}^1_{58}$ as octahedron:
\begin{equation}
\begin{gathered}
(\overline{m}_0,\overline{m}_1,\overline{m}_2,\overline{m}_3~-~1=\overline{m}_0^2=\overline{m}_1^2=\overline{m}_2^2=\overline{m}_3^2 \stackrel{1}{=} 
\overline{m}_0 s_2\overline{m}_2s_2^{-1}\stackrel{2}{=}\overline{m}_0 s_1^{-1}\overline{m}_3 s_1 
\stackrel{3}{=} \\ \stackrel{3}{=}
\overline{m}_0 s_2^{-1}\overline{m}_3s_2 \stackrel{4}{=}\overline{m}_1 s_1\overline{m}_3s_1^{-1}\stackrel{5}{=}\overline{m}_1 s_1^{-1}\overline{m}_2s_1
\stackrel{6}{=}\overline{m}_1 s_2^{-1}\overline{m}_2s_2). \notag
\end{gathered}
\end{equation}
\item[(\bf{2})] 
The second extension will be derived from the two half-turns at $u^+$ and $u^-$ that fix the vertex triangles by $\overline{h}_0,\dots ,\overline{h}_3$
as follows for $\mathcal{O}_{58}^2$ in addition:
\begin{equation}
\begin{gathered}
(\overline{h}_0,\overline{h}_1,\overline{h}_2,\overline{h}_3~-~1=\overline{h}_0^2=\overline{h}_1^2=\overline{h}_2^2=\overline{h}_3^2 \stackrel{1}{=} 
\overline{h}_0 s_2\overline{h}_2 s_2^{-1} \stackrel{2}{=}\overline{h}_0 s_1^{-1}\overline{h}_3s_2  \stackrel{4}{=} \\ \stackrel{4}{=}
\overline{h}_1 s_1\overline{h}_3s_1^{-1}\stackrel{5}{=}\overline{h}_1 s_1^{-1}\overline{h}_2s_2). \notag
\end{gathered}
\end{equation}
\item[(\bf{3})] The 3rd extension will be by rotatory reflection at $v^-$ (and $v^+$ as well) rotating $A_3^{s_1} \rightarrow A_2^{s_2^{-1}}$, 
$A_0 \rightarrow A_1^{s_2^{-2}}$ and changing the half-spaces at $F^0$. This leads to glide reflection $\overline{g}_1,\overline{g}_2$ and additional presentation:
\begin{equation}
\begin{gathered}
(\overline{g}_1,\overline{g}_2~-~1\stackrel{1}{=} 
\overline{g}_1 s_1\overline{g}_2^{-1} s_2^{-1} \stackrel{2}{=}\overline{g}_1 s_2^{-1} g_2 s_1 \stackrel{3}{=} \overline{g}_1 s_1^{-1} \overline{g}_2 s_2). \notag
\end{gathered}
\end{equation}
\item[(\bf{4})] The 4th extension is derived by a half screw motions through \textcircled{5} and \textcircled{5'} $A_2^{s_2^{-1}} \rightarrow A_3^{s_1}$, 
$A_0 \rightarrow A_1^{s_1^{2}}$ 
changing the half-spaces at $F^0$. This leads to screw motions $\overline{s}_1,\overline{s}_2$ and additional presentation:
\begin{equation}
\begin{gathered}
(\overline{s}_1,\overline{s}_2~-~1\stackrel{1}{=} 
\overline{s}_1 s_1 \overline{s}_2^{-1} s_2^{-1} \stackrel{2}{=}\overline{s}_1 \overline{s}_1^{-1} \overline{s}_2 s_1 \stackrel{3}{=} 
\overline{s}_1 \overline{s}_2^{-1} \overline{s}_2 s_2)^2). \notag
\end{gathered}
\end{equation}
Here we meet with fixed-point-free actions beyond the $u-$ and $v-$ rotations whose axes close. Thus we found a double $(u,v)$ link as an interesting topological phenomenon.
(The previous and next ones are also remarkable). E.g. the case $u=3$, $v=2$ 
leads to hyperbolic metric.

The case $u = 2 = v$ leads to ideal vertices at the absolute 
(with $\mathbb{E}^2$ metric to $F^0$ domain, Fig.~15), i.e. we get $\mathbb{H}^3$ 
realization with so-called cusp phenomenon. Imagine a fundamental tetrahedron in 
$\mathbb{H}^3$ with two rectangles at opposite $u = 2$ rotational edges and $\pi/4$ 
angles at the remaining four $v = 2$ rotational edges. 
We also cite from \cite{MPS06} the $(u, 1)$ link cases with spherical $\mathbb{S}^3$ 
realization. The ``pure topological´´ $(1, 1)$ link, as manifold with ``lens 
realization´´, moreover, the $4-$truncated lens as $\mathbb{S}^3$ manifold are 
extremely remarkable, also to our context and Fig.~16.
\item[(\bf{5})]
The 5th extension is obtained by point-reflection in \textcircled{5} at edge $1$ (or in \textcircled{5'} at $4$): $A_1^{s_2^{-2}} \rightarrow A_3^{s_2}$, 
$A_0 \rightarrow A_2^{s_2^{-1}}$. These lead $\widehat{g}_1$, $\widehat{g}_2$ pairings: 
\begin{equation}
\begin{gathered}
(\widehat{g}_1,\widehat{g}_2~-~1\stackrel{1}{=} 
(\widehat{g}_1 s_2^{-1})^2 \stackrel{2}{=}\widehat{g}_1 s_1 \widehat{g}_2 s_1 \stackrel{3}{=} 
\widehat{g}_1 s_2 \widehat{g}_2 s_2 )\stackrel{4}{=} (\widehat{g}_2 s_1^{-1})^2). \notag
\end{gathered}
\end{equation}
\item[(\bf{6})]
The 6th extension is by half turn about edge $1$ (also about $4$) $A_3^{s_1} \rightarrow A_1^{s_2^{-2}}$, 
$A_0 \rightarrow A_2^{s_2^{-1}}$, leading to screw motions $\widehat{s}_1$, $\widehat{s}_2$: 
\begin{equation}
\begin{gathered}
(\widehat{s}_1,\widehat{s}_2~-~1\stackrel{1}{=} 
(\widehat{s}_1 s_2^{-1})^2 \stackrel{2}{=}\widehat{s}_1 s_2\widehat{s}_2 s_1 \stackrel{3}{=} 
\widehat{s}_1 s_1 \widehat{s}_2 s_2 )\stackrel{4}{=} (\widehat{s}_2 s_1^{-1})^2). \notag
\end{gathered}
\end{equation}
\item[(\bf{7})] The 7th extension will be $A_0 \rightarrow A_3^{*}$, $A_1^* \rightarrow A_2^{*}$ by rotatory reflection at $u^-$ and $u^+$ ($1 \leftrightarrow 4$, $2 \leftrightarrow 3$,
$5 \leftrightarrow 6$) yielding glide reflections $\widetilde{g}_1,\widetilde{g}_2$:
\begin{equation}
\begin{gathered}
(\widetilde{g}_1,\widetilde{g}_2~-~1\stackrel{1}{=} 
\widetilde{g}_1 s_1^{-1} \widetilde{g}_2 s_2^{-1} \stackrel{2}{=}\widetilde{g}_1 s_2 \widetilde{g}_1 s_1 \stackrel{5}{=} 
\widetilde{g}_2 s_1 \widetilde{g}_2 {s}_2 ). \notag
\end{gathered}
\end{equation}
\item[(\bf{8})] The 8th extension is again $A_0 \rightarrow A_3^{*}$, $A_1^* \rightarrow A_2^{*}$, but by half turns about edges $2,3,5,6$ with $1 \leftrightarrow 4$ inducing
screw motions $\widetilde{s}_1,\widetilde{s}_2$:
\begin{equation}
\begin{gathered}
(\widetilde{s}_1,\widetilde{s}_2~-~1\stackrel{1}{=} 
\widetilde{s}_1 s_1^{-1} \widetilde{s}_2 s_2^{-1} \stackrel{2}{=}(\widetilde{s}_1 s_1)^2 \stackrel{3}{=} 
(\widetilde{s}_1 s_2)^2 \stackrel{5}{=} (\widetilde{s}_2 {s}_1)^2 \stackrel{6}{=} (\widetilde{s}_2 {s}_2)^2 ). \notag
\end{gathered}
\end{equation}
\end{enumerate}

\begin{figure}[htbp]
	\centering
\includegraphics[width=130mm]{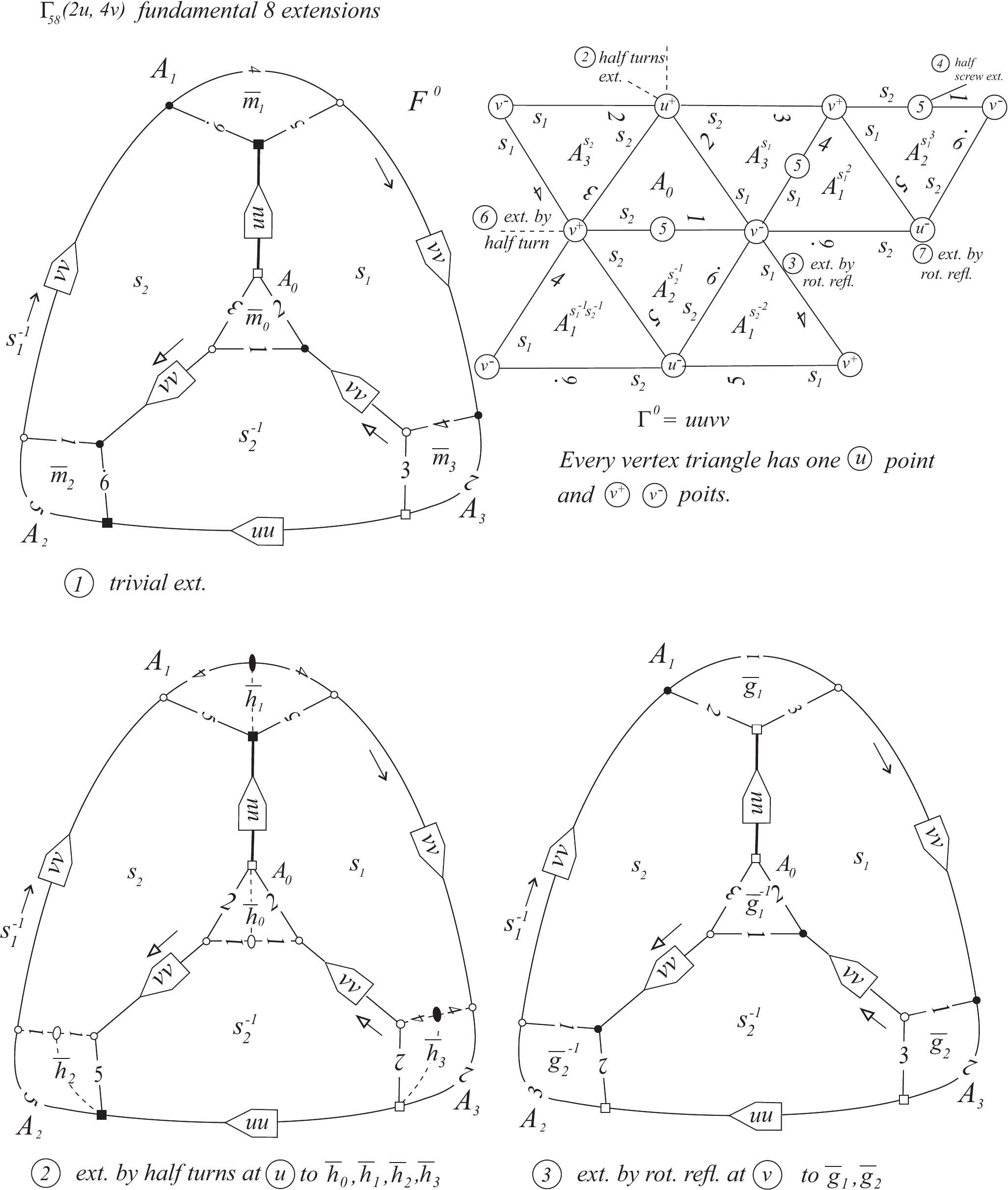}
\caption{F2-8+ for $\Gamma_{58}(2u,4v)$ the 8 fundamental extensions (repetition)}
\label{fig:F2-8+}
\end{figure}

\begin{figure}[htbp]
	\centering
\includegraphics[width=110mm]{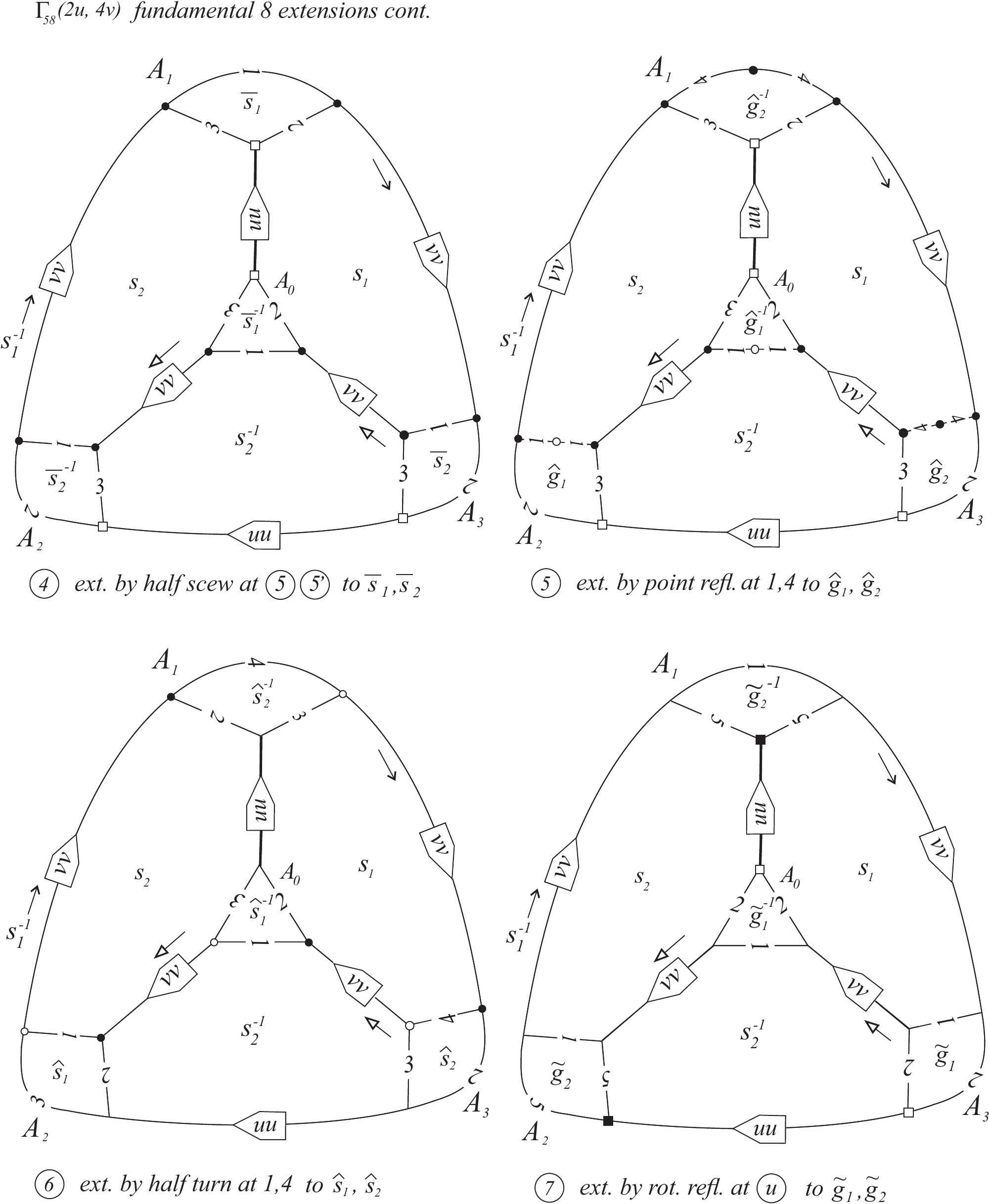} \includegraphics[width=90mm]{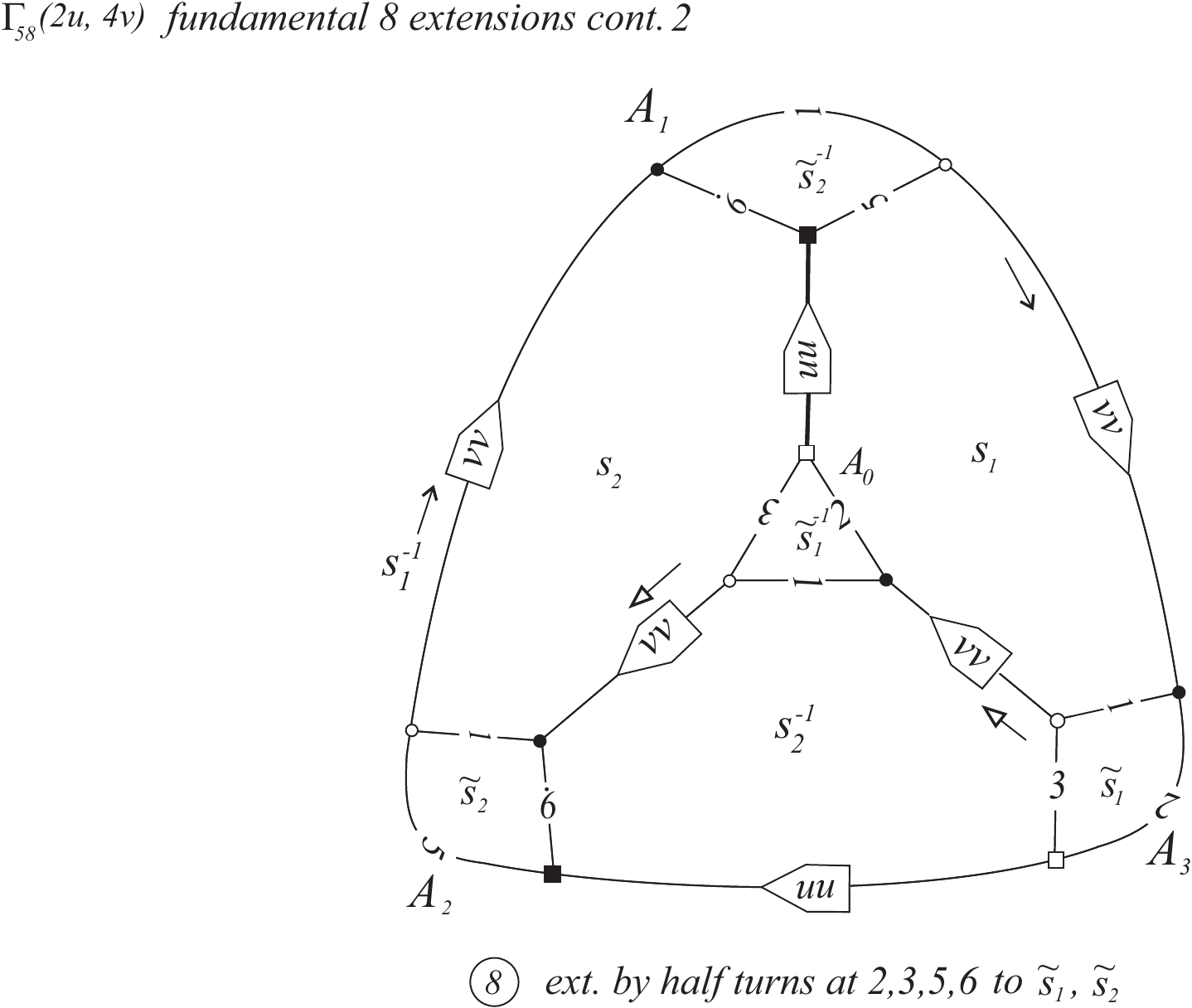}
\caption{F2-8+ for $\Gamma_{58}(2u,2v)$, important topological phenomena}
\label{fig:F2-8+cont-cont2}
\end{figure}


\section{Family F3} \label{sec:4}

\begin{figure}[htbp]
	\centering
\includegraphics[width=110mm]{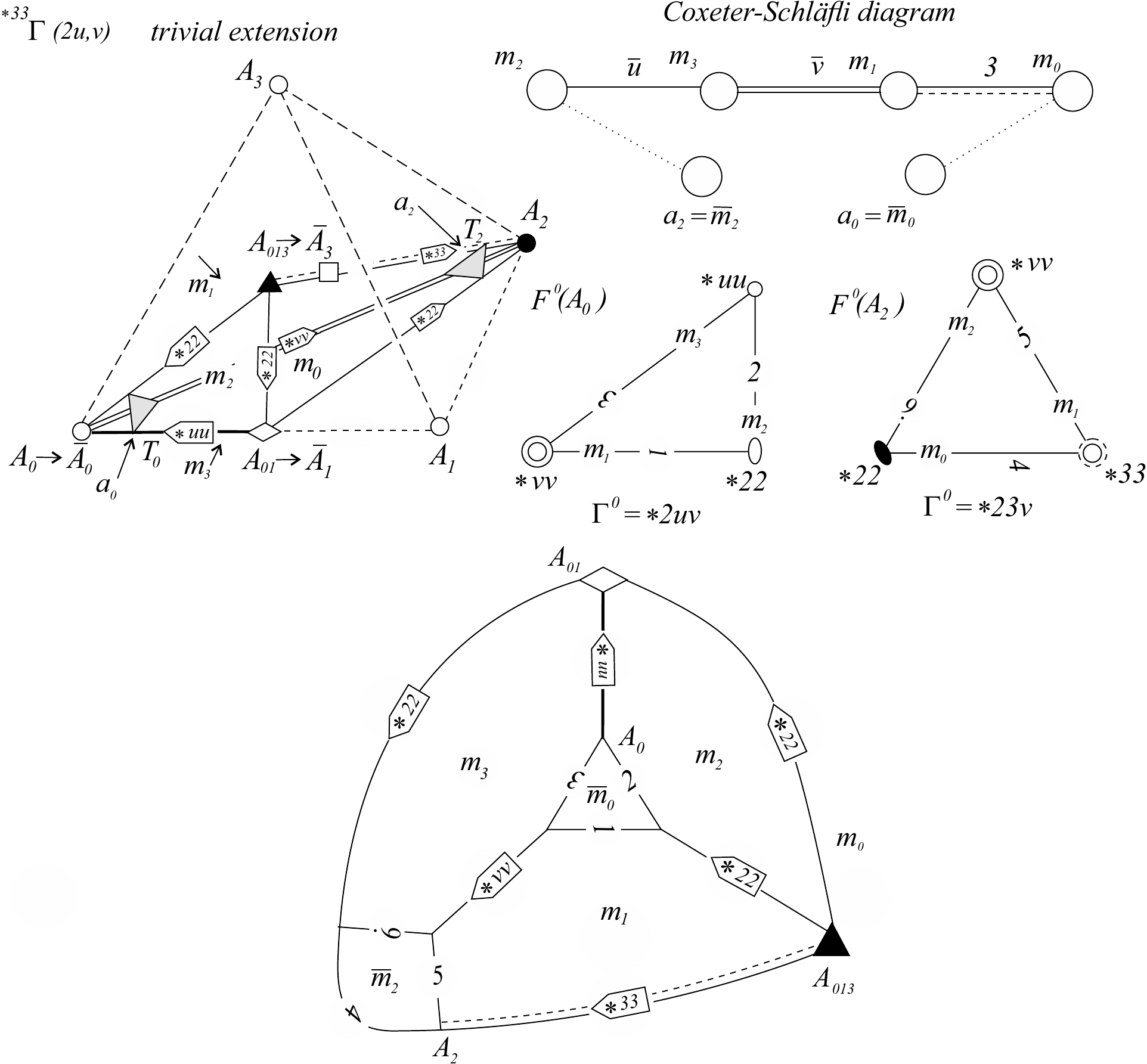} \includegraphics[width=110mm]{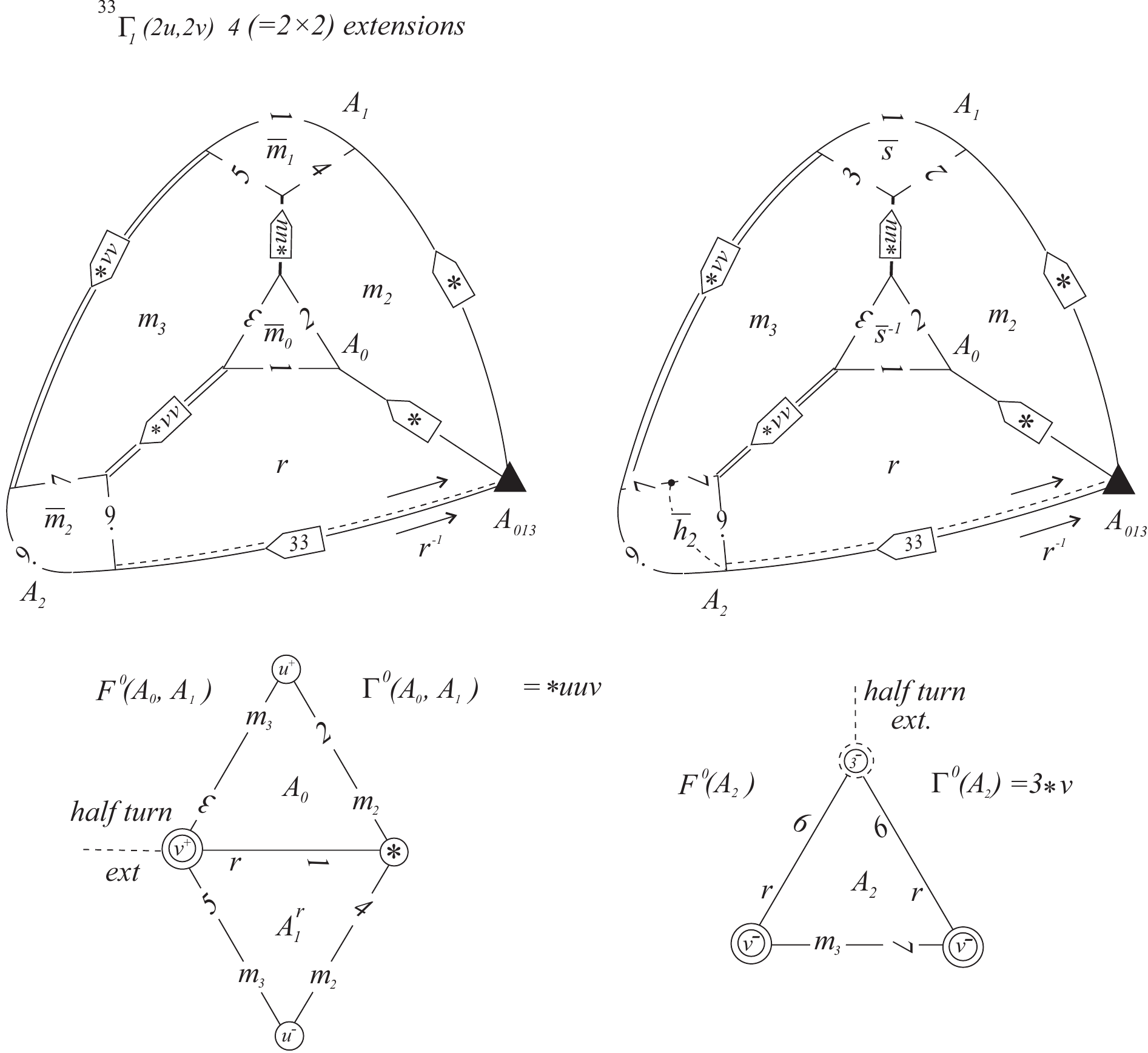}
\caption{F3-1 to $ ^{*33}\Gamma(2u,v)$ and F3-2 to $ ^{33}\Gamma_1(2u,2v)$}
\label{fig:F3-1-2}
\end{figure}

\begin{figure}[htbp]
	\centering
\includegraphics[width=110mm]{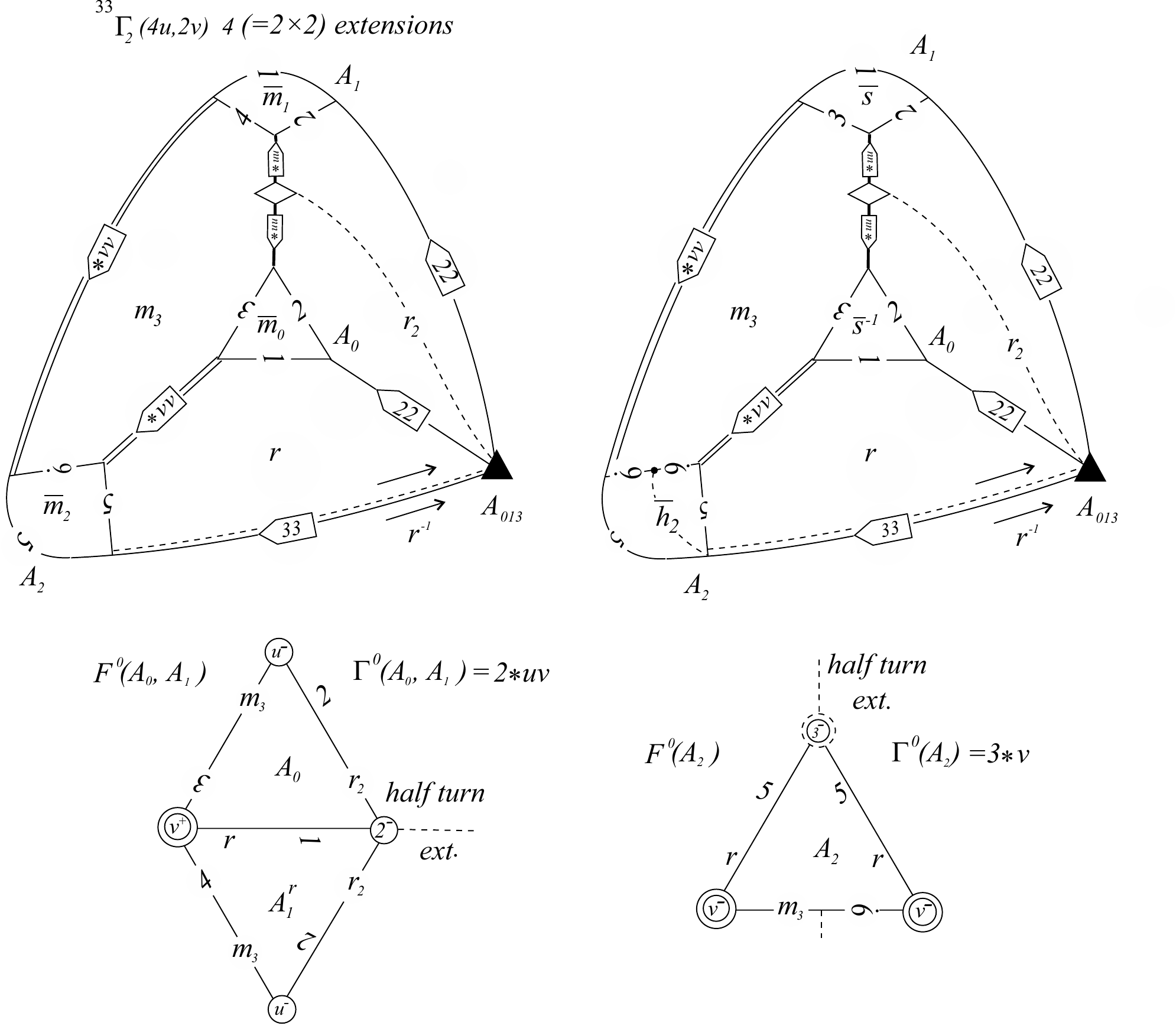} \includegraphics[width=110mm]{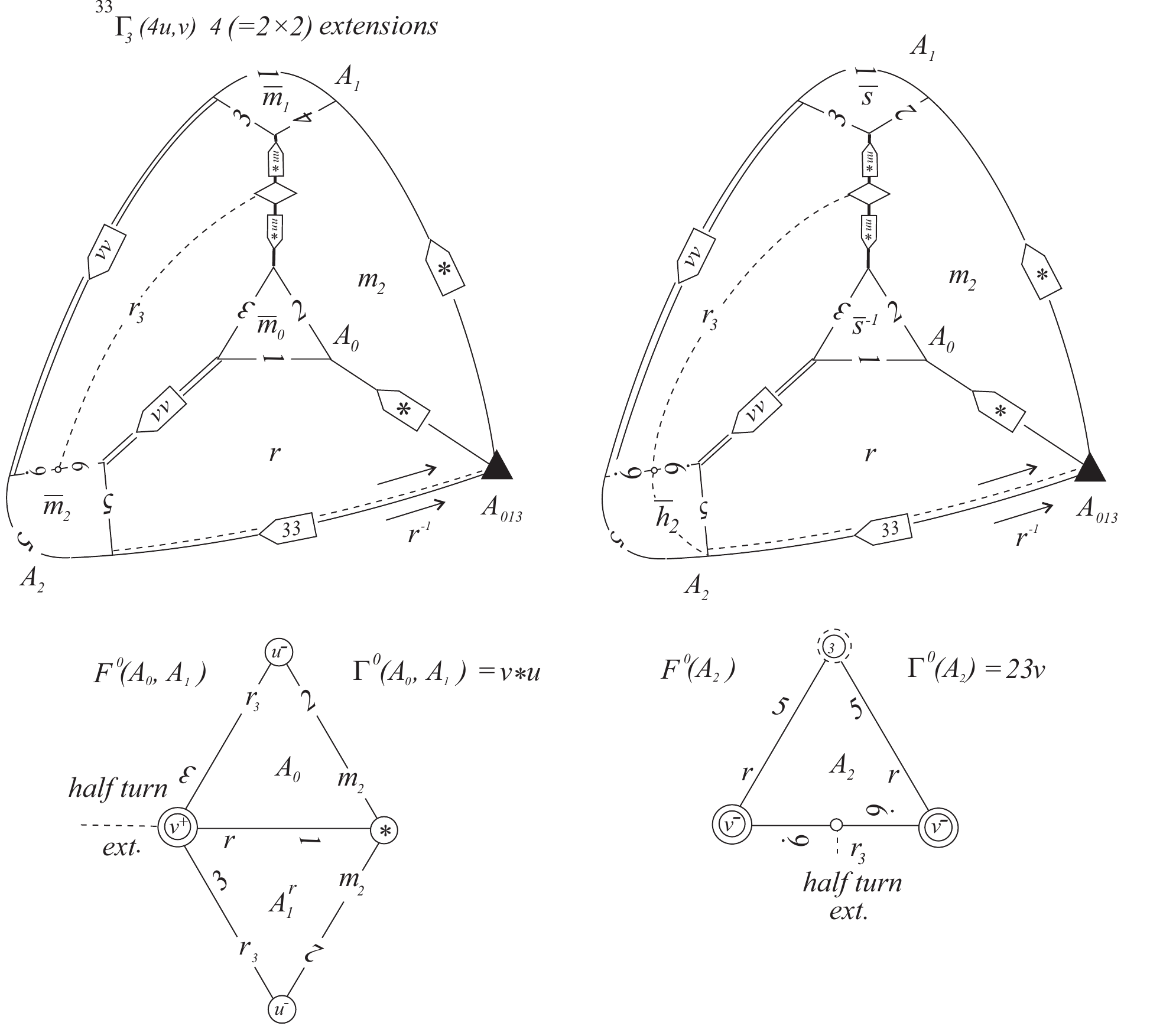}
\caption{F3-3 to $ ^{33}\Gamma_2(4u,2v)$ and F3-4 to $ ^{33}\Gamma_3(4u,v)$}
\label{fig:F3-3-4}
\end{figure}

This family is also characterized by a Coxeter reflection group 
$\{\bar{u}, \bar{v}, 3\}$ with Coxeter Schl\"{a}fli diagram 
in Fig.~\ref{fig:F3-1-2}: F3-1 with doubly truncated simplex (trunc-orthoscheme).

\noindent $\bullet \  ^{*33} \Gamma (2u, v) = ( m_0, m_1, m_2, m_3 - 1 = m_0^2 = m_1^2 = m_2^2 = m_3^2 = (m_0 m_1)^3 = \\ = (m_0 m_2)^2 =  (m_0 m_3)^2 = (m_1 m_2)^2 = (m_1 m_3)^v = (m_2 m_3)^u; \ 2u \neq v \in \mathbb{N}, \  \dfrac{1}{u} + \dfrac{1}{v} < \dfrac{1}{2}, \ \dfrac{1}{v} + \dfrac{1}{3} + \dfrac{1}{2} < 1 ) = \ ^{3m}_{\ \ 6} \Gamma (2u, v)$, with trivial extensions both for vertices $A_0$ and $A_2$

\noindent {\bf (1)} for $A_0$: $( \bar{m}_0 - 1 = \bar{m}_0^2 \stackrel{1}{=} 
(\bar{m}_0 m_1)^2 \stackrel{2}{=} (\bar{m}_0 m_2)^2 \stackrel{3}{=} (\bar{m}_0 m_3)^2 ).$

\noindent {\bf (1)} for $A_2$: $( \bar{m}_2 - 1 = \bar{m}_2^2 \stackrel{4}{=} 
(\bar{m}_2 m_0)^2 \stackrel{5}{=} (\bar{m}_2 m_1)^2 \stackrel{6}{=} (\bar{m}_2 m_3)^2 ).$

The pictures are given in Fig.~\ref{fig:F3-1-2}: F3-1
with $\Gamma ^0(A_0) = *2uv$ and $\Gamma ^0(A_2) = *23v$ as $\mathbb{H}^2$ groups. 

\vspace{3mm}

\noindent $\bullet \  ^{33} \Gamma _1 (2u, 2v) = (r, m_2, m_3 - 1 = m_2^2 = m_3^2 = r^3 = (m_2 m_3)^u = (m_3 r m_3 r^{-1})^v = m_2 r m_2 r^{-1}, \ u \neq v \in \mathbb{N}, \ \dfrac{2}{u} + \dfrac{1}{v} < 1, \ \dfrac{2}{3} + \dfrac{1}{v} < 1) = \ ^3_3 \Gamma _1 (2u, 2v)$, with 2 extensions both for vertices $A_0$, $A_1$ and $A_2$ 

\noindent {\bf (1)} for $A_0$, $A_1$: $(\bar{m}_0, \bar{m}_1 - 1 = \bar{m}_0^2 = \bar{m}_1^2 \stackrel{1}{=} \bar{m}_1 r \bar{m}_0 r^{-1} \stackrel{2}{=} (\bar{m}_0 m_2)^2 \stackrel{3}{=} (\bar{m}_0 m_3)^2 \stackrel{4}{=} (\bar{m}_1 m_2)^2 \stackrel{5}{=} (\bar{m}_1 m_3)^2 ) $

\noindent {\bf (1)} for $A_2$: $(\bar{m}_2 - 1 = \bar{m}_2^2 \stackrel{6}{=} \bar{m}_2 r \bar{m}_2 r^{-1} \stackrel{7}{=} (\bar{m}_2 m_3)^2 )$

\noindent {\bf (2)} for $A_0$, $A_1$: $(\bar{s} - 1 \stackrel{1}{=} (s r)^2 \stackrel{2}{=} \bar{s} m_2 \bar{s}^{-1} m_2 \stackrel{3}{=} \bar{s} m_3 \bar{s}^{-1} m_3 )$

\noindent {\bf (2)} for $A_2$: $(\bar{h}_2 - 1 = \bar{h}_2^2 \stackrel{6}{=} 
(\bar{h}_2 r)^2 \stackrel{7}{=} \bar{h}_2 m_3 \bar{h}_2 m_3 ).$

The pictures are given in Fig.~\ref{fig:F3-1-2}: F3-2 with $\Gamma ^0(A_0, A_1) = *uuv$ and $\Gamma ^0(A_2) = 3*v$ as $\mathbb{H}^2$ groups. 

\vspace{3mm}

\noindent $\bullet \  ^{33} \Gamma _2 (4u, 2v) = (r, r_2, m_3 - 1 = r_2^2 = m_3^2 = r^3 = (m_3 r_2 m_3 r_2)^u = (m_3 r m_3 r^{-1})^v = (r_2 r)^2, \ 2u \neq v \in \mathbb{N}, \ \dfrac{1}{u} + \dfrac{1}{v} < 1, \ \dfrac{2}{3} + \dfrac{1}{v} < 1) = \ ^3_3 \Gamma _2 (4u, 2v)$, with 2 extensions both for vertices $A_0$, $A_1$ and $A_2$

\noindent {\bf (1)} for $A_0$, $A_1$: $(\bar{m}_0, \bar{m}_1 - 1 = \bar{m}_0^2 = \bar{m}_1^2 \stackrel{1}{=} \bar{m}_1 r \bar{m}_0 r^{-1} \stackrel{2}{=} \bar{m}_0 r_2 \bar{m}_1 r_2 \stackrel{3}{=} (\bar{m}_0 m_3)^2 \stackrel{4}{=} (\bar{m}_1 m_3)^2 ) $

\noindent {\bf (1)} for $A_2$: $(\bar{m}_2 - 1 = \bar{m}_2^2 \stackrel{5}{=} \bar{m}_2 r \bar{m}_2 r^{-1} \stackrel{6}{=} (\bar{m}_2 m_3)^2 )$

\noindent {\bf (2)} for $A_0$, $A_1$: $(\bar{s} - 1 \stackrel{1}{=} (s r)^2 \stackrel{2}{=} (s r_2)^2 \stackrel{3}{=} \bar{s} m_3 \bar{s}^{-1} m_3 )$

\noindent {\bf (2)} for $A_2$: $(\bar{h}_2 - 1 = \bar{h}_2^2 \stackrel{5}{=} 
(\bar{h}_2 r)^2 \stackrel{7}{=} \bar{h}_2 m_3 \bar{h}_2 m_3 ).$

The pictures are given in Fig.~\ref{fig:F3-3-4}: F3-3 with $\Gamma ^0(A_0, A_1) = 2*uv$ and $\Gamma ^0(A_2) = 3*v$ as $\mathbb{H}^2$ groups. 
\begin{figure}[]
	\centering
\includegraphics[width=100mm]{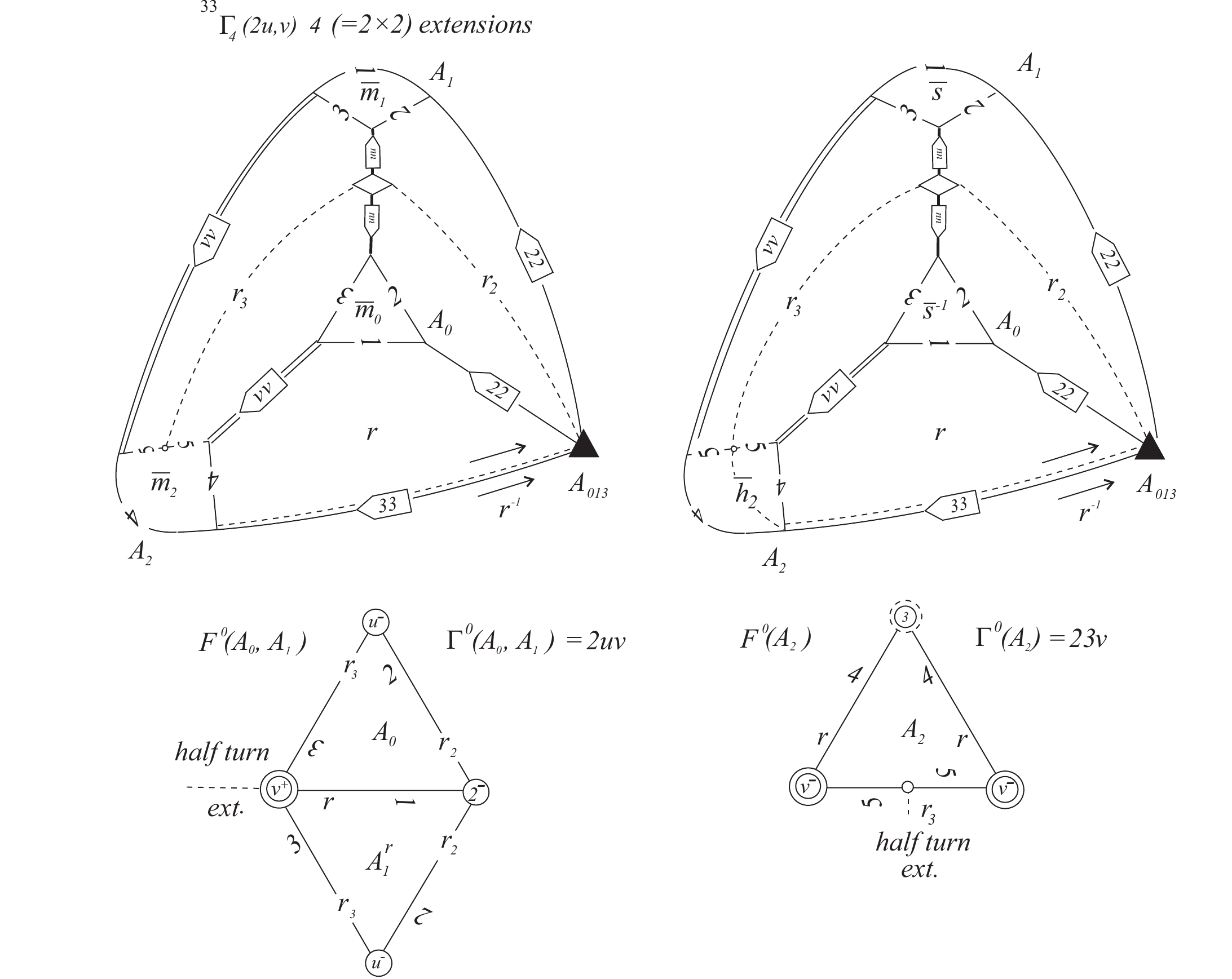} 
\caption{F3-5 to $ ^{33}\Gamma_4(2u,v)$}
\label{fig:F3-5}
\end{figure}

\noindent $\bullet \  ^{33} \Gamma _3 (4u, v) = (r, m_2, r_3 - 1 = m_2^2 = r_3^2 = r^3 = (m_2 r_3 m_2 r_3)^u = (r_3 r)^v = m_2 r m_2 r^{-1}, \ 4u \neq v \in \mathbb{N}, \ \dfrac{1}{u} + \dfrac{2}{v} < 1, \ \dfrac{2}{3} + \dfrac{2}{v} < 1)  = \ ^3_3 \Gamma _3 (4u, v)$, with 2 extensions both for vertices $A_0$, $A_1$ and $A_2$

\noindent {\bf (1)} for $A_0$, $A_1$: $(\bar{m}_0, \bar{m}_1 - 1 = \bar{m}_0^2 = \bar{m}_1^2 \stackrel{1}{=} \bar{m}_1 r \bar{m}_0 r^{-1} \stackrel{2}{=} (\bar{m}_0 m_2)^2 \stackrel{3}{=} \bar{m}_0 r_3 \bar{m}_1 r_3 \stackrel{4}{=} (\bar{m}_1 m_2)^2 ) $

\noindent {\bf (1)} for $A_2$: $(\bar{m}_2 - 1 = \bar{m}_2^2 \stackrel{5}{=} \bar{m}_2 r \bar{m}_2 r^{-1} \stackrel{6}{=} \bar{m}_2 r_3 \bar{m}_2 r_3 )$

\noindent {\bf (2)} for $A_0$, $A_1$: $(\bar{s} - 1 \stackrel{1}{=} (s r)^2 \stackrel{2}{=} \bar{s} m_2 \bar{s}^{-1} m_2 \stackrel{3}{=} (\bar{s} r_3)^2 )$

\noindent {\bf (2)} for $A_2$: $(\bar{h}_2 - 1 = \bar{h}_2^2 \stackrel{5}{=} 
(\bar{h}_2 r)^2 \stackrel{6}{=} (\bar{h}_2 r_3)^2 ).$

The pictures are given in Fig.~\ref{fig:F3-3-4}: F3-4 with $\Gamma ^0(A_0, A_1) = v*u$ and $\Gamma ^0(A_2) = 23v$ as $\mathbb{H}^2$ groups. 

\vspace{3mm}

\noindent $\bullet \  ^{33} \Gamma _4 (2u, v) = (r, r_2, r_3 - 1 = r_2^2 = r_3^2 = r^3 = (r_2 r_3)^u = (r_3 r)^v = (r r_2)^2, \ 2u \neq v \in \mathbb{N}, \ \dfrac{2}{u} + \dfrac{2}{v} < 1, \ \dfrac{2}{3} + \dfrac{2}{v} < 1) = \ ^3_3 \Gamma _4 (2u, v)$, with 2 extensions both for vertices $A_0$, $A_1$ and $A_2$

\noindent {\bf (1)} for $A_0$, $A_1$: $(\bar{m}_0, \bar{m}_1 - 1 = \bar{m}_0^2 = \bar{m}_1^2 \stackrel{1}{=} \bar{m}_1 r \bar{m}_0 r^{-1} \stackrel{2}{=} \bar{m}_0 r_2 \bar{m}_1 r_2 \stackrel{3}{=} \bar{m}_0 r_3 \bar{m}_1 r_3 ) $

\noindent {\bf (1)} for $A_2$: $(\bar{m}_2 - 1 = \bar{m}_2^2 \stackrel{4}{=} \bar{m}_2 r \bar{m}_2 r^{-1} \stackrel{5}{=} \bar{m}_2 r_3 \bar{m}_2 r_3 )$

\noindent {\bf (2)} for $A_0$, $A_1$: $(\bar{s} - 1 \stackrel{1}{=} (s r)^2 \stackrel{2}{=} (\bar{s} r_2)^2 \stackrel{3}{=} (\bar{s} r_3)^2 )$

\noindent {\bf (2)} for $A_2$: $(\bar{h}_2 - 1 = \bar{h}_2^2 \stackrel{4}{=} 
(\bar{h}_2 r)^2 \stackrel{5}{=} (\bar{h}_2 r_3)^2 ).$

The pictures are given in Fig.~\ref{fig:F3-5}: F3-5 with $\Gamma ^0(A_0, A_1) = 2uv$ and $\Gamma ^0(A_2) = 23v$ as $\mathbb{H}^2$ groups.


\section{Family F4} \label{sec:5}

The family F4 is characterized by

\noindent $^{*22} \Gamma_1 (u, 2v, w) = ( m_0, m_1, m_2, m_3 - 1 = m_0^2 = m_1^2 = m_2^2 = m_3^2 = (m_0 m_1)^w = \\ = (m_0 m_2)^2 =  (m_0 m_3)^2 = (m_1 m_2)^2 = (m_1 m_3)^v = (m_2 m_3)^u; \ u \neq w, v \in \mathbb{N},  \ \dfrac{1}{u} + \dfrac{1}{v} < \dfrac{1}{2}, \ \dfrac{1}{v} + \dfrac{1}{w} < \dfrac{1}{2} ) = \ ^{mm2}_{\ \ \ \ 4} \Gamma_1 (u, 2v, w)$, with trivial extensions both for vertices $A_0$ and $A_2$

\noindent {\bf (1)} for $A_0$: $( \bar{m}_0 - 1 = \bar{m}_0^2 \stackrel{1}{=} 
(\bar{m}_0 m_1)^2 \stackrel{2}{=} (\bar{m}_0 m_2)^2 \stackrel{3}{=} (\bar{m}_0 m_3)^2 ).$

\noindent {\bf (1)} for $A_2$: $( \bar{m}_2 - 1 = \bar{m}_2^2 \stackrel{4}{=} 
(\bar{m}_2 m_0)^2 \stackrel{5}{=} (\bar{m}_2 m_1)^2 \stackrel{6}{=} (\bar{m}_2 m_3)^2 ).$

\begin{figure}[htbp]
	\centering
\includegraphics[width=110mm]{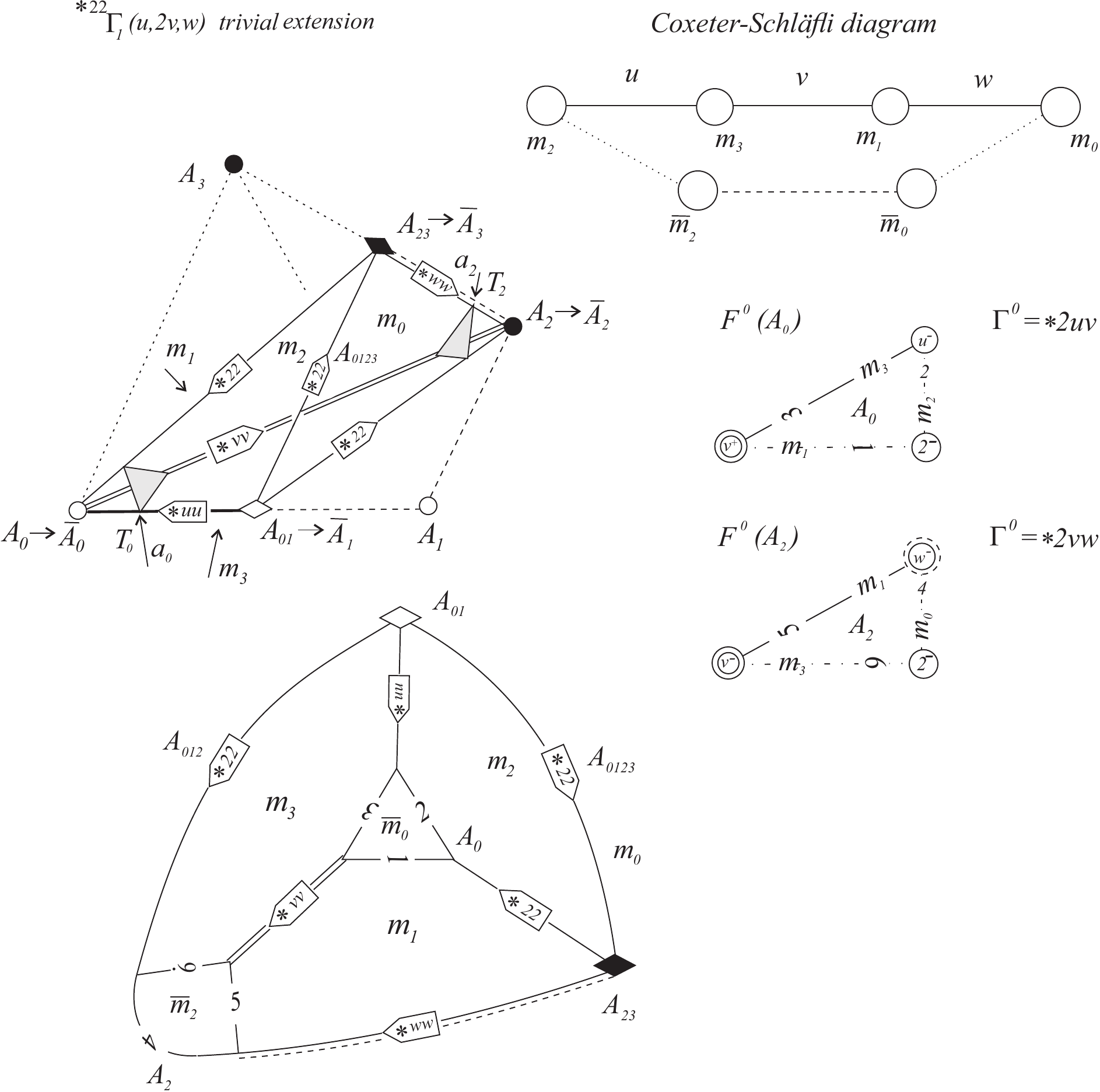} \includegraphics[width=110mm]{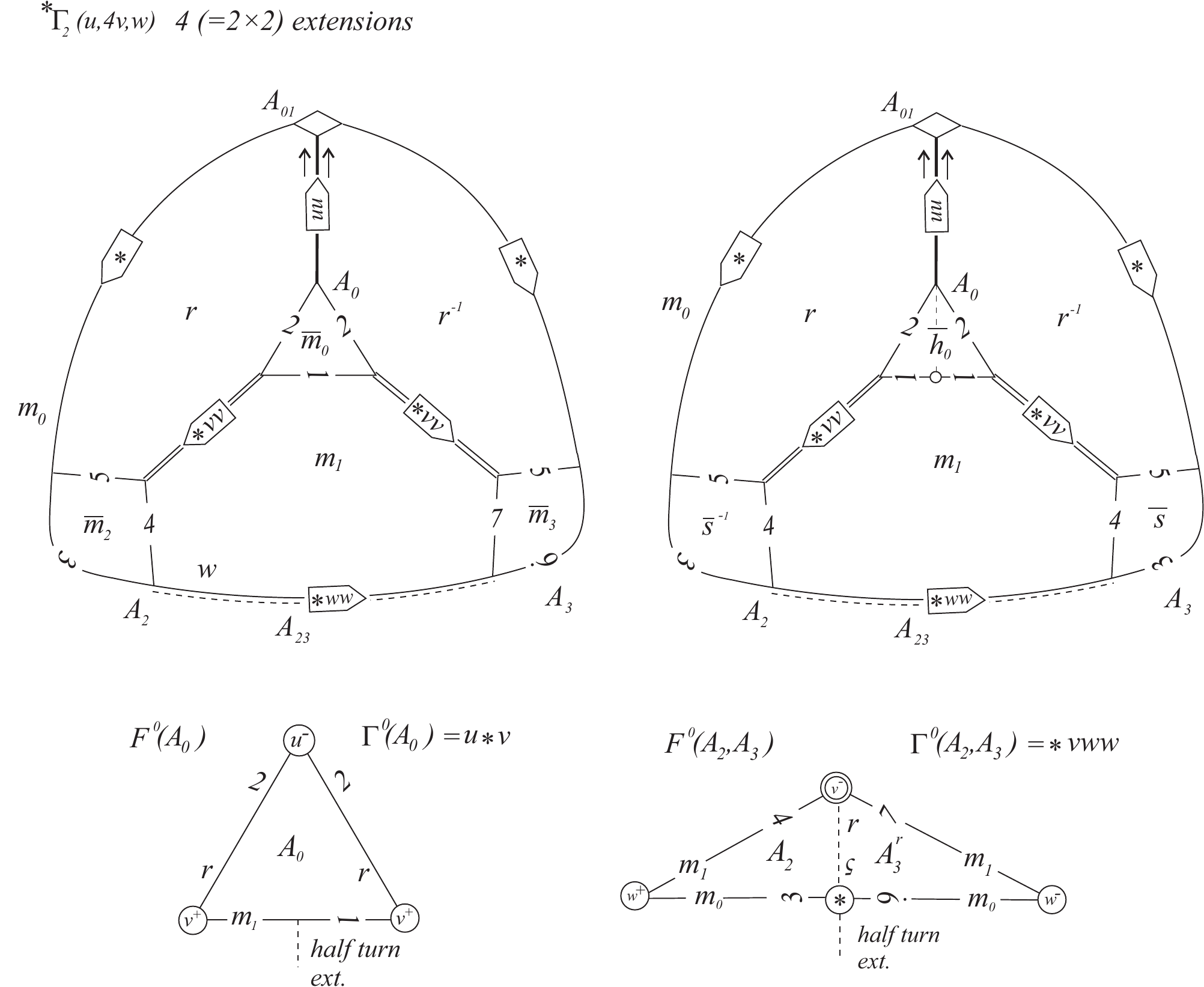} 
\caption{F4-1 to $ ^{*22}\Gamma_1(u,2v,w)$ and F4-2 to $ ^{*}\Gamma_2(u,4v,w)$}
\label{fig:F4-1-2}
\end{figure}

The pictures are given in Fig.~\ref{fig:F4-1-2}: F4-1 with 
$\Gamma ^0(A_0) = *2uv$ and $\Gamma ^0(A_2) = *2vw$ as $\mathbb{H}^2$ groups. 

{\bf Other non-fundamental cases from F4}

\noindent $\bullet \  ^* \Gamma _2 (u, 4v, w) = (r, m_0, m_1 - 1 = m_0^2 = m_1^2 = r^u = (m_1 r m_1 r^{-1})^v = (m_0 m_1)^w = m_0 r m_0 r^{-1}, \ u \neq w, v \in \mathbb{N},  \ \dfrac{2}{u} + \dfrac{1}{v} < 1, \ \dfrac{2}{w} + \dfrac{1}{v} < 1)  = \ ^m_2 \Gamma _2 (u, 4v, w)$, with 2 extensions both for vertices $A_0$ and $A_2$, $A_3$

\noindent {\bf (1)} for $A_0$: $(\bar{m}_0 - 1 = \bar{m}_0^2 \stackrel{1}{=} (\bar{m}_0 m_1)^2 \stackrel{2}{=} \bar{m}_0 r \bar{m}_0 r^{-1} ) $

\noindent {\bf (1)} for $A_2$, $A_3$: $(\bar{m}_2, \bar{m}_3 - 1 = \bar{m}_2^2 = \bar{m}_3^2 \stackrel{3}{=} (\bar{m}_2 m_0)^2 \stackrel{4}{=} (\bar{m}_2 m_1)^2 \stackrel{5}{=} \bar{m}_3 r \bar{m}_2 r^{-1} \stackrel{6}{=} (\bar{m}_3 m_0)^2 \stackrel{7}{=} (\bar{m}_3 m_1)^2   )$

\noindent {\bf (2)} for $A_0$: $(\bar{h}_0 - 1 = \bar{h}_0^2 \stackrel{1}{=} (m_1 \bar{h}_0)^2 \stackrel{2}{=} (\bar{h}_0 r)^2 )$

\noindent {\bf (2)} for $A_2$, $A_3$: $(\bar{s} - 1 \stackrel{3}{=} 
m_0 \bar{s} m_0 \bar{s}^{-1} \stackrel{4}{=} m_1 \bar{s} m_1 \bar{s}^{-1} 
\stackrel{5}{=} (r \bar{s})^2 ).$

The pictures are given in Fig.~\ref{fig:F4-1-2}: F4-2 with $\Gamma ^0(A_0) = u*v$ and $\Gamma ^0(A_2, A_3) = *vww$ as $\mathbb{H}^2$ groups. 

\vspace{3mm}

\noindent $\bullet \  ^* \Gamma _3 (2u, 4v, 2w) = (r_1, m_0, m_2, m_3 - 1 = r_1^2 = m_0^2 = m_2^2 = m_3^2 = (m_0 m_2)^2 = (m_0 m_3)^2 = (m_2 m_3)^u = (m_2 r_1 m_3 r_1)^v = (m_0 r_1 m_0 r_1)^w, \ u \neq w, v  \in \mathbb{N},  \ \dfrac{1}{u} + \dfrac{1}{v} < 1, \ \dfrac{1}{v} + \dfrac{1}{w} < 1) = \ ^m_2 \Gamma _3 (2u, 4v, 2w)$, with 2 extensions both for vertices $A_0$ and $A_2$, $A_3$

\noindent {\bf (1)} for $A_0$: $(\bar{m}_0 - 1 = \bar{m}_0^2 \stackrel{1}{=} (\bar{m}_0 r_1)^2  \stackrel{2}{=} (\bar{m}_0 m_2)^2 \stackrel{3}{=} (\bar{m}_0 m_3)^2 ) $

\noindent {\bf (1)} for $A_2$, $A_3$: $(\bar{m}_2, \bar{m}_3 - 1 = \bar{m}_2^2 = \bar{m}_3^2 \stackrel{4}{=} (\bar{m}_2 m_0)^2 \stackrel{5}{=} \bar{m}_2 r_1 \bar{m}_3 r_1 \stackrel{6}{=}
(\bar{m}_2 m_3)^2 \stackrel{7}{=} (\bar{m}_3 m_2)^2 \stackrel{8}{=} (\bar{m}_3 m_0)^2 )$

\noindent {\bf (2)} for $A_0$: $(\bar{h}_0 - 1 = \bar{h}_0^2 \stackrel{1}{=} (r_1 \bar{h}_0)^2 \stackrel{2}{=} m_2 \bar{h}_0 m_3 \bar{h}_0 )$

\noindent {\bf (2)} for $A_2$, $A_3$: $(\bar{s} - 1 \stackrel{4}{=} m_0 \bar{s} 
m_0 \bar{s}^{-1} \stackrel{5}{=} (r_1 \bar{s})^2 \stackrel{6}{=} m_3 \bar{s} 
m_2 \bar{s}^{-1} ).$

\begin{figure}[htbp]
	\centering
\includegraphics[width=110mm]{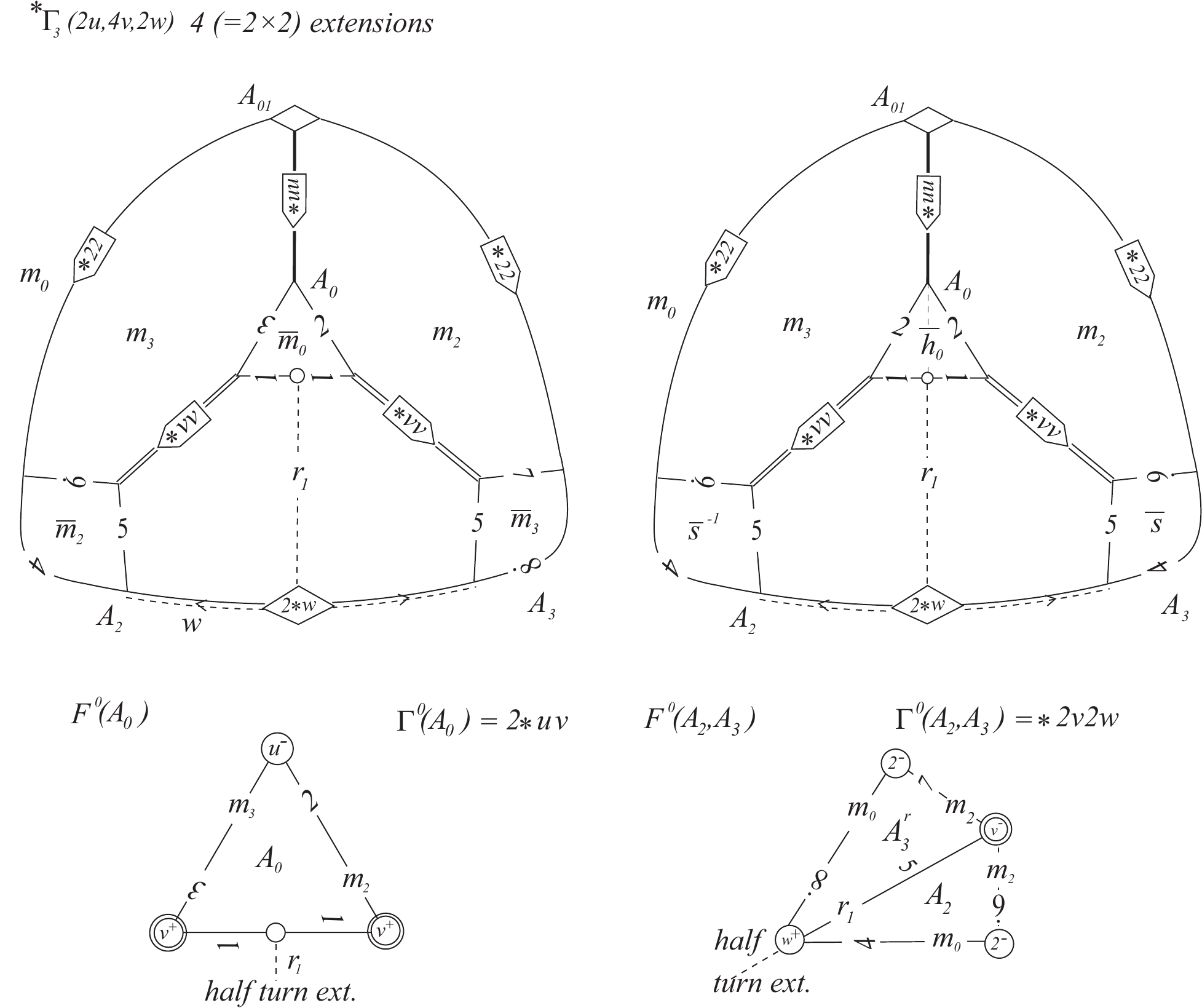} \includegraphics[width=110mm]{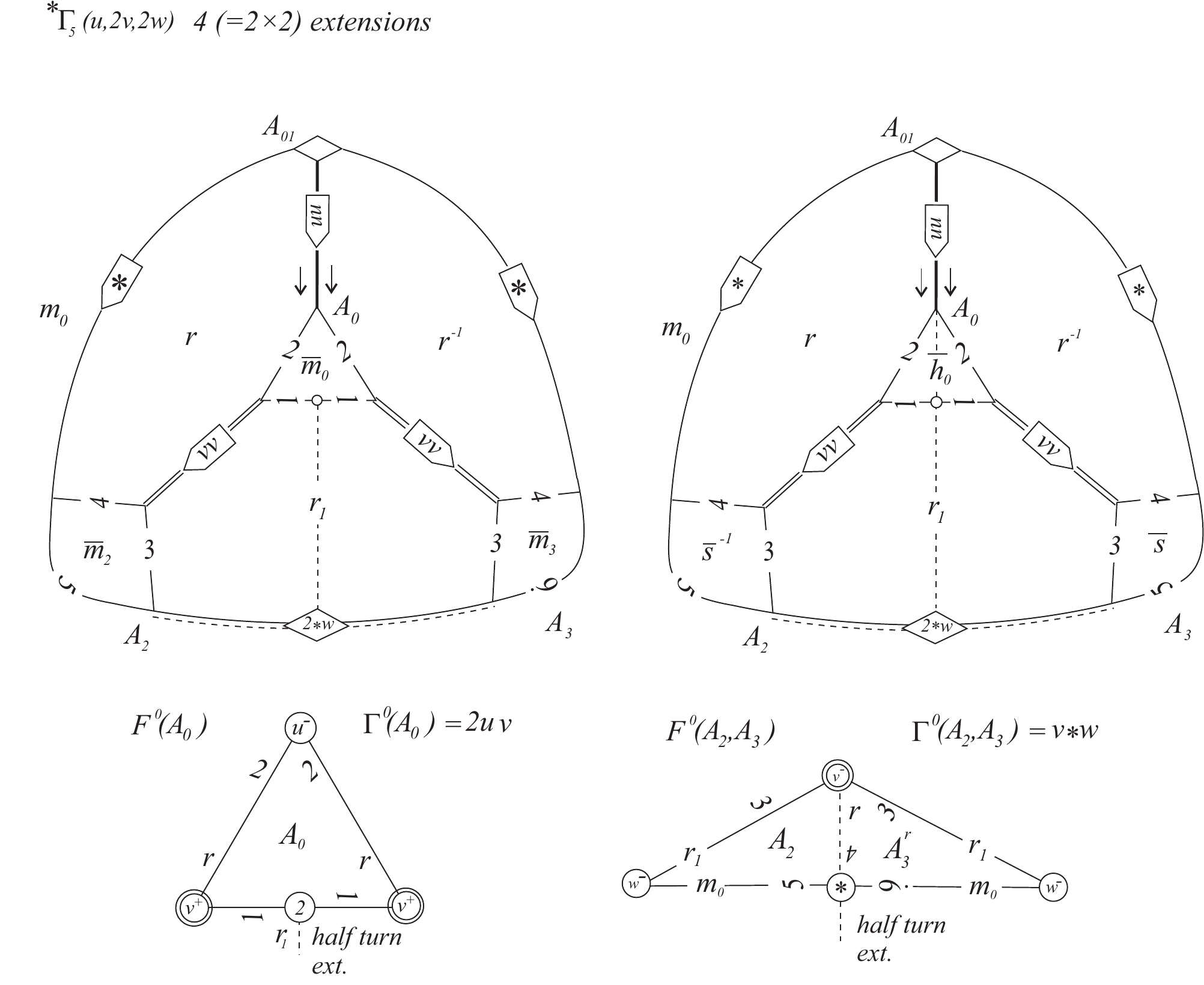}
\caption{F4-3 to $ ^{*}\Gamma_3(2u,4v,2w)$ and F4-4 to $ ^{*}\Gamma_5(u,2v,2w)$}
\label{fig:F4-3-4}
\end{figure}

The pictures are given in Fig.~\ref{fig:F4-3-4}: F4-3 with $\Gamma ^0(A_0) = 2*uv$ and $\Gamma ^0(A_2, A_3) = *2v2w$ as $\mathbb{H}^2$ groups. 

\vspace{3mm}

\noindent $\bullet \  ^* \Gamma _5 (u, 2v, 2w) = (m_0, r_1, r - 1 = m_0^2 = r_1^2 = r^u = (r_1 r)^v = (m_0 r_1 m_0 r_1)^w = m_0 r m_0 r^{-1}, \ u \neq 2w, v   \in \mathbb{N},  \ \dfrac{1}{u} + \dfrac{1}{v} < \dfrac{1}{2}, \ \dfrac{2}{v} + \dfrac{1}{w} < 1) = \ ^m_2 \Gamma _5 (u, 2v, 2w)$, with 2 extensions both for vertices $A_0$ and $A_2$, $A_3$ 

\noindent {\bf (1)} for $A_0$: $(\bar{m}_0 - 1 = \bar{m}_0^2 \stackrel{1}{=} (\bar{m}_0 r_1)^2 \stackrel{2}{=} \bar{m}_0 r \bar{m}_0 r^{-1} ) $

\noindent {\bf (1)} for $A_2$, $A_3$: $(\bar{m}_2, \bar{m}_3 - 1 = \bar{m}_2^2 = \bar{m}_3^2 \stackrel{3}{=} \bar{m}_2 r_1 \bar{m}_3 r_1  \stackrel{4}{=} \bar{m}_3 r \bar{m}_2 r^{-1} \stackrel{5}{=}  (\bar{m}_2 m_0)^2  \stackrel{6}{=} (\bar{m}_3 m_0)^2 )$

\noindent {\bf (2)} for $A_0$: $(\bar{h}_0 - 1 = \bar{h}_0^2 \stackrel{1}{=} (\bar{h}_0 r_1)^2 \stackrel{2}{=} (\bar{h}_0 r)^2 )$

\noindent {\bf (2)} for $A_2$, $A_3$: $(\bar{s} - 1 \stackrel{3}{=} (r_1 \bar{s})^2 
\stackrel{4}{=} (r \bar{s})^2 \stackrel{5}{=} m_0 \bar{s} m_0 \bar{s}^{-1}).$

The pictures are given in Fig.~\ref{fig:F4-3-4}: F4-4 with $\Gamma ^0(A_0) = 2uv$ and $\Gamma ^0(A_2, A_3) = v*w$ as $\mathbb{H}^2$ groups. 


\subsection{Truncated half-simplices from F4}


\noindent $\bullet \  ^{22} \Gamma _3 (u, 4v, 2w) = (m_1, r_3, h - 1 = m_1^2 = r_3^2 = h^2 = (r_3 h)^u = (m_1 r_3 h m_1 h r_3)^v = (m_1 h m_1 h)^w, \ u \neq 2w, v  \in \mathbb{N},  \ \dfrac{2}{u} + \dfrac{1}{v} < 1, \ \dfrac{1}{v} + \dfrac{1}{w} < 1) = ^2_2 \Gamma _3 (u, 4v, 2w)$, with 2 extensions both for vertices $A_0$, $A_1$ and $A_2$, $A_3$

\noindent {\bf (1)} for $A_0$, $A_1$: $(\bar{m}_0, \bar{m}_1 - 1 = \bar{m}_0^2 = \bar{m}_1^2 \stackrel{1}{=} (m_1 \bar{m}_0)^2 \stackrel{2}{=} \bar{m}_0 h \bar{m}_1 h \stackrel{3}{=} \bar{m}_0 r_3 \bar{m}_1 r_3 )$

\noindent {\bf (1)} for $A_2$, $A_3$: $(\bar{m}_2, \bar{m}_3 - 1 =  \bar{m}_2^2 = \bar{m}_3^2 \stackrel{4}{=} (\bar{m}_2 r_3)^2 \stackrel{5}{=} \bar{m}_2 h \bar{m}_3 h \stackrel{6}{=} (\bar{m}_2 m_1)^2 \stackrel{7}{=} (\bar{m}_3 m_1)^2 )$

\noindent {\bf (2)} for $A_0$, $A_1$ by half-turn $\bar{h}_0$: $(\bar{h}_0 - 1 = \bar{h}_0^2 \stackrel{1}{=} (m_1 \bar{h}_0)^2 \stackrel{3}{=} (\bar{h}_0 h r_3)^2 )$

\noindent {\bf (2)} for $A_2$, $A_3$ by half-turn $\bar{h}_2$: $(\bar{h}_2 - 1 
\stackrel{4}{=} (\bar{h}_2 r_3)^2 \stackrel{6}{=} m_1 \bar{h}_2 h m_1 h \bar{h}_2).$

\begin{figure}[htbp]
	\centering
\includegraphics[width=130mm]{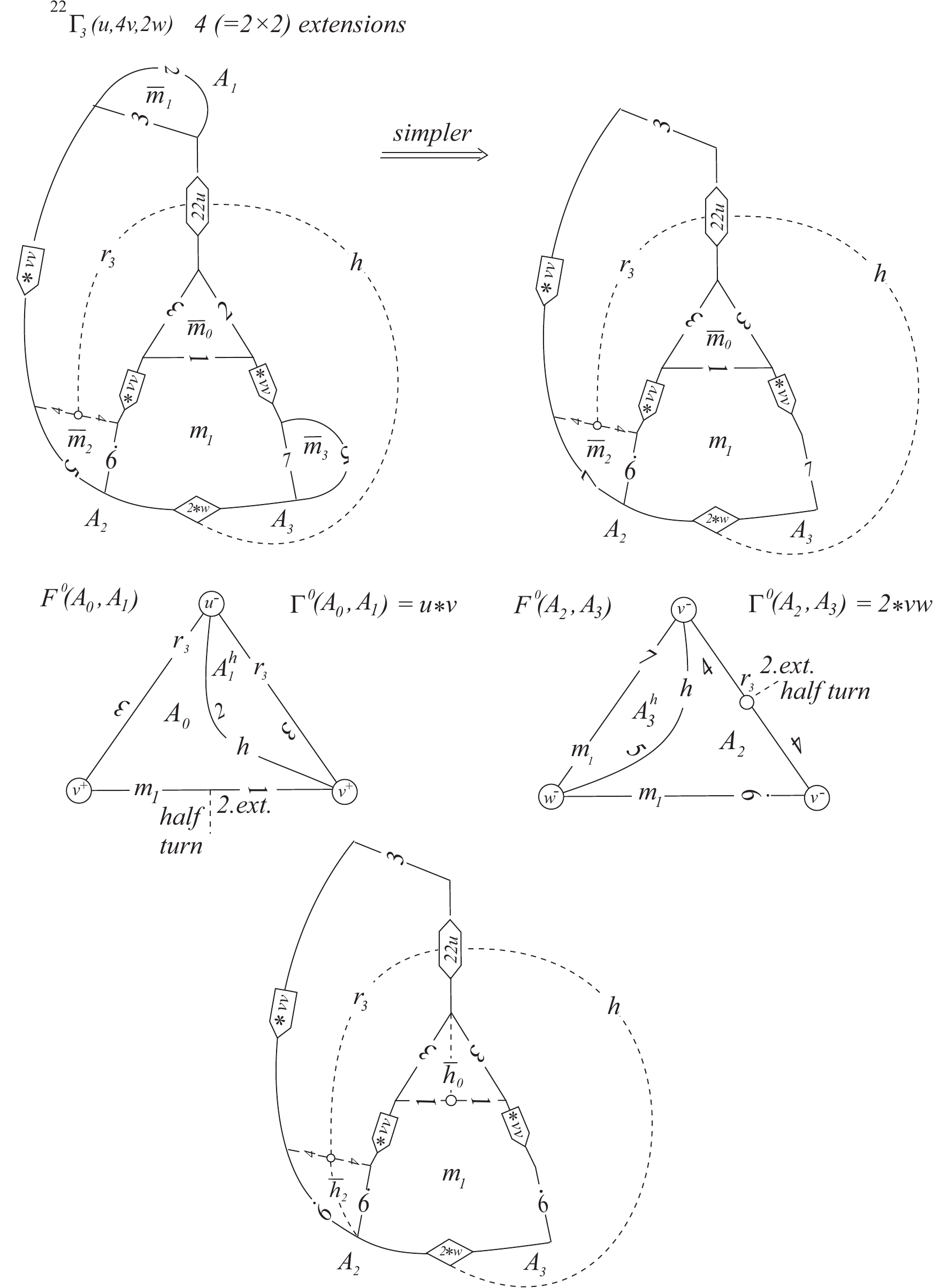} 
\caption{F4-5 to $ ^{22}\Gamma_3(u,4v,2w)$}
\label{fig:F4-5}
\end{figure}

The pictures are given in Fig.~\ref{fig:F4-5}: F4-5 with $\Gamma ^0(A_0, A_1) = u*v$ and $\Gamma ^0(A_2, A_3) = 2*vw$ as $\mathbb{H}^2$ groups. 

\vspace{3mm}

\noindent $\bullet \  ^{22} \Gamma _4 (u, 2v, w) = 
(r_1, r_3, h - 1 = r_1^2 = r_3^2 = h^2 = (r_3 h)^u = (r_1 r_3 h)^v = (r_1 h)^w, 
\ u \ne w, v \in \mathbb{N},  \ \dfrac{1}{u} + \dfrac{1}{v} < \dfrac{1}{2}, \ \dfrac{1}{v} + \dfrac{1}{w} < \dfrac{1}{2}) = ^2_2 \Gamma _4 (u, 2v, w)$, with 2 extensions both for vertices $A_0$, $A_1$ and $A_2$, $A_3$

\noindent {\bf (1)} for $A_0$, $A_1$: $(\bar{m}_0, \bar{m}_1 - 1 = \bar{m}_0^2 = \bar{m}_1^2 \stackrel{1}{=} (\bar{m}_0 r_1)^2 \stackrel{2}{=} \bar{m}_0 h \bar{m}_1 h \stackrel{3}{=} \bar{m}_0 r_3 \bar{m}_1 r_3 )$

\noindent {\bf (1)} for $A_2$, $A_3$: $(\bar{m}_2, \bar{m}_3 - 1 = \bar{m}_2^2 = \bar{m}_3^2 \stackrel{4}{=} (\bar{m}_2 r_3)^2 \stackrel{5}{=} \bar{m}_2 h \bar{m}_3 h \stackrel{6}{=} \bar{m}_2 r_1 \bar{m}_3 r_1)$

\noindent {\bf (2)} for $A_0$, $A_1$ by half-turn $\bar{h}_0$: 
$(\bar{h}_0 - 1 = \bar{h}_0^2 \stackrel{1}{=} (\bar{h}_0 r_1)^2 \stackrel{3}{=} (\bar{h}_0 h r_3)^2 )$

\noindent {\bf (2)} for $A_2$, $A_3$ by half-turn $\bar{h}_2$: $(\bar{h}_2 - 1=
\bar{h}_2^2 \stackrel{4}{=} (\bar{h}_2 r_3)^2 \stackrel{6}{=} (\bar{h}_2 h r_1)^2).$

\begin{figure}[htbp]
	\centering
\includegraphics[width=130mm]{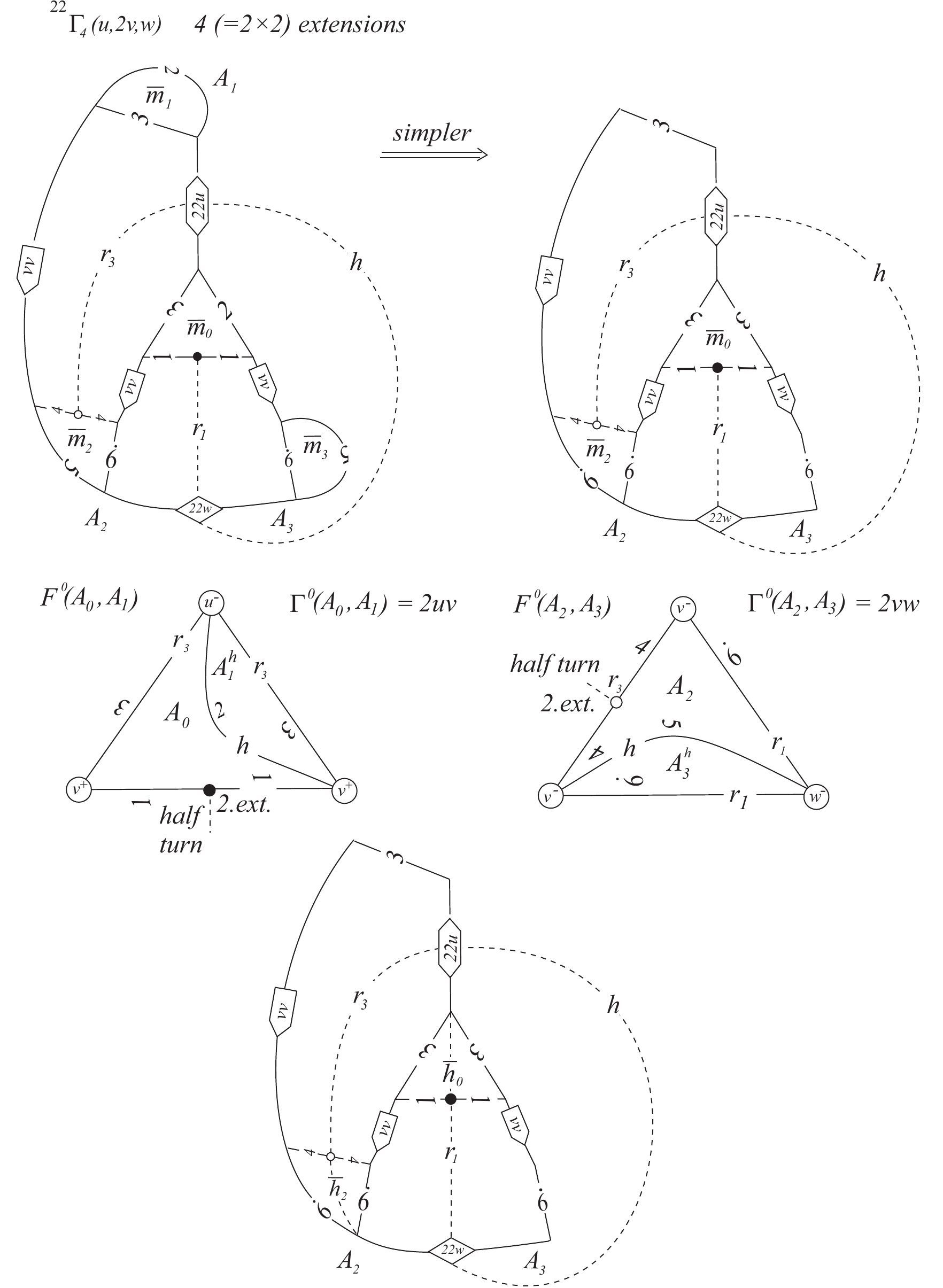}
\caption{F4-6 to $ ^{22}\Gamma_4(u,4v,w)$}
\label{fig:F4-6}
\end{figure}

The pictures are given in Fig.~\ref{fig:F4-6}: F4-6 with $\Gamma ^0(A_0, A_1) = 2uv$ and $\Gamma ^0(A_2, A_3) = 2vw$ as $\mathbb{H}^2$ groups. 
%


\section{On packings and coverings with congruent hyperballs in Family 1-4}
In this section we consider congruent hyperball packings and coverings for families 
F1-F2 and congruent hyperball packings for families F3-F4. 
The base planes of the hyperballs are always the truncation planes of the 
trunc-simplices.

\begin{enumerate}
\item The first we consider Family 1, with a simple truncated orthoscheme 
(e.g. $CHLA_1A_2A_3$ in Fig.~24.a and Fig.~1:F1-1) with Schl\"afli symbols 
$\{\overline{u},\overline{v},\overline{w}\}=$ $\{3,3,u\}$
$6< u\in\mathbb{N})$.
This simple truncated orthoscheme can be derived also by truncation from 
orthoscheme $A_0A_1A_2A_3=$ $b^0b^1b^2b^3$ with outer essential vertex $A_0$.
The truncating plane $a_0(\mbox{\boldmath$a$}_{0})=CLH$ is the polar plane of $A_0$, 
that {\it is the (ultraparallel to $b^0$) base plane  
of hyperball $\mathcal{H}^{h}_0$} with height $h_p$ for packing or $h_c$ for covering. 
\item For families F2-F4 (see Fig.~24.b and Fig.~8: F2-1, Fig.~17: F3-1, Fig.~20: F4-1, 
respectively) we consider a tiling $\mathcal{T}(\mathcal{O}(\overline{u},\overline{v},
\overline{w}))$ with Schl\"afli symbol $\{\overline{u},\overline{v},\overline{w}\}$, 
($\frac{1}{\overline{u}}+\frac{1}{\overline{v}} < 
\frac{1}{2},~\frac{1}{\overline{v}}+\frac{1}{\overline{w}} < \frac{1}{2}$, 
$\overline{u},\overline{v},\overline{w} \in \mathbb{N}$) whose fundamental 
domain is doubly truncated orthoscheme $\mathcal{O}(\overline{u},\overline{v},
\overline{w})=$ $CHLA_1A_2EJQ$ in Fig.~24.b. 

Let a truncated orthoscheme $\mathcal{O}(\overline{u},\overline{v},\overline{w})$ 
$\subset \HYP$ be a tile from 
the above tiling.  
This can be derived also by truncation from 
$A_0A_1A_2A_3=b^0b^1b^2b^3$ with essential outer vertices $A_0$ and $A_3$.
The truncating planes $a_0(\mbox{\boldmath$a$}_{0})=CLH$ and $a_3(\mbox{\boldmath$a$}_{3})=JEQ$ are the polar planes of $A_0$ and $A_3$, respectively,
that {\it will be (ultraparallel to $b^0$ and $b^3$) base planes 
of hyperballs $\mathcal{H}^{h}_i$} with height $h$ $(i=0,3)$. 
The distance between the two base planes 
$d(a_0(\mbox{\boldmath$a$}_{0}),a_3(\mbox{\boldmath$a$}_{3}))=d(H,J)$ ($d$ is the hyperbolic distance function) will be double of the height
of packing hyperball at most.
\end{enumerate}
\begin{figure}[htb]
\centering
\includegraphics[width=60mm]{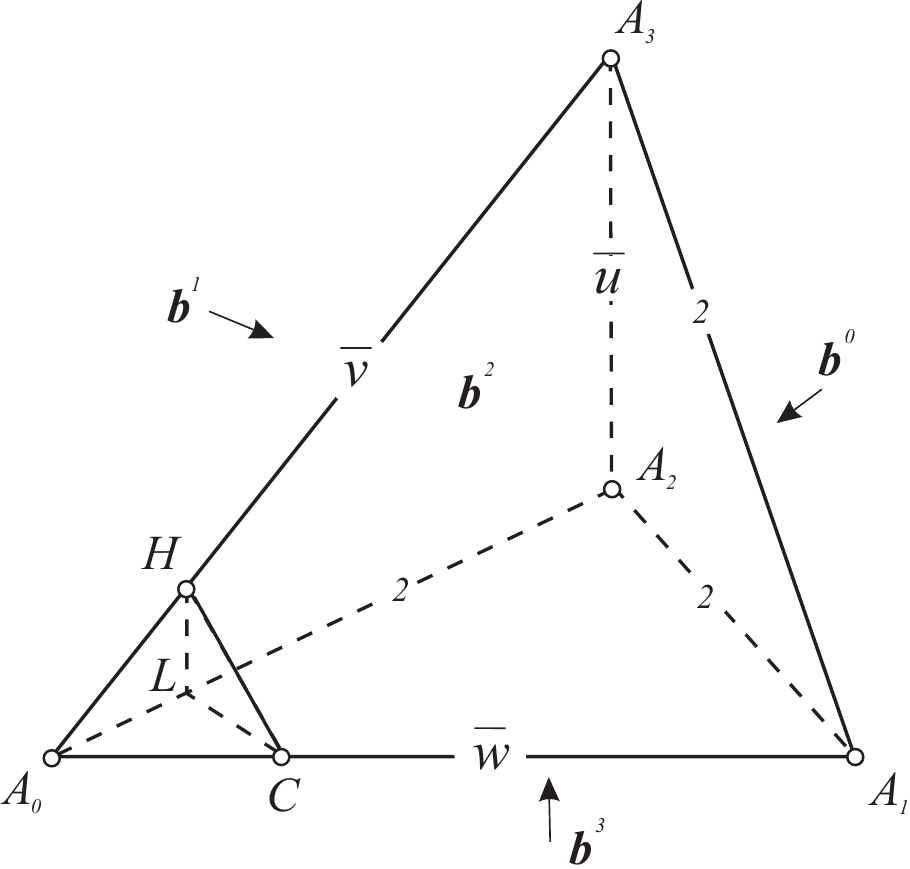} \includegraphics[width=60mm]{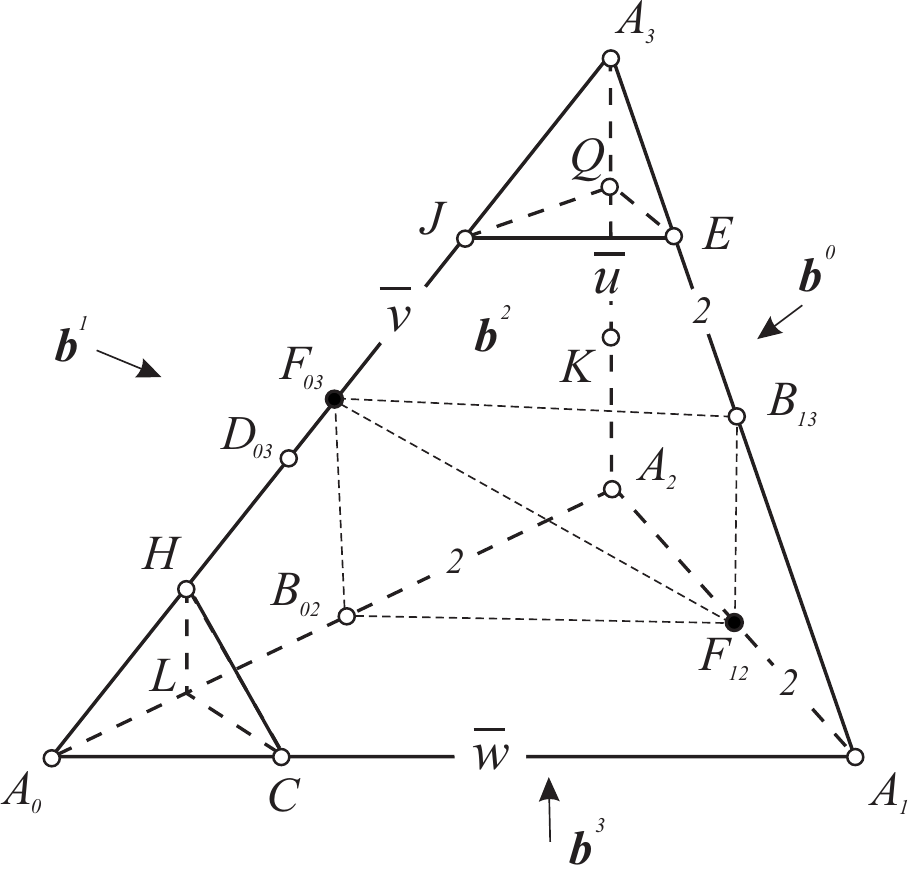}
\caption{Simple and double truncated complete orthoschemes (without the absolute) in standard coordinate simplex $b^0b^1b^2b^3=A_0A_1A_2A_3$, $(A_ib^j)=\delta_i^j$
with parameters $\beta^{01}=\frac{\pi}{\overline{u}}$, $\beta^{12}=\frac{\pi}{\overline{v}}$, $\beta^{23}=\frac{\pi}{\overline{w}}$. Sketchy pictures also with half-turn axis
$F_{03} F_{12}$ in case $\overline{u}=\overline{w}$}
\label{}
\end{figure}
The Coxeter-Schl\"afli symbol of the complete orthoscheme tiling $\mathcal{T}$ 
generated by reflections in the planes $b^i$ $(i\in \{0,1,2,3\})$ and $a_j$ $(j=0,3)$ 
of a complete orthoscheme $\mathcal{O}$.
To every scheme there is a corresponding
symmetric $4 \times 4$ matrix $(b^{ij})$ where $b^{ii}=1$ and, for $i \ne j\in \{0,1,2,3\}$,
$b^{ij}$ equals to $\cos(\pi-\beta^{ij})=-\cos{\beta^{ij}}$ for all dihedral 
angles $\beta^{ij}$ 
between the faces $b^i$,$b^j$ of $\mathcal{O}$.

For example, in formula (6.1) we see the 
Coxeter-Schl\"afli matrix with
parameters $(\overline{u},\overline{v},\overline{w})$, i.e. 
$\beta^{01}=\frac{\pi}{\overline{u}}$, 
$\beta^{12}=\frac{\pi}{\overline{v}}$, $\beta^{23}=\frac{\pi}{\overline{w}}$.
\[
(b^{ij})=(\langle \mbox{\boldmath$b^i$},\mbox{\boldmath$b^j$} \rangle):=\begin{pmatrix}
1& -\cos{\frac{\pi}{\overline{u}}}& 0 & 0 \\
-\cos{\frac{\pi}{\overline{u}}} & 1 & -\cos{\frac{\pi}{\overline{v}}}& 0 \\
0 & -\cos{\frac{\pi}{\overline{v}}} & 1 & -\cos{\frac{\pi}{\overline{w}}} \\
0 & 0 & -\cos{\frac{\pi}{\overline{w}}} & 1
\end{pmatrix}. \tag{6.1}
\]

This $3$-dimensional complete (truncated or frustum) orthoscheme $\mathcal{O}=\mathcal{O}(\overline{u},\overline{v},\overline{w})$ and 
its reflection group 
$\mathbf{G}_{\bar{u}\bar{v}\bar{w}}$ will be described in Fig.~24,1,8,17,20 and by the symmetric Coxeter-Schl\"afli matrix $(b^{ij})$ in formula (6.1), 
furthermore by its inverse matrix $(a_{ij})$ in formula
(6.2).
\[
\begin{gathered}
(a_{ij})=(b^{ij})^{-1}=\langle \ba_i, \ba_j \rangle:=\\
=\frac{1}{B} \begin{pmatrix}
\sin^2{\frac{\pi}{\overline{w}}}-\cos^2{\frac{\pi}{\overline{v}}}& \cos{\frac{\pi}{\overline{u}}}\sin^2{\frac{\pi}{\overline{w}}}& 
\cos{\frac{\pi}{\overline{u}}}\cos{\frac{\pi}{\overline{v}}} & 
\cos{\frac{\pi}{\overline{u}}}\cos{\frac{\pi}{\overline{v}}}\cos{\frac{\pi}{\overline{w}}} \\
\cos{\frac{\pi}{\overline{u}}}\sin^2{\frac{\pi}{\overline{w}}} & \sin^2{\frac{\pi}{w}} & \cos{\frac{\pi}{\overline{v}}}& \cos{\frac{\pi}{\overline{w}}}\cos{\frac{\pi}{\overline{v}}} \\
\cos{\frac{\pi}{\overline{u}}}\cos{\frac{\pi}{\overline{v}}} & \cos{\frac{\pi}{\overline{v}}} & \sin^2{\frac{\pi}{\overline{u}}}  & 
\cos{\frac{\pi}{w}}\sin^2{\frac{\pi}{\overline{u}}}  \\
\cos{\frac{\pi}{\overline{u}}}\cos{\frac{\pi}{\overline{v}}}\cos{\frac{\pi}{\overline{w}}}  & \cos{\frac{\pi}{w}}\cos{\frac{\pi}{\overline{v}}} & 
\cos{\frac{\pi}{\overline{w}}}\sin^2{\frac{\pi}{\overline{u}}}  
& \sin^2{\frac{\pi}{\overline{u}}}-\cos^2{\frac{\pi}{\overline{v}}}
\end{pmatrix}, \tag{6.2}
\end{gathered}
\]
where
$$
B=\det(b^{ij})=\sin^2{\frac{\pi}{\overline{u}}}\sin^2{\frac{\pi}{\overline{w}}}-\cos^2{\frac{\pi}{\overline{v}}} <0, \ \ \text{i.e.} \ 
\sin{\frac{\pi}{\overline{u}}}\sin{\frac{\pi}{\overline{w}}}-\cos{\frac{\pi}{\overline{v}}}<0.
$$
In the following the volume of $\mathcal{O}(\overline{u},\overline{v},\overline{w})$ 
is derived by the next
Theorem of {{R.~Kellerhals}} (\cite{K89}, by the ideas of N.~I.~Lobachevsky):
\begin{theorem}{\rm{(R.~Kellerhals)}} The volume of a three-dimensional hyperbolic
complete orthoscheme $\mathcal{O}=\mathcal{O}(\overline{u},\overline{v},\overline{w}) 
\subset \mathbb{H}^3$
can be expressed with the essential
angles $\beta^{01}=\frac{\pi}{\overline{u}}$, $\beta^{12}=\frac{\pi}{\overline{v}}$, $\beta^{23}=\frac{\pi}{\overline{w}}$, $(0 \le \beta^{ij}
\le \frac{\pi}{2})$ in the following form:

\begin{align}
&\mathrm{Vol}(\mathcal{O})=\frac{1}{4} \{ \mathcal{L}(\beta^{01}+\theta)-
\mathcal{L}(\beta^{01}-\theta)+\mathcal{L}(\frac{\pi}{2}+\beta^{12}-\theta)+ \notag \\
&+\mathcal{L}(\frac{\pi}{2}-\beta^{12}-\theta)+\mathcal{L}(\beta^{23}+\theta)-
\mathcal{L}(\beta^{23}-\theta)+2\mathcal{L}(\frac{\pi}{2}-\theta) \}, \notag
\end{align}
where $\theta \in (0,\frac{\pi}{2})$ is defined by:
$$
\tan{\theta}=\frac{\sqrt{ \cos^2{\beta^{12}}-\sin^2{\beta^{01}} \sin^2{\beta^{23}
}}} {\cos{\beta^{01}}\cos{\beta^{23}}},
$$
and where $\mathcal{L}(x):=-\int\limits_0^x \log \vert {2\sin{t}} \vert dt$ \ denotes the
Lobachevsky function (in J. Milnor's interpretation).
\end{theorem}
The hypersphere (or equidistant surface) is a quadratic surface at a constant distance
from a plane (base plane) in both halfspaces. The infinite body of the hypersphere, containing the base plane, is called {\it hyperball}.

The {\it half hyperball } (i.e., the part of the hyperball lying on one side of its base plane) with distance $h$ to a base plane $\beta$ is denoted by $\mathcal{H}^h_+$.
The volume of the intersection of $\mathcal{H}^h_+(\mathcal{A})$ and the right prism with
base a $2$-polygon $\mathcal{A} \subset \beta$, can be determined by the classical formula 
of J.~Bolyai.
\begin{equation}
\mathrm{Vol}(\mathcal{H}^h_+(\mathcal{A}))=\frac{1}{4}\mathrm{Area}(\mathcal{A})\left[k \sinh \frac{2h}{k}+
2 h \right]. \tag{6.3}
\end{equation}
The constant $k =\sqrt{\frac{-1}{K}}$ is the natural length unit in
$\mathbb{H}^3$, where $K$ denotes the constant negative sectional curvature. In the following we may assume that $k=1$.
The distance $d$ of two proper points
$X(\mathbf{x})$ and $Y(\mathbf{y})$ is calculated by the formula
\begin{equation}
\cosh{{d}}=\frac{-\langle ~ \mathbf{x},~\mathbf{y} \rangle }{\sqrt{\langle ~ \mathbf{x},~\mathbf{x} \rangle
\langle ~ \mathbf{y},~\mathbf{y} \rangle }}. \tag{6.4}
\end{equation}
\subsection{Congruent hyperballs to a simply truncated orthoscheme, Family 1}
The maximal series of family F1 is characterized by the simplex tiling whose group series is $^{\bar{4}3m}_{\ \ 24} \Gamma (u)$ 
with previous crystallographic notation from \cite{MPS06}, and now with the new orbifold notation $^{*233} \Gamma (u)$. 
It has the trivial extension by the reflection $a_0 = \bar{m}_0$, as in Sect.~2.

The corresponding unique fundamental simplex with parameters $(\overline{u},\overline{v}, \overline{w})$=$(3,3,u)$ is determined by its Coxeter-Schl\"{a}fli matrix which is 
derived by (6.1). Its determinant is $B = \det ( b^{ij}) = \frac{3}{4} \sin^2 \frac{\pi}{u} - \frac{1}{4} < 0$
now $6 < u \in \mathbb{N}$.

So, the signature is $(+++-)$ and the simplex is realizable in hyperbolic 
space $\mathbb{H}^3$. The inverse $(a_{ij})$ of the Coxeter-Schl\"afli  
matrix $(b^{ij})$ follows by (6.2) providing the distance metric (by the general theory of projective metrics \cite{MS18}). 

The optimal hyperball height is $h_p=d(a_0, m_0) = d(C, A_1)$ for packings and $h_c=d(H, A_3)$ is the optimal covering distance where
\begin{equation*}
\begin{gathered}
 h_p=\cosh d(C, A_1) =  \frac{- {a}_{00} {a}_{11} + {a}_{01}^2}{\sqrt{{a}_{00} {a}_{11} \left( - {a}_{01}^2 + {a}_{00} {a}_{11}\right)} } 
 = \sqrt{\frac{\frac{1}{4} - \frac{3}{4} \sin^2 \frac{\pi}{u}}{\frac{1}{4} - 
 \sin^2 \frac{\pi}{u}}},
\end{gathered}
\end{equation*}
\begin{equation*}
\begin{gathered}
h_c=\cosh d(H, A_3) = \frac{- {a}_{00} {a}_{33} + {a}_{03}^2}{\sqrt{{a}_{00} {a}_{33} \left( - {a}_{03}^2 + {a}_{00} {a}_{33}\right)} }=
\sqrt{\frac{\frac{1}{2} \left( \frac{1}{4} - \sin^2 \frac{\pi}{\overline{w}} \right) + \frac{1}{16} \cos^2 \frac{\pi}{u}}{\frac{1}{2} \left( \frac{1}{4} - 
\sin^2 \frac{\pi}{u} \right)} } \\
\text{equivalent to}~\sinh d(C, {A}_1) = \frac{\frac{1}{2} \sin \frac{\pi}{u}}{\sqrt{\frac{1}{4} - \sin^2 \frac{\pi}{u}}}~\text{and}~ 
\sinh d(H, A_3) =  \frac{\frac{1}{4} \cos \frac{\pi}{u}}{\sqrt{ \frac{1}{2} 
\left( \frac{1}{4} - \sin^2 \frac{\pi}{u} \right) }}.
\end{gathered}
\end{equation*}

Then we get the densities $\delta$ (packing), 
$\Delta$ (covering), respectively. E.g.
$$\delta(\mathcal{O}(3,3,u)) = \frac{\mathrm{Area}(\mathcal{A}) \frac{1}{4} \left[ \sinh{2h_p} + 2h_p \right] }{\mathrm{Vol}(\mathcal{O})},$$
where $\mathrm{Area}(\mathcal{A}) = \left( \dfrac{\pi}{2} - \dfrac{\pi}{3} - \dfrac{\pi}{u} \right)$. 

The $\mathrm{Vol}(\mathcal{O})$ can be calculated by Theorem 3.

In this case the maximal packing density is $\approx 0.82251$ with $u = 7$ (see in Table 1p). The optimal covering density of Family 1 is $\approx 1.33093$ for $u = 7$ in Table 1c.
\medbreak
{\scriptsize
\centerline{\vbox{
\halign{\strut\vrule~\hfil $#$ \hfil~\vrule
&\quad \hfil $#$ \hfil~\vrule
&\quad \hfil $#$ \hfil\quad\vrule
&\quad \hfil $#$ \hfil\quad\vrule
&\quad \hfil $#$ \hfil\quad\vrule
\cr
\noalign{\hrule}
\multispan5{\strut\vrule\hfill\bf Table 1p,~$(3,3,u),~$ $6 < u \in \mathbb{N}$ \hfill\vrule}%
\cr
\noalign{\hrule}
\noalign{\vskip2pt}
\noalign{\hrule}
u & h & \mathrm{Vol}(\mathcal{O}) & \mathrm{Vol}(\mathcal{H}^h_+(\mathcal{A}))& 
\delta \cr
\noalign{\hrule}
7 & 0.78871 & 0.08856 & 0.07284 & {\bf 0.82251} \cr
\noalign{\hrule}
8 & 0.56419 & 0.10721 & 0.08220 & 0.76673 \cr
\noalign{\hrule}
9 & 0.45320 & 0.11825 & 0.08474 & 0.71663 \cr
\noalign{\hrule}
\vdots & \vdots  & \vdots  & \vdots  & \vdots \cr
\noalign{\hrule}
20 & 0.16397 & 0.14636 & 0.06064 & 0.41431 \cr
\noalign{\hrule}
\vdots & \vdots  & \vdots  & \vdots  & \vdots \cr
\noalign{\hrule}
50 & 0.06325 & 0.15167 & 0.02918 & 0.19240 \cr
\noalign{\hrule}
\vdots & \vdots  & \vdots  & \vdots  & \vdots \cr
\noalign{\hrule}
100 & 0.03147 & 0.15241 & 0.01549 & 0.10165 \cr
\noalign{\hrule}
u \to \infty & 0 & 0.15266 & 0 & 0 \cr
\noalign{\hrule}}}}}
\medbreak
{\scriptsize
\centerline{\vbox{
\halign{\strut\vrule~\hfil $#$ \hfil~\vrule
&\quad \hfil $#$ \hfil~\vrule
&\quad \hfil $#$ \hfil\quad\vrule
&\quad \hfil $#$ \hfil\quad\vrule
&\quad \hfil $#$ \hfil\quad\vrule
\cr
\noalign{\hrule}
\multispan5{\strut\vrule\hfill\bf Table 1c,~$(3,3,u)$~$6 < u \in \mathbb{N}$ \hfill\vrule}%
\cr
\noalign{\hrule}
\noalign{\vskip2pt}
\noalign{\hrule}
u & h & \mathrm{Vol}(\mathcal{O}) & \mathrm{Vol}(\mathcal{H}^h_+(\mathcal{A}))& \Delta \cr
\noalign{\hrule}
7 & 1.06739 & 0.08856 & 0.11787 & {\mathbf{1.33093}} \cr
\noalign{\hrule}
8 & 0.89198 & 0.10721 & 0.15304 & 1.42747 \cr
\noalign{\hrule}
9 & 0.81696 & 0.11825 & 0.17882 & 1.51225 \cr
\noalign{\hrule}
\vdots & \vdots  & \vdots  & \vdots  & \vdots \cr
\noalign{\hrule}
20 & 0.68136 & 0.14636 & 0.29213 & 1.99596 \cr
\noalign{\hrule}
\vdots & \vdots  & \vdots  & \vdots  & \vdots \cr
\noalign{\hrule}
50 & 0.66193 & 0.15167 & 0.35361 & 2.33146 \cr
\noalign{\hrule}
\vdots & \vdots  & \vdots  & \vdots  & \vdots \cr
\noalign{\hrule}
100 & 0.65934 & 0.15241 & 0.37580 & 2.46566 \cr
\noalign{\hrule}
u \to \infty & 0.65848 & 0.15266 & 0.39911 & 2.61438 \cr
\noalign{\hrule}}}}}
\smallbreak
By completing Tables 1p and 1c we can read the optimal hyperball packings and coverings with their densities for other extensions in F1 as well.
E.g. to $^{233}\Gamma(2u)$ we get these for $2u=8$ with $\delta=\approx 0.76673$, $\Delta \approx 1.42747$.
To $^{2\times}\Gamma_3(3u)$ we get these for $3u=9$ with $\delta \approx 0.71663$, $\Delta \approx 1.51225$.
To $^{33}\Gamma_5(4u)$ we get these for $4u=8$ again with $\delta \approx 0.76673$, $\Delta \approx 1.42747$.

The argument is by our convention that the parameters of groups above in
parentheses are just the $m^{23}$ entries of the $D-$matrix function, $2m^{23}$ is the number
of 
subsimplices around simplex edges. 
Mostow's rigidity principle for co-compact hyperbolic space groups are 
also used (see e.g. \cite{V}).
\subsection{On hyperball packings and coverings in a doubly truncated orthoscheme, 
Family 2}
%
In this case we consider a doubly truncated orthoscheme with 
Schl\"afli symbol $\{\overline{u},\overline{v},\overline{w}=\overline{u} \}$, 
($\frac{1}{\overline{u}}+\frac{1}{\overline{v}} < \frac{1}{2}$, 
$3 \le \overline{u},\overline{v} \in \mathbb{N}$) whose fundamental 
domain is $CHLA_1A_2EJQ$ in Fig.~24.b. 
Its Coxeter Schl\"afli matrix $(b^{ij})$ and its inverse $(a_{ij})$ can be derived by formulas (6.1) and (6.2).

Volume formula of half hyperball $H^h_+(\mathcal{A}_i)$ $(i=0,3)$ of height $h$ can be calculated by (6.3) 
where $\mathrm{Area}(\mathcal{\cA}_0) = \mathrm{Area}(\mathcal{\cA}_3) = \pi \left( \dfrac{1}{2} - \dfrac{1}{\bar{u}} - \dfrac{1}{\bar{v}} \right)$ according to the above 
orthoscheme.

Both polar planes assigne hyperspheres that are congruent with each other therefore the heights $h^{0}=h^3$ of optimal hyperballs $\mathcal{H}_i^{h^{i}}$ $(i=0,3)$ can be computed by
Fig.~24.b.

{\bf The maximal height for optimal congruent hyperball packing} is $h_p = \min \left\{ d({A}_2, {a}_3), d(F_{03}, {a}_3)  \right\}$.

In the above doubly truncated orthoscheme 
$d({A}_2, {a}_3) = d({A}_2, Q \sim {a}_{23} {\ba}_3 - {a}_{33} {\ba}_2)$, $d(F_{03}, {a}_3) = d({\ba}_0 + {\ba}_3, {a}_{33} {\ba}_0 - {a}_{03} 
{\ba}_3 \sim J)$, (as in \cite{MS18}).
Thus, 
\begin{equation*}
\begin{gathered}
 \cosh \left( d({A}_2,{a}_3) \right) = \frac{-{a}_{33} {a}_{22} + {a}_{23}^2 } {\sqrt{{a}_{22} \left( {a}_{33}^2 {a}_{23} - {a}_{23}^2 {a}_{33} \right)}} = 
 \sqrt{1-\frac{{a}_{23}^2 }{a_{22} a_{33}}}, 
\end{gathered}
\end{equation*}
and
\begin{equation*}
\begin{gathered}
 \cosh \left( d (F_{03}, {a}_3) \right) = \sqrt{\frac{a_{33} - a_{03} }{2 a_{33}}} = 
 \sqrt{\frac{1}{2} \left( 1 + \frac{\cos^2 \frac{\pi}{\bar{u}} \cos \frac{\pi}{\bar{v}} }{\cos^2 \frac{\pi}{\bar{v}} - \sin^2 \frac{\pi}{\bar{u}} } \right)}.
\end{gathered}
\end{equation*}

The density of optimal hyperball packing will be determined, by the following formula:
\begin{equation}
\begin{gathered}
\delta(\mathcal{O}(\overline{u},\overline{v},\overline{u})) = \frac{ \mathrm{Area}(\cA) ( \sinh {2h_p} +2h_p)}{2\mathrm{Vol}(\mathcal{O})} 
= \frac{\pi \left( \frac{1}{2} - \frac{1}{\bar{u}} - \frac{1}{\bar{v}} \right) \left( \sinh {2h_p} + 2h_p \right)}{2 \mathrm{Vol}(\mathcal{O})}. \tag{6.5} 
\end{gathered}
\end{equation}

{\bf For the minimal covering hyperball height}:
\begin{equation*}
\begin{gathered}
h_c = d(F_{12}, {a}_0) = d(F_{12}, {a}_3) = d({\ba}_1 + \ba_2, \ba_1 + \ba_2 + c \ba_3),
\end{gathered}
\end{equation*}
where $a_{13} + a_{23} + c a_{33} = 0$, i.e. $c = \dfrac{- (a_{13} + a_{23})}{a_{33}}$. Thus, by $a_{11} = a_{22}$, we get  
\begin{equation*}
\begin{gathered}
\cosh \left( d(F_{12}, {a}_3 \right)  =  \frac{- \left( 2 \left( a_{22} + a_{12} \right) - \frac{1}{a_{33}} \left( a_{13} + a_{23}  \right)^2 \right)} 
{\sqrt{ 2 \left( a_{22} + a_{12} \right) 
\left( 2 a_{22} + 2 a_{12} - \frac{1}{a_{33}} \left( a_{13} + a_{23} \right)^2 \right)}}. 
\end{gathered}
\end{equation*}

The covering density $\Delta(\mathcal{O}(\overline{u},\overline{v},\overline{u}))$ 
can be defined similarly to the packing density (see (6.5)).

$Vol(\mathcal{O}(\overline{u},\overline{v},\overline{u}))$ can be calculated by Theorem 3. The maximal volume sum of the hyperball 
pieces lying in ${\mathcal{O}(\overline{u},\overline{v},\overline{u})}$ can be computed by the formulas (6.1-4) and by the above described computation method 
for each given possible parameters $\overline{u},\overline{v},\overline{u}$. 
Therefore, the maximal packing density and the thinnest covering density of the congruent hyperball packings and coverings can be computed for each possible parameters.

We see the optimal packing density for F2 is $\approx 0,81335$ with $\overline{u}=\overline{w} = 7$, $\overline{v} = 3$ in Table 2p and the top minimal covering density is 
$\approx 1,26869$ with $\overline{u} = 7$, $\overline{v} = 3$ in Table 2c.
\medbreak
{\scriptsize
\centerline{\vbox{
\halign{\strut\vrule~\hfil $#$ \hfil~\vrule
&\quad \hfil $#$ \hfil~\vrule
&\quad \hfil $#$ \hfil\quad\vrule
&\quad \hfil $#$ \hfil\quad\vrule
&\quad \hfil $#$ \hfil\quad\vrule
\cr
\noalign{\hrule}
\multispan5{\strut\vrule\hfill\bf Table 2p,~$(\overline{u},\overline{v},\overline{u})$ \hfill\vrule}%
\cr
\noalign{\hrule}
\noalign{\vskip2pt}
\noalign{\hrule}
\{\overline{u},\overline{v},\overline{u}\} & h_p  & Vol({\mathcal{O}(\overline{u},\overline{v},\overline{u})}) & \sum_{i=0,3} {Vol}(\mathcal{H}^h(\mathcal{A}_i))\\
& \delta({\mathcal{O}(\overline{u},\overline{v},\overline{u})})\cr
\noalign{\hrule}
\{7,3,7\} & 1.23469 & 0.38325 & 0.31172 & {\mathbf{0.81335}} \cr
\noalign{\hrule}
\{6,4,6\} & 0.69217; & 0.55557 & 0.42610 & 0.76696 \cr
\noalign{\hrule}
\{8,3,8\} & 0.94946 & 0.44383 & 0.33794 & 0.76143 \cr
\noalign{\hrule}
\{8,4,8\} & 0.56419 & 0.64328 & 0.49322 & 0.76673 \cr
\noalign{\hrule}
\{5,4,5\} & 0.88055 & 0.46190 & 0.36007 & 0.77955 \cr
\noalign{\hrule}
\{4,5,4\} & 0.80846 & 0.43062 & 0.31702 & 0.73620 \cr
\noalign{\hrule}
\{4,6,4\} & 0.57311 & 0.50192 & 0.33516 & 0.66775 \cr
\noalign{\hrule}
\{3,7,3\} & 0.98399 & 0.27899 & 0.20481 & 0.73411 \cr
\noalign{\hrule}}}}}
\medbreak
{\scriptsize
\centerline{\vbox{
\halign{\strut\vrule~\hfil $#$ \hfil~\vrule
&\quad \hfil $#$ \hfil~\vrule
&\quad \hfil $#$ \hfil\quad\vrule
&\quad \hfil $#$ \hfil\quad\vrule
&\quad \hfil $#$ \hfil\quad\vrule
\cr
\noalign{\hrule}
\multispan5{\strut\vrule\hfill\bf Table 2c,~$(\overline{u},\overline{v},\overline{u})$  \hfill\vrule}%
\cr
\noalign{\hrule}
\noalign{\vskip2pt}
\noalign{\hrule}
\{\overline{u},\overline{v},\overline{u}\} & h_c  & Vol({\mathcal{O}(\overline{u},\overline{v},\overline{u})}) & \sum_{i=0,3} {Vol}(\mathcal{H}^h(\mathcal{A}_i))\\
& \delta({\mathcal{O}(\overline{u},\overline{v},\overline{u})})\cr
\noalign{\hrule}
\{\mathbf{7,3,7}\} & 1.49903 & 0.38325 & 0.48607 & {\mathbf{1.26829}} \cr
\noalign{\hrule}
\{6,4,6\} & 1.01481 & 0.55557 & 0.75523 & 1.35938 \cr
\noalign{\hrule}
\{8,3,8\} & 1.26595 & 0.44383 & 0.57470 & 1.29487 \cr
\noalign{\hrule}
\{5,4,5\} & 1.19095 & 0.46190 & 0.60856 & 1.31751 \cr
\noalign{\hrule}
\{8,4,8\} & 0.89198 & 0.64328 & 0.91826 & 1.42747 \cr
\noalign{\hrule}
\{4,5,4\} & 1.16974 & 0.43062 & 0.58741 & 1.36411 \cr
\noalign{\hrule}
\{4,6,4\} & 0.99583 & 0.50192 & 0.73137 & 1.45714 \cr
\noalign{\hrule}
\{3,7,3\} & 1.36406 & 0.27899 & 0.38699 & 1.38713 \cr
\noalign{\hrule}}}}}
\medbreak
By completing Tables 2p and 2c we can determine the other densities in Family 2. E.g. to $^{*22}\Gamma_2(2u,2v)$ (i.e. $u=4,v=3$) we get
$\delta \approx 0.76143$. $\Delta \approx 1.29487$. Or to $\Gamma_{58}(2u,4v)$ we get $\{\overline{u},\overline{v},\overline{w}\}=\{6,4,6\}$ 
for optimal arrangement, $\delta \approx 0.76696$. $\Delta \approx 1.35938$.
\subsection{Packings with congruent hyperballs in doubly truncated orthoscheme, Family 3-4}
In this subsection we consider congruent hyperball packings to truncated orthoscheme tilings.
Both polar planes assign hyperballs with equal height (see Fig.~24.b). It is clear, 
that the height of optimal hyperballs $\mathcal{H}_i^{h^{i}}$ $(i=0,3)$ is
\begin{equation}
\begin{gathered}
h_p(\overline{u},\overline{v},\overline{w}) 
=\min\{d(H,J)/2,d(Q,A_2),d(C,A_1)\}, \notag
\end{gathered}
\end{equation}
where $\overline{u},\overline{v},\overline{w}$ are integer parameters. 
In this case the volume sum of the hyperball pieces 
lying in the orthoscheme has to divide with the volume of trunc-orthoscheme $\mathcal{O}(\overline{u},\overline{v},\overline{w})$ as usual.
Segments $A_1C$, $A_2Q$ and $JH$ 
can be determined by the machinery of the projective metrics (see subsection 6.1 or \cite{MS18}).

The volume 
of the orthoscheme $\mathcal{O}(\overline{u},\overline{v},\overline{w})$ can be determined by Theorem 3.
We note here, that the role of the parameters $u$ and $w$ is symmetrical, 
therefore we can assume that $\overline{u} > \overline{w}$ in F3 or $\overline{u}< \overline{w}$ in F4.
\medbreak
{\scriptsize
\centerline{\vbox{
\halign{\strut\vrule~\hfil $#$ \hfil~\vrule
&\quad \hfil $#$ \hfil~\vrule
&\quad \hfil $#$ \hfil\quad\vrule
&\quad \hfil $#$ \hfil\quad\vrule
&\quad \hfil $#$ \hfil\quad\vrule
\cr
\noalign{\hrule}
\multispan5{\strut\vrule\hfill\bf Table 3p,~$(\overline{u},\overline{v},\overline{w}=3)$, congruent hyperballs  \hfill\vrule}%
\cr
\noalign{\hrule}
\noalign{\vskip2pt}
\noalign{\hrule}
\{\overline{u},\overline{v},\overline{w}\} & h  & Vol({\mathcal{O}(\overline{u},\overline{v},\overline{w})}) & \sum_{i=0,3} {Vol}(\mathcal{H}^h(\mathcal{A}_i))\\
& \delta^1({\mathcal{O}(\overline{u},\overline{v},\overline{w})})\cr
\noalign{\hrule}
\{4,7,3\} & 0.59710 & 0.39274 & 0.27700 & 0.70529  \cr
\noalign{\hrule}
\{5,7,3\} & 0.41812; & 0.43216 & 0.25203 & 0.58320 \cr
\noalign{\hrule}
\vdots & \vdots  & \vdots  & \vdots  & \vdots  \cr
\noalign{\hrule}
\{50,7,3\} & 0.03492 & 0.49140 & 0.03962 & 0.08062 \cr
\noalign{\hrule}
\vdots & \vdots  & \vdots  & \vdots  & \vdots  \cr
\noalign{\hrule}
\{4,8,3\} & 0.56419 & 0.42885 & 0.32881 & {\mathbf{0.76673}} \cr
\noalign{\hrule}
\{5,8,3\} & 0.67409& 0.47536 & 0.26747 & 0.56266 \cr
\noalign{\hrule}
\vdots & \vdots  & \vdots  & \vdots  & \vdots  \cr
\noalign{\hrule}
\{50,8,3\} & 0.03405 & 0.52378 & 0.04245 & 0.08105 \cr
\noalign{\hrule}
\vdots & \vdots  & \vdots  & \vdots  & \vdots  \cr
\noalign{\hrule}
\{4,9,3\} & 0.46841 & 0.45130 & 0.30800 & 0.68247 \cr
\noalign{\hrule}
\{5,9,3\} & 0.39083 & 0.48771 & 0.31589 & 0.64771 \cr
\noalign{\hrule}
\vdots & \vdots  & \vdots  & \vdots  & \vdots  \cr
\noalign{\hrule}
\{50,9,3\} & 0.03348& 0.54384 & 0.04466 & 0.08212 \cr
\noalign{\hrule}}}}}
\medbreak
{\scriptsize
\centerline{\vbox{
\halign{\strut\vrule~\hfil $#$ \hfil~\vrule
&\quad \hfil $#$ \hfil~\vrule
&\quad \hfil $#$ \hfil\quad\vrule
&\quad \hfil $#$ \hfil\quad\vrule
&\quad \hfil $#$ \hfil\quad\vrule
\cr
\noalign{\hrule}
\multispan5{\strut\vrule\hfill\bf Table 4p,~$(\overline{u},\overline{v},\overline{w})$, congruent hyperballs  \hfill\vrule}%
\cr
\noalign{\hrule}
\noalign{\vskip2pt}
\noalign{\hrule}
\{\overline{u},\overline{v},\overline{w}\} & h  & Vol({\mathcal{O}(\overline{u},\overline{v},\overline{w})}) & \sum_{i=0,3} {Vol}(\mathcal{H}^h(\mathcal{A}_i))\\
& \delta^1({\mathcal{O}(\overline{u},\overline{v},\overline{w})})\cr
\noalign{\hrule}
\{7,3,8\} & 0.93100 & 0.41326 & 0.25726 & 0.62251  \cr
\noalign{\hrule}
\{7,3,9\} & 0.76734 & 0.43171 & 0.23355 & 0.54099 \cr
\noalign{\hrule}
\vdots & \vdots  & \vdots  & \vdots  & \vdots  \cr
\noalign{\hrule}
\{7,3,50\} & 0.11380 & 0.49016 & 0.06121 & 0.12488 \cr
\noalign{\hrule}
\vdots & \vdots  & \vdots  & \vdots  & \vdots  \cr
\noalign{\hrule}
\{8,3,9\} & 0.78366 & 0.46266 & 0.29474 & 0.63704 \cr
\noalign{\hrule}
\{8,3,10\} & 0.67409& 0.47536 & 0.26747 & 0.56266 \cr
\noalign{\hrule}
\vdots & \vdots  & \vdots  & \vdots  & \vdots  \cr
\noalign{\hrule}
\{8,3,50\} & 0.11668 & 0.52248 & 0.06935 & 0.13274 \cr
\noalign{\hrule}
\vdots & \vdots  & \vdots  & \vdots  & \vdots  \cr
\noalign{\hrule}
\{5,4,6\} & 0.73969 & 0.50747 & 0.37287 & 0.73476 \cr
\noalign{\hrule}
\{5,4,7\} & 0.59326 & 0.53230 & 0.32974 & 0.61947 \cr
\noalign{\hrule}
\vdots & \vdots  & \vdots  & \vdots  & \vdots  \cr
\noalign{\hrule}
\{5,4,50\} & 0.07206 & 0.59291 & 0.06350 & 0.10710 \cr
\noalign{\hrule}
\vdots & \vdots  & \vdots  & \vdots  & \vdots  \cr
\noalign{\hrule}
\{4,5,5\} & 0.69129 & 0.49789 & 0.38284 & {\mathbf{0.76893}} \cr
\noalign{\hrule}
\{4,5,6\} & 0.53064 & 0.52971 & 0.33597 & 0.63426 \cr
\noalign{\hrule}
\vdots & \vdots  & \vdots  & \vdots  & \vdots  \cr
\noalign{\hrule}
\{4,5,50\} & 0.05502 & 0.59318 & 0.05710 & 0.096256 \cr
\vdots & \vdots  & \vdots & \vdots  & \vdots  \cr
\noalign{\hrule}
\{4,6,5\} & 0.50625 &  0.55992 &  0.37558 & 0.67078 \cr
\noalign{\hrule}
\{4,6,6\} & 0.48121 & 0.58850 & 0.40850 & 0.69414 \cr
\noalign{\hrule}
\vdots & \vdots  & \vdots  & \vdots  & \vdots  \cr
\noalign{\hrule}
\{4,6,50\} & 0.05138 & 0.64697 & 0.06409 & 0.09906 \cr
\noalign{\hrule}}}}}
\begin{figure}[ht]
\centering
\includegraphics[width=120mm]{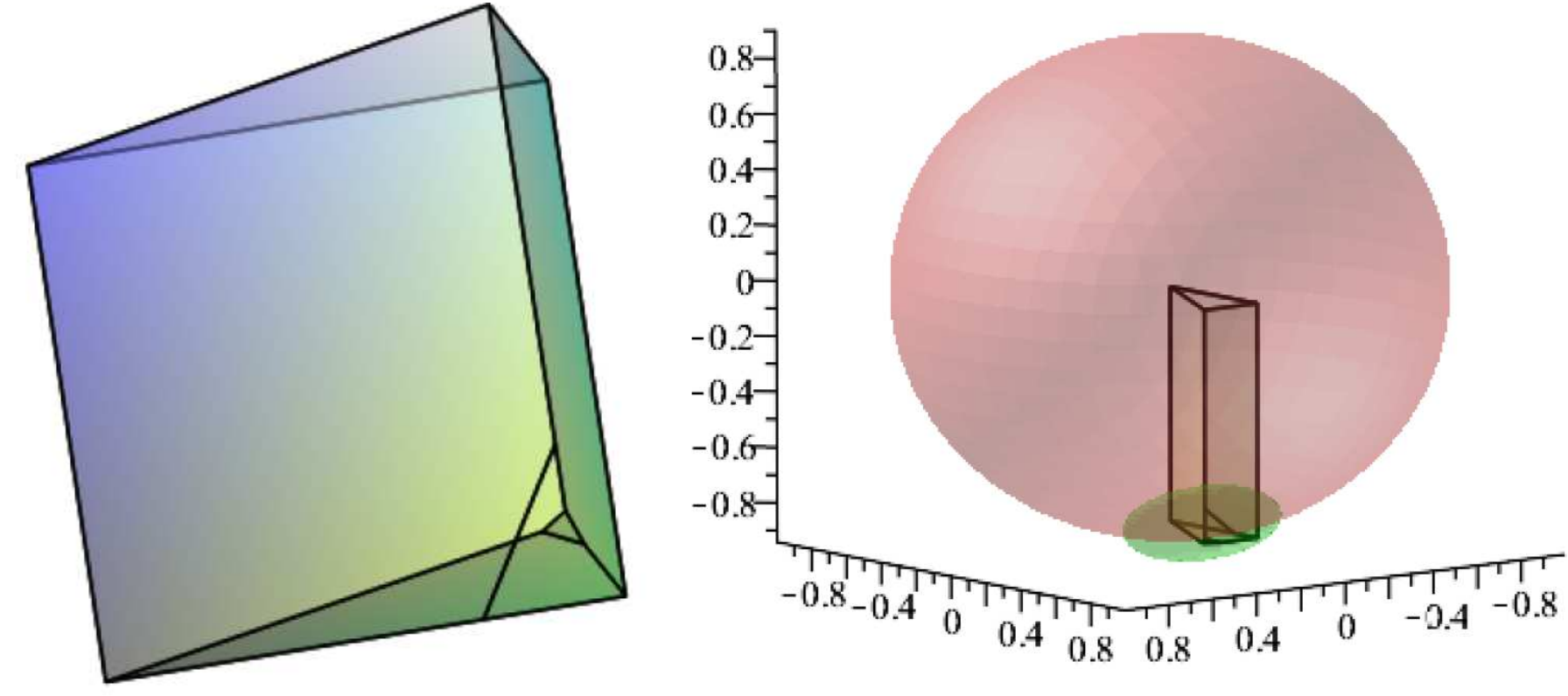} 
\caption{The least dense congruent hyperball covering arrangement 
to parameters $\{7,3,7\}$ 
with density $\approx 1.26829$ and the corresponding Coxeter orthoscheme in Euclidean B-C-K model.}
\label{}
\end{figure}


\begin{thebibliography}{33}
%
\bibitem{FTL} Fejes~T\'oth~L., Regular Figures, \textit{Macmillian (New York)}, 1964.
%
\bibitem{K89} Kellerhals R., On the volume of hyperbolic polyhedra, {\it Math. Ann.}, {\bf 285} (1989) 541--569.
%
\bibitem{LMV18} Lu\v{c}i\'{c} Z., Moln\'{a}r E., Vasiljevi\'{c} N., An algorithm for classification of fundamental polygons for a plane discontinuous group, {\it Geom-Sym Conf. 2015: Discrete Geometry and Symmetry} {\bf 234} (2018) 257--278, DOI: 10.1007/978-3-319-78434-2\_14
%
\bibitem{BM} Maskit B., On Poincar\'{e}'s theorem for fundamental polygons, {\it Adv. in Math.}, {\bf 7} (1971), 219--230.
%
\bibitem{M90} Moln\'{a}r E., Eine Klasse von hyperbolischen Raumgruppen, {\it Beitr\"{a}ge zur Algebra und Geometrie}, {\bf 30} (1990), 79--100.
%
\bibitem{M92} Moln\'{a}r E., Polyhedron complexes with simply transitive group actions and their realizations, {\it Acta Math. Hung.}, {\bf 59} (1-2) (1992), 175--216.
%
\bibitem{MP88} Moln\'{a}r E., Prok I., A polyhedron algorithm for finding space groups, {\it Proceedings of Third Int. Conf. On Engineering Graphics and Descriptive Geometry}, Vienna 1988, {\bf 2}, 37--44.
%
\bibitem{MP94} Moln\'{a}r E., Prok I., Classification of solid transitive simplex tilings in simply connected 3-spaces, Part I, {\it Colloquia Math. Soc. J\'{a}nos Bolyai 63. Intuitive Geometry}, Szeged (Hungary), 1991. North-Holland (1994), 311--362.
%
\bibitem{MPS97} Moln\'{a}r E., Prok I., Szirmai  J., Classification of solid transitive simplex tilings in simply connected 3-spaces, Part II, {\it Periodica Math. Hung.} {\bf 35} (1-2) (1997), 47--94.
%
\bibitem{MPS06} Moln\'{a}r E., Prok I., Szirmai J., Classification of tile-transitive 3-simplex tilings and their realizations in homogeneous spaces, {\it Non-Euclidean Geometries, J\'{a}nos Bolyai Memorial Volume, Editors: A. Pr\'{e}kopa and E. Moln\'{a}r, Mathematics and Its Applications}, {\bf 581}, Springer (2006), 321--363.
%
\bibitem{Mo11} Moln\'{a}r E.,  On $D-$symbols and orbifolds in an algorithmic way, 
{\it Atti Semin. Mat. Fis. Univ. Modena Reggio Emilia}, {\bf 58} (2011), 263--276.
%
\bibitem{MSW11} Moln\'{a}r  E., Szirmai J.,  Weeks J., 3-simplex tilings, splitting orbifolds and manifolds, {\it Symmetry: Culture and Science}, {\bf 22} (3-4) (2011), 435-458.
%
\bibitem{MS18} Moln\'{a}r E., Szirmai J., Top dense hyperbolic ball packings and coverings for complete Coxeter orthoscheme groups, {\it Publications de l'Institut Mathematique}, {\bf 103} (117) (2018), 129--146, DOI: 10.2298/PIM1817129M
%
\bibitem{S97} Stojanovi\'{c} M., Some series of hyperbolic space groups, {\it Annales Univ. Sci. Budapest}, Sect. Math. {\bf 36} (1993), 85--102.
%
\bibitem{S10} Stojanovi\'{c} M., Fundamental simplices with outer vertices for hyperbolic groups, {\it Filomat}, {\bf 24} (1) (2010), 1--19. 
%
\bibitem{S11} Stojanovi\'{c} M., Four series of hyperbolic space groups with simplicial domains, and their supergroups, {\it Krag. J.Math.}, {\bf 35} (2) (2011), 303--315.
%
\bibitem{S13} Stojanovi\'{c} M., Supergroups for six series of hyperbolic simplex groups, {\it Periodica Math. Hung.}, {\bf 67} (1) (2013), 115--131.
%
\bibitem{S14} Stojanovi\'{c} M., Coxeter groups as automorphism groups
of solid transitive 3-simplex tilings, {\it Filomat}, {\bf 28} (3) (2014), 557--577,
DOI: 10.2298/FIL1403557S
%
\bibitem{S17} Stojanovi\'{c} M., Hyperbolic space groups and their supergroups for fundamental simplex tilings, {\it Acta Math. Hung.}, {\bf 153} (2) (2017), 276--288, DOI: 10.1007/s10474-017-0761-z 
%
\bibitem{S19}	Stojanovi\'{c} M., Hyperbolic space groups with truncated simplices as fundamental domains, {\it Filomat}, {\bf 33} (4) (2019), 1107-1116, DOI: 10.2298/FIL1904107S
%
\bibitem{Sz13-4} Szirmai~J., The least dense hyperball covering to the regular prism tilings in the hyperbolic $n$-space,
\textit{Ann. Mat. Pur. Appl.} (2016) {\bf 195}, 235-248, DOI: 10.1007/s10231-014-0460-0.
%
\bibitem{Sz17-1} Szirmai~J., Density upper bound of congruent and non-congruent hyperball packings generated by truncated regular simplex tilings,
\textit{Rendiconti del Circolo Matematico di Palermo Series 2}, {\bf 67} (2018), 307--322, DOI: 10.1007/s12215-017-0316-8.
%
\bibitem{Sz17-2} Szirmai~J., Decomposition method related to saturated hyperball packings,
\textit{Ars Mathematica Contemporanea}, {\bf{16}} (2019), 349--358, DOI: 10.26493/1855-3974.1485.0b1.
%
\bibitem{Sz18} Szirmai~J., Hyperball packings in hyperbolic 3-space, {\it Ma\-te\-ma\-ti\-\v{c}ki Vesnik}, {\bf 70} (3) (2018), 211--221.
%
\bibitem{Sz19-4} Szirmai~J.,
Hyperball packings related to cube and octahedron tilings in hyperbolic space,
\emph{Contributions to Discrete Mathematics} (to appear) (2019), arXiv:1803.04948.
%
\bibitem{Sz19} Szirmai~J., Upper bound of density for packing of congruent hyperballs in hyperbolic $3-$space,
\textit{Submitted manuscript}, (2019).
%
\bibitem{V}
Vinberg~E.~B. (Ed.),
\textit{Geometry II. Spaces of Constant Curvature}
Spriger Verlag Berlin-Heidelberg-New York-London-Paris-Tokyo-Hong Kong-Barcelona-Budapest, 1993.
%
\bibitem{Z} Zhuk~I.~K., Fundamental tetrahedra in Euclidean and Lobachevsky spaces, {\it Soviet Math. Dokl.}, {\bf 270} (3) (1983), 540--543.

\end{thebibliography}
\end{document}